\documentclass{amsart}

\usepackage{amsmath}    
\usepackage{amssymb}
\usepackage{amsfonts}
\usepackage{amsthm}
\usepackage[all,arc]{xy}
\usepackage{enumerate}
\usepackage{mathrsfs}
\usepackage{hyperref}
\hypersetup{
	pdftitle={},
	pdfsubject={},
	pdfauthor={Saba Aliyari},
	pdfkeywords={}
}
\usepackage{bm}
\usepackage{pgfplots}
\usepackage{enumitem}
\usepackage{relsize}
\usepackage{version}
\includeversion{fullversion}
\excludeversion{shortversion}
\newtheorem{thm}{Theorem}[section]

\newtheorem{prop}[thm]{Proposition}
\newtheorem{lem}[thm]{Lemma}

\newtheorem{claim}[thm]{Claim}

\theoremstyle{definition}
\newtheorem{defn}[thm]{Definition}

\newtheorem{exmp}[thm]{Example}

\newtheorem{notn}[thm]{Notation}

\newtheorem{assum}[thm]{Assumption}

\theoremstyle{remark}
\newtheorem{rem}[thm]{Remark}

\newcommand{\eqdef}{=\joinrel=}
\DeclareMathOperator{\mathoptrop}{trop}
\DeclareMathOperator{\mathopin}{in}
\DeclareMathOperator{\mathopsupp}{supp}
\numberwithin{equation}{section}

\title{The Difference Kapranov Theorem}

\author{Saba Aliyari}

\pgfplotsset{compat=1.18}
\begin{document}
	\begin{abstract}
		This paper creates a link between \textit{Tropical Geometry} and \textit{Difference Algebra}. The main result is a difference version of \textit{Kapranov's Theorem}. In this theorem, we extend Kapranov's Theorem to the case of a Laurent difference polynomial with coefficients from a multiplicative valued difference field, where the residue field is an algebraically closed field with a generic automorphism (ACFA). A result of this paper that plays a critical role in the proof of the Difference Kapranov Theorem is a difference version of \textit{Newton's Lemma}.    
		
		\noindent \textbf{Keywords.}
		\noindent \textbf{Multiplicative valued difference fields; ACFA; Difference Newton Lemma; Kapranov's Theorem } 
	\end{abstract}

	\maketitle
	\section{Introduction}
	In this paper, we establish a \textit{difference} version of \textit{Tropical Geometry} that connects \textit{Tropical Geometry} to \textit{Difference Algebra}. To accommodate the readers from both areas, we start with a brief introduction to both.
	\par \textit{Tropical Geometry} is a field that studies polynomials and their geometric properties in the \textit{tropical semiring} $\left(\mathbb{R} \cup \lbrace \infty\rbrace ,\bigoplus,\bigodot\right)$, where $\bigoplus$ is the minimum and $\bigodot$ is the classical addition. Tropical arithmetic and a procedure called \textit{tropicalization} enable us to find tropical analogues of classical mathematical objects, for instance polynomials. Polynomials may have complicated graphs, but tropicalization turns them into piecewise linear graphs, which are much easier objects to study.
	Therefore, Tropical Geometry has been a useful tool. For an overview of this field, see, for example, \cite{Katz}.
	\par One of the important theorems that builds this connection is \textit{Kapranov's theorem}. It forges a connection between algebraic hypersurfaces and \textit{tropical hypersurfaces} in $\mathbb{R}^n$. This connection can be generalized to arbitrary varieties in the \textit{Fundamental Theorem of Tropical Algebraic Geometry} which is a central theorem in this field. In fact, Kapranov's Theorem is considered as a critical step in the proof of the Fundamental Theorem. 
	\begin{shortversion}
		\textcolor{red}{As suggested by Immi, I omitted the historical part of Tropical Geometry, although I think it is better to treat the readers of both areas equally. }
		\par The name of this field was chosen in honor of Imre Simon, who first introduced this theory. He is a Hungarian-born computer scientist from the tropical country of Brazil.
		Efforts to consolidate the definitions of the theory began in the late 1990s. Imre Simon, Grigory Mikhalkin, and Bernd Sturmfels, along with many other mathematicians, made significant contributions in this area. 
	\end{shortversion}
	For further information on this topic, helpful references include \cite{Diane} and \cite{Mikhalkin}.
	\par \textit{Difference Algebra} is an area in mathematics which studies \textit{difference field}s and \textit{difference polynomial}s. A difference field is a field $K$ together with a field automorphism $\sigma$. A difference polynomial in $n$ variables $x_1, \dots, x_n$ with coefficients from $K$, is a polynomial in infinitely many formal variables $\sigma^j(x_i)$, for $i\in \lbrace 1, \dots,n\rbrace$ and $j\in \mathbb{N}$ with $\sigma^j$ being the $j$-th iteration of $\sigma$. In this case, we say that $f$ is an element of the \textit{ring of difference polynomials} and use the notation $f\in K_{\sigma}[x_1, \dots,x_n]$. 
	\par \textit{Difference Algebra} is considered as an analogue to \textit{Differential Algebra}, but in this area difference equations are studied rather than differential equations. Difference Algebra was first introduced as an independent field of study by Joseph Ritt and Richard Cohn. Later, Hrushovski applied it in proving \textit{The Manin-Mumford Conjecture} in 2001 \cite{Manin}. Another interesting application appears in the connection between \textit{algebraic dynamical systems} and difference fields. In fact the knowledge about difference fields leads to a better understanding of algebraic dynamical systems, \cite{Scanlon} provides valuable insights on this topic. 
	\par Difference polynomials are complex objects to study, and determining their roots, if possible, is not an easy task. For example, even a very simple difference polynomial, such as $x-x^{\sigma}$, defines a variety that is an infinite field that contains $\mathbb{Q}$. This aspect served as a motivation for our project, in which we utilized tropical tools to gain a deeper understanding of difference polynomials and their roots.
	\par A \textit{valued difference field} $K$, is a valued field together with an automorphism that takes the valuation ring to the valuation ring. For an example, see Example \ref{Hahn}. The valuation and the automorphism of $K$ can interact in different ways, which are explained in Remark \ref{cases}. In this work, we assume the value group $\Gamma$ to be a divisible subgroup of $\mathbb{R}$. We work with the case of \textit{multiplicative valued difference fields} in which  we have 
	\begin{equation*}
		\forall x\in K^{\times}:\,\,\,\,\,v(\sigma(x))= \rho\cdot v(x),
	\end{equation*}
	and the value group is a $\mathbb{Z}[\rho]$-module where $\rho$ is a fixed positive real number. 
	
	We call $\rho$ the \textit{scaling exponent} of $\sigma$. In addition, throughout this paper we assume that $\Gamma$ is a $\mathbb{Q}(\rho)$-module for $\rho$ being transcendental over $\mathbb{Q}$.
	\par For further exploration of Difference Algebra, we recommend referring to \cite{Wibmer}, \cite{Cohn} and \cite{Levin}.
	\par The goal of this paper is to extend Kapranov's Theorem to the case where $f$ is a Laurent difference polynomial with coefficients from a multiplicative valued difference field.  Perhaps, in the future, the present work could serve as the base case for establishing a difference version of the Fundamental Theorem. We should also emphasize that the extension of Kapranov's Theorem to other cases of valued difference fields has not been studied yet.
	\par For a Laurent polynomial $f$ in $K[x_1^{\pm1},\dots,x_n^{\pm1}]$, Kapranov's Theorem builds a bridge between its associated classical hypersurface $V(f)$ and the associated tropical hypersurface $\mathoptrop(V(f))$. More precisely, it states
	\begin{thm}(Kapranov's Theorem)
		If $(K,v)$ is an algebraically closed valued field, with nontrivial valuation, and if $f \in K[x_1^{\pm1},\dots,x_n^{\pm1}] $, then the following sets coincide:
		\begin{enumerate}
			\item $\mathoptrop(V(f))$ which is a subset of $\mathbb{R}^n$;
			\item the set of all points in $\mathbb{R}^n$ for which the initial form of $f$ is not a monomial;
			\item the topological closure of $\lbrace \left(v(y_1),\dots v(y_n)\right) |\,\, \left(y_1,\dots, y_n\right) \in V(f) \rbrace$ in $\mathbb{R}^n$. 
		\end{enumerate}
	\end{thm}
	The proof of this theorem can be found, for instance, in \cite{Kapranov}, and in \cite{Diane}.
	\par The equality of the sets appearing in $(1)$ and $(3)$ provides information about the roots of $f$. In this paper, we obtain a corresponding result for the case that $f\in K_{\sigma}[x_1^{\pm1}, \dots ,x_n^{\pm1}]$, where $K$ is a multiplicative valued difference field. 
	\par This is the main result of this paper, which is the Difference Kapranov Theorem (Theorem \ref{diffkap}), the statement of which is provided below.
	The necessary definitions required to comprehend the statement of this theorem are given on the following page.
	For the definition of $\mathopin_w(f)$, see Definition \ref{diffinitial}. 
	In addition, in this theorem, we assume $\rho$ to be transcendental over $\mathbb{Q}$. This appears in the tropicalization of $f$.
	In fact, if $\rho$ is algebraic, for an exponent $u(\sigma)$ appearing in $f$, $u(\rho)$ that appears in the corresponding tropical monomial might annihilate. This affects the number of tropical roots, and as a result 
	the accuracy of the Difference Kapranov Theorem. 
	\begin{thm}{( The Difference Kapranov Theorem)}
		
		Let $K$ be a multiplicative valued difference field of characteristic zero that is spherically complete.
		Assume its difference residue field is an ACFA of characteristic zero and the scaling exponent $\rho$ is transcendental.
		Let the  difference value group $\Gamma$ of $K$ be a subgroup of $\mathbb{R}$ which is a $\mathbb{Q}(\rho)$-module.\\
		Suppose $f \in K_{\sigma}\left[x_1^{\pm1},\dots,x_n^{\pm1}\right]$ is a Laurent difference polynomial. The following sets coincide:
		\begin{enumerate}
			\item $\mathoptrop(V(f)) \subseteq \mathbb{R}^n$ which is the difference tropical hypersurface associated to $f$;
			\item the set of all the points $w \in \mathbb{R}^n$ for which the initial form $\mathopin_w(f)$ is not a monomial;
			\item the topological closure of the set
			\begin{equation*}
				A=\left\lbrace (v(y_1),\dots,v(y_n)):(y_1,\dots,y_n)\in V(f) \right\rbrace,   
			\end{equation*}
			in $\mathbb{R}^n$.
		\end{enumerate} 
		
	\end{thm}
	\par In \cite{Diane}, in order to prove that the set in $(1)$ is included in the set in $(3)$, induction is used. In the base case, the assumption that $K$ is an algebraically closed field plays an important role. Using this assumption, it is possible to decompose $f$ into linear factors, which finally results in the existence of a root with the desired conditions. In contrast, in our context, concerning Laurent \textit{difference} polynomials, we do not have such a decomposition. As a solution, we proved the \textit{Difference Newton Lemma}(\ref{DiffNewton}), which is another result of this paper, to guarantee the existence of a root of a difference polynomial in one variable.
	\par As a natural next step, one may reflect on the Difference version of the Fundamental Theorem of Tropical Algebraic Geometry. Based on our current knowledge, we are uncertain whether the difference version of the fundamental theorem exists. At first glance, it seems difficult to prove, and if possible, it may need a significant amount of effort. 
	\par \textbf{Section 2} consists of two subsections. The first is devoted to \textit{Difference Tropical Geometry}.
	We begin this subsection with some preliminaries in Difference Algebra. Then we explain how they can be connected to Tropical Geometry, and we define difference tropical objects.
	\par One of the main assumptions in this paper is that the residue field is an \textit{ACFA}. This concept comes from \textit{Model Theory} and model theorists say that a field is a model of ACFA, which is an abbreviation for \textit{Algebraically Closed Field with a generic Automorphism}. A field $k$ is a model of ACFA if $\sigma$ is an automorphism of $k$, and it satisfies the conditions of Theorem \ref{axiomatization} which are the axioms presented in $3.1$ of \cite{Zoe}.
	In this paper, we do not use the model theoretic approach. Therefore, we simply say that a field is an ACFA. This notion can be interpreted as a difference version of algebraically closed fields.
	For a precise definition, see Definition \ref{ACFA}. In this subsection, we prove Theorem \ref{alpharoot} which states
	\begin{thm}
		Let $k$ be an ACFA. Suppose $f$ is in $k_{\sigma}[x_1^{\pm 1},\dots,x_n^{\pm 1}]$, and is not a monomial. Then $f$ has a root in $(k^*)^n$.
	\end{thm}
	The significance of this theorem becomes evident in the proof of the Difference Kapranov Theorem.
	The key step in the proof of Theorem \ref{alpharoot} is the one variable case, which is Lemma \ref{ACFAroot}. 
	\par An important assumption on $K$ is that it is \textit{spherically complete}. This means that if we take any chain of balls in $K$, their intersection is nonempty.
	\par Given a Laurent difference polynomial, we can define its tropicalization and its initial form (Definition \ref{difftrop} and Definition \ref{diffinitial}). These are the concepts which are introduced at the end of this subsection. Moreover, a \textit{difference tropical hypersurface} is defined to be the set of all tropical roots of a Laurent difference polynomial.
	\par Since \textit{Polyhedral Geometry} is the main tool we will be using in \textbf{Section 3}, the second subsection of \textbf{Section 2} is dedicated to introducing some preliminaries of this topic. Most readers familiar with Tropical Geometry are already acquainted with the concepts covered in this subsection. 
	The only new object introduced here is a $(\Gamma, \mathbb{Q}(\rho))$-\textit{polyhedral complex}. It is defined in Definition \ref{Gamma}. This is in fact the difference analogue of a $\Gamma$-rational polyhedral complex.
	\par In \textbf{Section 3}, we prove the first result of this paper, namely Proposition \ref{combi}. This proposition is about the combinatorics of a difference tropical hypersurface. It states
	\begin{prop}
		If $f=\underset{u(\sigma)}{\sum}c_{u(\sigma)}x^{u(\sigma)}$ is a Laurent difference polynomial, then its associated difference tropical hypersurface, $\mathoptrop(V(f))$, is the support of a pure $\left(\Gamma,\mathbb{Q}(\rho)\right)-$polyhedral complex of dimension $(n-1)$.\\
		More precisely, it is the $(n-1)$-skeleton of the polyhedral complex dual to the regular subdivision $\Sigma_{val}$ of $P$, where $P$ is defined as in Definition \ref{Newpolytope}, and $val$ is the weight vector given by $v\left( c_{u(\sigma)}\right) $ for $c_{u(\sigma)}\neq 0$.
	\end{prop}
	\par In \textbf{Section 4}, we prove the Difference Newton Lemma. As explained before, this theorem plays a critical role in our proof of the Difference Kapranove Theorem. This theorem states the following.
	\begin{thm}(Difference Newton Lemma)
		Let $K$ be a multiplicative valued difference field of characteristic zero that is spherically complete.
		Assume $\mathbf{k}$, the difference residue field of $K$, is an ACFA and of characteristic zero. Suppose that the difference value group $\Gamma$ is a subgroup of $\mathbb{R}$ that is a $\mathbb{Q}(\rho)$-module, where $\rho$, the scaling exponent of $\sigma$, is transcendental.
		Given $f\in K_{\sigma}[x]$ is not constant and suppose $b\in K$ such that $f(b)\neq 0$.\\
		For some $\varepsilon$, there exists a root $a \in K $ of $f$ such that $v(a-b)=\varepsilon$.
	\end{thm}
	For the precise definition of $\varepsilon$ see Theorem \ref{DiffNewton}, Definition \ref{f_J} and Notation \ref{multindex}.
	\par To prove this theorem, we follow two main steps. Firstly, in Lemma \ref{l1}, we prove that we can have a better estimation of a possible root of $f\in K_{\sigma}[x]$ around a nonroot $b\in K$. Finally, in Proposition \ref{l2}, we use the assumption of spherical completeness to find a root for $f$. In a spherically complete field, any pseudocauchy sequence has a pseudolimit, see \cite{Kaplansky}. In Proposition \ref{l2}, we construct a pseudocauchy sequence whose pseudolimite is a root of $f$. 
	\par The main result of this paper is proved in \textbf{Section 5}, namely the Difference Kapranov Theorem. The main tool to prove this theorem is Proposition \ref{main}, whose statement is given below.
	\begin{prop}\label{improp}
		Let $f \in K_{\sigma}[x_1^{\pm1},\dots,x_n^{\pm1}]$ be a Laurent difference polynomial, and $\underline{w}=(w_1, \dots, w_n) \in \Gamma^n$ such that $\mathopin_{\underline{w}}(f)$ is not a monomial. Suppose $\bar{\alpha}$ is a root of $\mathopin_{\underline{w}}(f)$ in $(\mathbf{k^*})^n$ where $\mathbf{k}$ is the difference residue field of $K$. Then there exists an element $y$ in $(K^*)^n$ which is a root of $f$, and satisfies the following conditions:
		\begin{itemize}
			\item $v(y)=\underline{w}$,
			\item $\forall i,\,\,\,1\leq i \leq n:\,\, \overline{t^{-w_i}\cdot y_i}=\bar{\alpha}_i$.
		\end{itemize}
	\end{prop}
	To prove this proposition, we start by proving the same statement for a Laurent difference polynomial in one variable. This is done in Lemma \ref{1case} and Lemma \ref{n=1}. In this case, Difference Newton Lemma plays a highlighted role. 
	Then with the same strategy as in the proof from \cite{Diane}, we prove Proposition \ref{main} which is Proposition \ref{improp} above.
	In the proof of the Difference Kapranov Theorem, the tricky part is to prove that the set appearing in $(1)$ is included in the set appearing in $(3)$ and this is done by Proposition \ref{main}.
	\par \textbf{Acknowledgements:} This paper is derived from my Ph.D. thesis supervised by Immanuel Halupczok at Heinrich-Heine-Universität Düsseldorf. I would like to thank him for all his guidance and the fruitful discussions we had. I am also grateful for his feedback and comments on this work.
	I wish to express my thanks to David Bradley-Williams for his valuable comments and suggestions. 
	I would also like to thank Zoé Chatzidakis for her contribution to Lemma \ref{ACFAroot}.
	\begin{fullversion}
		\section{Preliminaries}
		\subsection{Difference Tropical Geometry }
		\subsubsection{Difference Algebra}
		In this part, we provide the difference algebra needed in this paper. The main references here are \cite{Wibmer}, \cite{Pal}, and \cite{Zoe}. \cite{Wibmer} is a well-written lecture note on this topic and is used as a reference in this section. Any material presented from this source can also be found in \cite{Cohn}. In addition, we establish our general assumptions and present all necessary definitions.       
		
		\begin{defn}
			A \textit{difference field} is a field $K$ together with an automorphism $\sigma:K \longrightarrow K$. It is denoted by $\left(K, \sigma \right)$.
		\end{defn}
		\begin{defn}
			A \textit{valued difference field} is a valued field $K$ together with an automorphism $\sigma$ satisfying $\sigma \left(\mathcal{O}_K\right)=\mathcal{O}_K$, where $\mathcal{O}_K$ is the valuation ring.
		\end{defn}
		\begin{rem} The automorphism of a valued difference field induces two important automorphisms as follows:
			\begin{itemize}
				\item Let $\Gamma$ be the value group of $K$. Then $\sigma$ induces an automorphism $\sigma_{\Gamma}$ on $\Gamma$ as below:
				\begin{align*}
					\sigma_{\Gamma}:& \Gamma \longrightarrow \Gamma,\\
					&v(a) \longmapsto v \left( \sigma (a) \right),
				\end{align*}
				where $a$ is an element of $K$.\\
				In this case, $\Gamma$ is called a \textit{difference value group}.
				\par The property $\sigma \left(\mathcal{O}_K\right)=\mathcal{O}_K$ guarantees that $\sigma_{\Gamma}$ is well-defined.
				\begin{shortversion}
					To see this, assume $a$ and $a^{\prime}$ are two elements of $K$ such that $v(a)=\gamma=v(a^{\prime})$. We want to prove that $v\left(\sigma(a)\right)=v\left(\sigma(a^{\prime})\right)$. $v(a)=v(a^{\prime})$ means $v\left(\frac{a}{a^{\prime}}\right)=0$, or equivalently $\frac{a}{a^{\prime}}\in \mathcal{O}_K$. From $\sigma \left(\mathcal{O}_K\right)=\mathcal{O}_K$, we have $\sigma\left(\frac{a}{a^{\prime}}\right) \in \mathcal{O}_K$. This means that $v\left(\sigma \left(\frac{a}{a^{\prime}}\right)\right)\geq 0$ which means $v\left(\sigma(a)\right)\geq v\left(\sigma(a^{\prime})\right) $. Similarly, we can show $v\left(\sigma(a^{\prime})\right)\geq v\left(\sigma(a)\right)$. Hence, $v\left(\sigma(a^{\prime})\right)= v\left(\sigma(a)\right)$ and $\sigma_{\Gamma}$ is well-defined.\\
					$\sigma_{\Gamma}$ is also order preserving. Assume $\gamma > \gamma^{\prime}$ with $\gamma=v(a)$ and $\gamma^{\prime}=v(a^{\prime})$ for $a, a^{\prime} \in K$. This means $v\left(\frac{a}{a^{\prime}}\right)>0$ or equivalently $\frac{a}{a^{\prime}} \in \mathcal{M}$, where by $\mathcal{M}$ we mean the maximal ideal of $\mathcal{O}_K$. As $\mathcal{M}$ contains all noninvertible elements, $\frac{a}{a^{\prime}} $ is not invertible. Therefore, $\sigma\left(\frac{a}{a^{\prime}} \right)$ is not invertible, which means $\sigma\left(\frac{a}{a^{\prime}} \right) \in \mathcal{M}$ or equivalently we have
					$v\left(\sigma\left(\frac{a}{a^{\prime}} \right)\right)>0$. This gives $v\left(\sigma(a)\right)>v\left(\sigma(a^{\prime})\right)$, in fact $\sigma_{\Gamma}(\gamma)>\sigma_{\Gamma}(\gamma^{\prime})$.
					
				\end{shortversion}
				
				\item 
				$\sigma$ also induces an automorphism $\bar{\sigma}$ on $\mathbf{k}$, the residue field of $K$, as follows:
				\begin{align*}
					\bar{\sigma}:&\mathbf{k} \longrightarrow \mathbf{k},\\
					& \bar{a} \longmapsto \overline{\sigma(a)}.
				\end{align*}
				
				Then $\left( \mathbf{k}, \bar{\sigma}\right)$ is called the \textit{difference residue field} of $K$.
				Again, from the property $\sigma \left(\mathcal{O}_K\right)=\mathcal{O}_K$, $\bar{\sigma}$ is well defined.
				\begin{shortversion}
					To see this, assume that for two elements $\bar{a}$ and $\bar{b}$ of $\mathbf{k}$, we have $\bar{a}=\bar{b}$. This means $a-b \in \mathcal{M}$. In other words $v\left(a-b\right)>0$. As we discussed for $\sigma_{\Gamma}$, we have $v\left(\sigma(a-b)\right)>0$. Thus, $\sigma(a)-\sigma(b)\in \mathcal{M}$. Hence, $\overline{\sigma(a)}=\overline{\sigma(b)}$.\\
				\end{shortversion}
			\end{itemize} 
		\end{rem}
		
		\begin{notn}
			We fix the above notations for the rest of this work. Explicitly, from now on, for a valued field $K$, we denote the value group of $K$ by $\Gamma$, its residue field by $\mathbf{k}$, and the valuation ring by $\mathcal{O}_K$. 
		\end{notn} 
		As a clarification, we look at an example of a Hahn-field. In general, if a field $\mathbf{k}$ and an ordered abelian group $\Gamma$ are given, then one can define a valued field whose residue field is $\mathbf{k}$ and whose value group is $\Gamma$. This field is called a \textit{Hahn-field} and it is defined as follows:
		\begin{equation*}
			K=\mathbf{k}\left(\left(t^{\Gamma}\right)\right):=\lbrace f=\underset{\gamma \in \Gamma}{\textstyle\sum}a_{\gamma}t^{\gamma}\,\,\,|\,\,\, a_{\gamma}\in \mathbf{k},\,\,\ \text{and $\mathopsupp(f)$ is well-ordered}\rbrace,
		\end{equation*}
		where $\mathopsupp(f)=\lbrace \gamma \in \Gamma\,\,\,|\,\,\, a_{\gamma} \neq 0 \rbrace $.
		$K$ is a field with a natural addition and multiplication, and it is a valued field with the following valuation:
		\begin{align*}
			v: K &\longrightarrow \Gamma,\\
			\underset{\gamma \in \Gamma}{\textstyle\sum}a_{\gamma}t^{\gamma} &\longmapsto \min \left\lbrace \gamma\,\,|\,\, a_{\gamma} \neq 0 \right\rbrace.
		\end{align*}
		
		Suppose $\mathbf{k}$ is a field and $\bar{\sigma}$ is an automorphism on $\mathbf{k}$. Also assume that an ordered abelian group $\Gamma$ with an order preserving automorphism $\sigma_{\Gamma}$ on it are given. Using these automorphisms, an automorphism on the corresponding Hahn-field is defined. This automorphism $\sigma$ of $K$ is defined as follows: 
		\begin{align*}
			\sigma:K &\longrightarrow K,\\
			\underset{\gamma \in \Gamma} {\textstyle\sum}a_{\gamma}t^{\gamma} &\longmapsto \underset{\gamma \in \Gamma}{\textstyle\sum}\bar{\sigma}(a_{\gamma})t^{\sigma_{\Gamma}(\gamma)}. 
		\end{align*}
		The following is a concrete example of a Hahn-field.
		\begin{exmp}\label{Hahn}
			Let $\mathbf{k}$ be the field of complex numbers, $\mathbb{C}$, and $\Gamma$ be $\mathbb{R}$, regarded as an ordered abelian group. As it is defined above, the corresponding Hahn-field is 
			\begin{equation*}
				K=\mathbb{C}\left(\left(t^{\mathbb{R}}\right)\right)=\lbrace f=\underset{\gamma \in \mathbb{R}}{\textstyle\sum}a_{\gamma}t^{\gamma}\,\,\,|\,\,\, a_{\gamma}\in\mathbb{C},\,\,\ \text{and $\mathopsupp(f)$ is well-ordered}\rbrace.
			\end{equation*}
			Assume $\mathbb{C}$ is equipped with the identity automorphism; in this case, $\left(\mathbb{C},id \right)$ is called a constant difference field. If we consider $\mathbb{R}$ as an ordered abelian group, then any automorphism $\sigma_{\Gamma}$ on $\mathbb{R}$ is of the following form:
			\begin{equation*}
				x \mapsto \sigma_{\Gamma}(x)=\rho \cdot x \,\,\, \text{for some fixed}\,\,\, \rho > 0.
			\end{equation*}
			Using these two automorphisms, we can define an automorphism $\sigma$ on the Hahn-field $K$ as follows:
			\begin{align*}
				\sigma:K &\longrightarrow K,\\
				\underset{\gamma \in \mathbb{R}}{\textstyle\sum}a_{\gamma}t^{\gamma} &\longmapsto \underset{\gamma \in \mathbb{R}}{\textstyle\sum}a_{\gamma}t^{\rho\cdot\gamma}. 
			\end{align*}
			Then $\left(K,v, \sigma \right)$ is a valued difference field.
		\end{exmp} 
		\begin{rem}\label{cases}
			The valuation and the automorphism of a valued difference field can interact in different ways. This interaction results in different cases of valued difference fields. Below, we present three interesting cases:
			\begin{itemize}
				\item The \textit{isometric} case in which the following condition holds:
				\begin{equation*}
					\forall x\in K^{\times}:\,\,\,\,\, v(\sigma(x))=v(x).
				\end{equation*}
				The eager reader can see \cite{BMS}  to learn more about this case.
				\item The \textit{contractive} case in which $\forall x\in K^{\times}$ with $v(x)>0$, the following condition is satisfied:
				\begin{equation}\label{contra}
					\forall n\in \mathbb{N}:\,\,\,\,\,v(\sigma(x))> nv(x). 
				\end{equation}
				In this case, we have $rank(\Gamma)>1$. This means that it can not be embedded in $\mathbb{R}$, see \cite{Prestel}.  
				\begin{shortversion}
					In this case, if the value group is $\mathbb{R}$, then condition \eqref{contra} implies $\forall x\in K^{\times}$ with $v(x)>0$, we have $v(\sigma(x))=\infty$, or equivalently $\sigma(x)=0$. This is not an automorphism. To fix this, we add a nonstandard element to $\mathbb{R}$, which is greater than all real numbers. Therefore, $v(\sigma(x))$ can be larger than all real numbers and also not $\infty$.\\
				\end{shortversion}
				This issue is far from the main subject of this paper, but the eager reader can look at \cite{Azgin} to learn more about contractive valued difference fields. 
				
				\item The \textit{multiplicative} case in which the difference value group is a $\mathbb{Z}[\rho]$-module with $\rho$ being a positive real number, and we have
				\begin{equation*}
					\forall x\in K^{\times}:\,\,\,\,\,v(\sigma(x))= \rho\cdot v(x).
				\end{equation*}
				This case is well studied in \cite{Pal}.\\
				Note that the isometric case is a special case of the multiplicative case.
			\end{itemize}
			\par In this paper, we assume that the difference value group $\Gamma$ of the field is a divisible subgroup of $\mathbb{R}$.
			Therefore, any automorphism $\hat{\sigma}$ of ordered abelian groups on $\Gamma$ is of the following form:
			\begin{equation*}
				\forall x\in \Gamma: \,\,\, \hat{\sigma}(x)= \rho\cdot x\,\, \text{where}\,\, \rho\,\, \text{is a fixed positive real number}.
			\end{equation*}
			So, if we consider the induced automorphism on $\Gamma$,
			for a fixed $\rho$, we have
			\begin{equation*}
				\forall x\in K^{\times}:\,\,\,\,\, v\left(\sigma(x)\right)=\sigma_{\Gamma}\left(v(x)\right)=\rho\cdot v(x).
			\end{equation*}
			This means that in this paper, we work with multiplicative valued difference fields. See Assumption \ref{assump}. We call $\rho$ the \textit{scaling exponent} of the automorphism $\sigma$.
		\end{rem}
		\begin{rem}
			In Assumption \ref{assump}, we will make further assumptions about $\Gamma$. Specifically, we assume that $\Gamma$ is a $\mathbb{Q}(\rho)$-module.
		\end{rem}
		\begin{defn}
			Let $K$ be a difference field. The \textit{difference polynomial ring} over $K$, in difference variables $x=(x_1,\dots,x_n)$, is denoted by $K_{\sigma}\left[x\right]$. It is the polynomial ring over $K$ in  formal variables $\sigma^{i}(x_j)$ for $i\in \mathbb{N}$ and $j\in \lbrace1,\dots,n\rbrace$ where $\sigma^0(x_j):=x_j$. In other words, we have
			\begin{equation*}
				K_{\sigma}\left[x\right]=K\left[\sigma^{i}(x_j)\,\, |\,\, i\in \mathbb{N},\,\,\,\, j\in \lbrace1,\dots ,n\rbrace\right].
			\end{equation*}
			Any element of a difference polynomial ring is called a \textit{difference polynomial}.
			\par Note that here by $x$ we mean $(x_1, \dots ,x_n)$, even though in some parts of this paper, $x$ may refer to a single variable, which will be clear from the context.
			\par Similarly, we can define the \textit{ring of Laurent difference polynomials} in $n$ variables over $K$, which is denoted by $K_{\sigma}[x_1^{\pm},\dots,x_n^{\pm}]$. Any element of this ring is called a \textit{Laurent difference polynomial}. 
		\end{defn}
		\begin{notn}\label{sigma}
			We commonly use the notation $\sigma(x)=x^{\sigma}$. Using this notation, the difference monomial $x^{a_0}\sigma(x)^{a_1} \dots \sigma^m(x)^{a_m}$ in one variable $x$ can be written as $x^{a_0+a_1\sigma+ \dots +a_m\sigma^{m}}$.
			\par If $\mathbb{Z}[\sigma]$ denotes the set
			\begin{equation*}
				\lbrace \sum_{i=0}^{m}a_i\sigma^{i}|\,\, \forall i,\,\, a_i\in \mathbb{Z}\,\, \text{and}\,\, \sigma^{i}\,\, \text{is the i-th iteration of}\,\, \sigma \rbrace, 
			\end{equation*}
			by using this notation, all the exponents appearing in a Laurent difference polynomial in $n$ variables are elements of $\left(\mathbb{Z}[\sigma]\right)^n$, which are called \textit{$\sigma$-powers}.  Note that in this case, we use the notation $x^{u(\sigma)}:=x_1^{u_1(\sigma)}\cdots x_n^{u_n(\sigma)}$. This means, if $f$ is a Laurent difference polynomial in variables $x_1, \dots x_n$, it can be written as   
			\begin{equation*}
				f(x)=\underset{u(\sigma)\in \Lambda}{\sum}c_{u(\sigma)}x^{u(\sigma)},
			\end{equation*} 
			where $\Lambda$ is a finite subset of $\left(\mathbb{Z}[\sigma]\right)^n$.
			\par Take the subset $\mathbb{N}[\sigma]$ of $\mathbb{Z}[\sigma]$ in which for each $i$, $a_i$ is an element of $\mathbb{N}$, then the $\sigma$-powers appearing in a difference polynomial in $n$ variables are elements of $\left(\mathbb{N}[\sigma]\right)^n$, meaning that any difference polynomial $f$, in variables $x_1,\dots,x_n$ can be written as follows:
			\begin{equation*}
				f(x)=\underset{u(\sigma)\in \Lambda}{\sum}c_{u(\sigma)}x^{u(\sigma)},
			\end{equation*}
			with $\Lambda$ being a finite subset of $\left(\mathbb{N}[\sigma]\right)^n$.
		\end{notn}
		\begin{shortversion}
			The following example clarifies Notation \ref{sigma}.
			\begin{exmp}
				Consider 
				\begin{equation*}
					f(x_1,x_2)=x_1^{2}x_2+\sigma(x_1)^4+x_1\sigma^3(x_2)x_2.
				\end{equation*}
				If we use Notation \ref{sigma}, it can be rewritten as:
				\begin{equation*}
					f(x_1,x_2)=x_1^{2}x_2+x_1^{4\sigma}+x_1x_2^{\sigma^3+1}.
				\end{equation*}
				This is a difference polynomial in $\mathbb{C}_{\sigma}[x]$, with $x=(x_1,x_2)$.    
			\end{exmp}
		\end{shortversion}     
		In the case of a difference polynomial in one variable, sometimes it is more convenient to use the following notations (namely Remark \ref{sigmaJ}) rather than the one defined in Notation \ref{sigma}. 
		\begin{notn}\label{multindex}
			$J$ is called a \textit{multi-index}, if it is an element of $\mathbb{N}^{n+1}$. For $J=(j_0,j_1, \dots, j_n)$, its \textit{length} which is denoted by $|J|$ is defined as follows:
			\begin{equation*}
				|J|=j_0+j_1+\dots+j_n.
			\end{equation*}
			For a positive real number $\rho$, the \textit{$\rho$-length} of $J$ which is denoted by $|J|_{\rho}$, is defined as follows:
			\begin{equation*}
				|J|_{\rho}=\rho^0\cdot j_0+\rho^1\cdot j_1+\dots+\rho^n\cdot j_n.  
			\end{equation*}
			Throughout the rest of this paper, when we write $|J|_{\rho}$, $\rho$ refers to the scaling exponent of $\sigma$.
		\end{notn}
		\begin{notn}\label{sigmabold}
			For an automorphism $\sigma$, and an $n$ which is clear from the context, by $\boldsymbol{\sigma}(x)$ we mean the following tuple:
			\begin{equation*}
				\boldsymbol{\sigma}(x)=\left(\sigma^0(x),\sigma(x), \dots, \sigma^n(x)\right).
			\end{equation*}
			For a multi-index $J=(j_0,j_1, \dots j_n)$, by $\boldsymbol{\sigma}^{J}(x)$ we mean 
			\begin{align*}
				\boldsymbol{\sigma}^{J}(x)&=x^{j_0}\cdot \sigma(x)^{j_1}\cdots \left(\sigma^n(x)\right)^{j_n}\\
				&=x^{j_0+j_1\sigma+\dots+j_n\sigma^n}.
			\end{align*}
		\end{notn}
		\begin{rem}\label{sigmaJ}
			Let $f$ be a difference polynomial in one variable. Using Notation \ref{sigmabold}, $f$ is of this form
			\begin{equation*}
				f(x)=\underset{J\in \Lambda }{\mathlarger{\sum}}c_J \boldsymbol{\sigma}^J(x),
			\end{equation*}
			where $\Lambda$ is a finite subset of $\mathbb{N}^{n+1}$.
		\end{rem}
		\begin{rem}
			As an example, the difference monomial $x\sigma(x)$ can be written as $\boldsymbol{\sigma}^{(1,1)}(x)$ using Notation \ref{sigmabold}, and it can also be expressed as $x^{1+\sigma}$ using Notation \ref{sigma}. Throughout this paper, for difference polynomials in one variable, we will switch between these two notations. 
		\end{rem}
		\begin{rem}\label{Jmax}
			Suppose $f \in K_{\sigma}[x^{\pm1}]$ is a Laurent difference polynomial. Then it is of the form $f(x)= \underset{J \in \Lambda}{\sum}c_J\boldsymbol{\sigma}^{J}(x)$ where $\Lambda$ is a finite subset of $\mathbb{Z}^{n+1}$.\\
			Given $\Lambda$ as above, define $J_{\max}$ to be the multi-index such that 
			\begin{equation*}
				\forall i,\,\,\,\,\, 0\leq i \leq n, \,\,\,\,\, \left(J_{\max} \right)_i=	\begin{cases}
					0  & \text{ if }\,\, \forall J \in \Lambda\,\, j_i\geq 0 ,\\
					\max \left\lbrace \lvert j_i \rvert\,\,|\,\, j_i<0 \right\rbrace & \text{ if }\,\,   j_i<0 \text{ for some }\,\, J\in \Lambda.
				\end{cases}
			\end{equation*}
			Multiplying $f(x)$ by $\boldsymbol{\sigma}^{J _{\max}}(x)$ gives a difference polynomial $g(x)$. In other words
			\begin{equation*}
				f(x)\cdot \boldsymbol{\sigma}^{J_{\max}}(x)=g(x)\in K_{\sigma}[x].    
			\end{equation*}
			\par
			If $f$ is a Laurent difference polynomial in difference variables $x=\left(x_1, \dots, x_n\right)$, using Notation \ref{sigma}, it is of the following form:
			\begin{equation*}
				f(x)=\underset{u(\sigma)\in \Lambda }{\sum}c_{u(\sigma)}x^{u(\sigma)},
			\end{equation*} 
			where $\Lambda$ is a finite subset of $\left(\mathbb{Z}[\sigma]\right)^n$.
			\par Suppose $u(\sigma)=\left(u_1(\sigma), \dots, u_n(\sigma)\right)$ is one of the exponents appearing in $f$. For each $i$, $1\leq i \leq n$, we have 
			\begin{equation*}
				u_i(\sigma)= \sum_{j_i=0}^{m_i}a_{j_i}\sigma^{j_i},
			\end{equation*}
			where for each $j_i$, $a_{j_i}\in \mathbb{Z}$.
			Define $|u_i|_{\circ}(\sigma)=\sum\limits _{j_i=0}^{m_i}\tilde{a}_{j_i}\sigma^{j_i}$ such that $\tilde{a}_{j_i}=|a_{j_i}|$ if $a_{j_i}$ is negative and  $\tilde{a}_{j_i}=0$ otherwise.
			Set $|u|_{\circ}(\sigma)=\left(|u_1|_{\circ}(\sigma), \dots, |u_n|_{\circ}(\sigma)\right)$.
			Multiplying $f(x)$ by $\underset{u(\sigma)\in \Lambda}{\prod}x^{|u|_{\circ}(\sigma)}$ gives a difference polynomial $g(x)$. In other words, we have 
			\begin{equation*}
				f(x)\cdot \underset{u(\sigma)\in \Lambda}{\prod}x^{|u|_{\circ}(\sigma)}=g(x)\in K_{\sigma}[x_1, \dots,x_n].
			\end{equation*}
		\end{rem}
		\begin{defn}\label{f_J}
		\sloppy	To any difference polynomial in a single variable, $f(x)=\underset{J\in \Lambda }{\mathlarger{\sum}}c_J \boldsymbol{\sigma}^J(x)$, a polynomial $P(\boldsymbol{x})=\underset{J\in \Lambda}{\mathlarger{\sum}}c_J \boldsymbol{x}^J$ is associated, where $\boldsymbol{x}=(x_0,\dots,x_n)$, so that $f(x)=P(\boldsymbol{\sigma}(x))$.\\
			We use the notation $I!:=i_0!i_1!\dots i_n!$  with $I=\left(i_0, \dots, i_n\right)$. Then for a multi-index $I$ and a point $a$, $f_{(I)}(a)$ is defined as follows:
			\begin{equation*}
				f_{(I)}(a)=P_{(I)}(\boldsymbol{\sigma}(a))=\frac{\partial^{|I|}P\left(a,\sigma(a),\dots ,\sigma^n(a) \right) }{\partial x_0^{i_0} \partial x_1^{i_1}\dots \partial x_n^{i_n}}\cdot \frac{1}{I!}.
			\end{equation*}
			Note that for any multi-index $I$, we have $f_{(I)}(0)=c_{I}$.
			\par Similarly, to any Laurent difference polynomial $f$, a Laurent polynomial $P$ is associated.
		\end{defn}
		
		\begin{defn}\label{ACFA}
			A difference field $\left(k, \sigma \right)$  is called an \textit{ACFA}, if for any finite system of difference polynomial equations over $k$ with a solution in an extension $k'$ of $k$, this system has a solution in $k$. 
			This concept can be interpreted as a difference version of algebraically closed fields.
		\end{defn}
		\par To know more on this topic, see \cite{Zoe}.
		
		\begin{rem}
			In the same way as one constructs the algebraic closure of a field, one
			constructs a difference algebraic closure of a difference field;
			this is an ACFA.
			\par Unlike the algebraic closure of a field, a difference algebraic closure of a difference field is not unique. 
			\par It is clear from the definition of an ACFA that any ACFA is an algebraically closed field.
		\end{rem}
		
		\begin{defn}\label{axiom}
			Let $(K,\sigma)$ be a difference field. Then $\sigma$ canonically extends to an automorphism of $K[x_1, \dots, x_n]$ ( for all $i,\,\, 1\leq i \leq n\,\,\ ;\,\,\, \sigma(x_i)=x_i$). It is denoted by the same notation as $\sigma$.
			\par Let $K$ be an algebraically closed field. By a variety, we mean an irreducible Zariski closed subset of $K^n$. 
			\par Suppose $U$ is a variety over $K$. For $I(U)=\lbrace f\in K[x_1, \dots, x_n]\,\,|\,\, f(U)=0 \rbrace$, $\sigma \left( I(U)\right)$ is defined as follows:
			\begin{equation*}
				\sigma \left( I(U)\right)=\lbrace \sigma(f)\,\,|\,\, f \in I(U) \rbrace.
			\end{equation*}
			In this case, the variety $V\left(\sigma\left(I(U)\right)\right)$ is denoted by $U^{\sigma}$.  			
			\par Suppose $V$ is another variety over $K$, such that $V \subseteq U \times U^{\sigma} $, and two projection maps are given as follows:
			\begin{align*}
				& \pi_1:U \times U^{\sigma}\rightarrow U\\
				& \pi_2: U \times U^{\sigma}\rightarrow U^{\sigma}.
			\end{align*}
			Then the projection of $V$ to $U$ (to $U^{\sigma}$) is called \textit{generically onto}, if $\pi_1(V)$ ( if $\pi_2(V)$) is Zariski dense in $U$ (in $U^{\sigma}$). 
		\end{defn}
		The following theorem is not a result of this paper. It is derived from the axioms in $3.1$ and Theorem $3.2$ of \cite{Zoe}.
		\begin{thm}\label{axiomatization}
			Let $(k, \sigma)$ be a difference field. Then it is an ACFA if and only if for any two varieties $U$ and $V$ over $k$, with $V\subseteq U \times U^{\sigma}$, $k$ satisfies the following conditions:
			\begin{enumerate}
				\item $k$ is algebraically closed;
				\item If $V$ projects generically onto $U$ and $U^{\sigma}$, then there exists a point $a=(a_1,\dots,a_n)\in k^n$, such that $(a,\sigma(a))\in V$.
			\end{enumerate}
		\end{thm}
		\begin{proof}
			For the proof, see Theorem $3.2$ in \cite{Zoe}.
		\end{proof}
		\begin{lem}\label{ACFAroot}
			Let $k$ be an ACFA. Suppose $f$ is a Laurent difference polynomial in one variable with coefficients from $k$ that is not a monomial. Then $f$ has a nonzero root in $k$.
		\end{lem}
		\begin{proof}Let $f$ be a Laurent difference polynomial, such that $\sigma^n(x)$ is the greatest iteration of $\sigma$ appearing in $f$. If $n=0$, $f$ is a Laurent polynomial in one variable with coefficients from an ACFA. Since any ACFA is an algebraically closed field and $f$ is not a monomial, it has a nonzero root in $k$. Therefore, we assume that $n \neq 0$. In this case, we say that the order of $f$ is $n$. From Definition \ref{f_J}, there exists a Laurent polynomial $P$ in $k[y^{\pm1}_0, \dots, y^{\pm1}_n]$, for which we have
			\begin{equation*}
				f(x)=P(x,\sigma(x),\dots, \sigma^n(x)).
			\end{equation*}
			By Remark \ref{Jmax}, we can multiply $f$ by $\boldsymbol{\sigma}^{J _{\max}}(x)$, and convert it to a difference polynomial in $k_{\sigma}[x]$. Therefore, from now on, we assume that $f$ is in $k_{\sigma}[x]$, and is irreducible.\\
			We define the three following sets:
			\begin{align*}
				U&=\lbrace(\underline{r}, s_0)\,\,|\,\,P(\underline{r})=0, r_0\cdots r_ns_0-1=0\rbrace;\\
				\sigma(U)&=\lbrace (\underline{t},w_0)\,\,|\,\,\sigma(P)(\underline{t})=0, t_0\cdots t_nw_0-1=0\rbrace;\\
				V&=\lbrace (\underline{r},s_0,\underline{t},w_0)\in U\times \sigma(U)\,\,|
				\,\,t_i=r_{i+1}\,\, \text{for}\,\, 0\leq i \leq n-1 \rbrace,
			\end{align*}
			where $\underline{r}=(r_0, \dots , r_n)$, and $\underline{t}=(t_0, \dots,t_n)$.
			\par To find a nonzero root of $f$, we apply Theorem \ref{axiomatization}. If we prove that $V$ projects generically onto $U$ and $\sigma(U)$, then from this theorem, there exists a point $a=(a_0,\dots,a_n,b_0)$ such that $\left(a,\sigma(a)\right)\in V$.\\
			Finally, we prove that $a_0$ is a root of $f$. 
			\par Now, we prove that $V$ projects generically onto $U$.
			Since $f$, and consequently $P$, is irreducible, and as the order of $f$ is $n\neq 0$, there are at least two monomials in $\sigma(P)$ with different powers of $y_n$. Suppose $l$ and $m$ are two different powers of $y_n$ appearing in $\sigma(P)$. Regarding $\sigma(P)$ as a polynomial in $y_n$ with coefficients in $k[y_0, \dots, y_{n-1}]$, it can be written as
			\begin{equation*}
				\sigma(P)= \underset{i\in I}{\sum}g_i(y_0,\dots, y_{n-1})y_n^{i},
			\end{equation*}
			where $I$ is a finite subset of $\mathbb{N}_0$. Define
			\begin{equation}\label{g_lg_m}
				h(y_0, \dots ,y_{n-1})=g_l(y_0, \dots ,y_{n-1})\cdot g_m(y_0, \dots ,y_{n-1}).
			\end{equation}
			Using the following claim, we prove that $V$ projects generically onto $U$. Similarly it can be proven that $V$ projects generically onto $\sigma(U)$.
			\begin{claim}
				Let $(\underline{r}, s_0)$ be in $U$, such that $h(r_1, \dots, r_n)\neq 0$. Then $(\underline{r}, s_0)$ is in $\pi_1(V)$. 
			\end{claim}
			\renewcommand\qedsymbol{$\blacksquare$}
			\begin{proof} 
				In order to show that $(\underline{r}, s_0)\in \pi_1(V)$, we look for a point $(\underline{t}, w_0)$, such that $(\underline{r}, s_0,\underline{t}, w_0)$ is in $V$. For any $i$, $0 \leq i \leq n-1$, set $t_i=r_{i+1}$. Therefore, we have 
				\begin{equation*}
					h(t_0, \dots, t_{n-1})=h(r_1, \dots, r_n) \neq 0.
				\end{equation*}
				From \eqref{g_lg_m}, we obtain $g_l(t_0, \dots, t_{n-1})\neq 0$, and $g_m(t_0, \dots, t_{n-1})\neq 0$.
				This means that $\sigma(P)(t_0,\dots, t_{n-1},y_n)$ is a polynomial in $y_n$ with coefficients in $k$, with at least two monomials. Since $k$ is an ACFA, it is an algebraically closed field. Hence, this polynomial has a nonzero root $t_n$ in $k$. This gives the nonzero root $(t_0, \dots, t_n)$ of $\sigma(P)$. Also set $w_0=\dfrac{r_0s_0}{t_n}$. Then, we have
				\begin{align*}
					&t_0\cdots t_{n-1}t_nw_0-1=r_1\cdots r_n t_n\dfrac{r_0s_0}{t_n}-1\\
					&=r_0 \cdots r_n s_0 -1=0.
				\end{align*}
				Thus, $(\underline{t}, w_0)\in \sigma(U)$, and consequently $(\underline{r}, s_0,\underline{t}, w_0)\in V$.

			\end{proof}
			
			\begin{rem}
				Let $D(h)=k^{n+2}\setminus V(h)$ be the Zariski open set defined by $h$. Then $D(h) \cap U$ is nonempty.
			\end{rem}
			\begin{proof}
				As we discussed on the previous page, there are at least two monomials in $\sigma(P)$ with different powers of $y_n$. Similarly, $P$ has at least two monomilas with different powers of $y_n$. 
				Suppose $l$ and $m$ are two different powers of $y_n$ appearing in $P$. Regard $P$ as a polynomial in one variable $y_n$ with coefficients in $k[y_0, \dots, y_{n-1}]$. Then $P$ is of the following form:
				\begin{equation*}
					P= \underset{i\in I}{\sum}q_{i}(y_0,\dots, y_{n-1})y_n^{i},
				\end{equation*}
				where $I$ is a finite subset of $\mathbb{N}_0$.
				Define 
				\begin{equation*}
					h'(y_0, \dots y_{n-1})=q_{l}(y_0, \dots y_{n-1}) \cdot q_{m}(y_0, \dots y_{n-1})\cdot y_0 \cdots y_{n-1}.
				\end{equation*} 	
				Since $h$ and $h'$ are nonzero polynomials, $h\cdot h'$ has a nonroot $(r_0,\dots r_{n-1})$. Therefore, we have $h(r_0, \dots,r_{n-1} )\neq 0$, which means that this is a point in $D(h)$. On the other hand, we have $h'(r_0, \dots,r_{n-1} )\neq 0$. This means that $P(r_0, \dots r_{n-1},y_n)$ is a nonzero polynomial with coefficients in $k$. Since $k$ is an algebraically closed field, and  $P$ (regarded as a polynomial in one variable $y_n$) has at least two monomials, it has a nonzero root $r_n$ in $k$. Hence, $(r_0, \dots, r_{n-1},r_n)$ is a root of $P$. From the definition of $h'$, it is clear that for all $i$, $1 \leq i \leq n-1$, we have $r_i\neq 0$. Set $s_0=\dfrac{1}{r_0 \cdots r_n}$, then $(r_0, \dots ,r_n,s_0)$ is a point in $U$, which is also a point of $D(h)$. Thus, $D(h) \cap U$ is nonempty. 
				
			\end{proof}
			\renewcommand\qedsymbol{$\square$}
			From the previous claim, we have
			\begin{equation*}
				D(h) \cap U\subseteq \pi_1(V).
			\end{equation*}
			By taking the Zariski closure, we obtain the following relation:
			\begin{equation*}
				\overline{D(h)}\cap\overline{U}\subseteq \overline{\pi_1(V)}.
			\end{equation*}
			Since Zariski open sets are dense, we have $U\subseteq \overline{\pi_1(V)}$. Moreover, from $\overline{\pi_1(V)}\subseteq U$, we obtain $\overline{\pi_1(V)}=U$. This means that $\pi_1(V)$ is Zariski dense in $U$, and $V$ projects generically onto $U$.  
			With a similar argument, $V$ projects generically onto $\sigma(U)$. Hence, from Theorem \ref{axiomatization}, there exists a point $a=(a_0,\dots,a_n,b_0)$, such that $\left(a,\sigma(a)\right)\in V$, where $\sigma(a)=\left(\sigma(a_0), \dots \sigma(a_n),\sigma(b_0)\right)$, and $a_{i+1}=\sigma(a_i)$. This means that
			\begin{equation*}
				f(a_0)=P(a_0,\sigma(a_0), \dots \sigma^n(a_0))=P(a_0,a_1, \dots, a_n)=0.
			\end{equation*}
			Since $a\in U$, we have $a_0\neq 0$. Hence, $f$ has a nonzero root $a_0$.
			
		\end{proof}
		
		\begin{thm}\label{alpharoot}
			Let $k$ be an ACFA. Suppose $f$ is in $k_{\sigma}[x_1^{\pm 1},\dots,x_n^{\pm 1}]$, and is not a monomial. Then $f$ has a root in $(k^*)^n$.
		\end{thm}
		\begin{proof}
			Remark \ref{Jmax} enables us to convert a Laurent difference polynomial to a difference polynomial. Therefore, we assume that $f$ is an element of $k_{\sigma}[x_1, \dots, x_n]$.\\
			\sloppy Since $f$ is not a monomial, we assume that $c_{u(\sigma)}x_1^{u_1(\sigma)}\cdots x_n^{u_n(\sigma)}$ and $c_{u^{\prime}(\sigma)}x_1^{u^{\prime}_1(\sigma)}\cdots x_n^{u^{\prime}_n(\sigma)}$ are two distinct monomials of $f$. Therefore, for some $i$, $1\leq i \leq n$, we have $u_i(\sigma)\neq u^{\prime}_i(\sigma)$. Without loss of generality, we assume that $i=n$, and we have at least two monomials with different $\sigma$-powers of $x_n$.\\
			Regard $f$ as a difference polynomial in one variable $x_n$, with coefficients in $k_{\sigma}[x_1, \dots, x_{n-1}]$. We write $f$ as follows:
			\begin{equation*}
				f=\sum_{k=1}^{N}g_k\left(x_1,\dots,x_{n-1}\right)x^{u_{(n,k)}(\sigma)}_n,
			\end{equation*}
			where $u_{(n,k)}(\sigma)$ denotes distinct $\sigma$-powers of $x_n$ appearing in $f$. \\
			Define $h\left(x_1,\dots,x_{n-1}\right)=g_1\left(x_1,\dots,x_{n-1}\right)\cdots g_N\left(x_1,\dots,x_{n-1}\right)$. The difference polynomial $h$ is nonzero. In the following claim, we want to find $\left(a_1, \dots,a_{n-1}\right)\in (k^*)^{n-1}$ such that $h\left(a_1, \dots,a_{n-1}\right)\neq 0$. 
			In this case, $f\left(a_1,\dots,a_{n-1},x_n\right)$ is a difference polynomial in one variable with coefficients from an ACFA which is not a monomial. Therefore, from Lemma\ref{ACFAroot} it has a nonzero root $a_n$. Hence, $\left(a_1, \dots,a_{n-1},a_n\right)$ is a root of $f$ in $(k^*)^n$, and the proof is complete.
			\begin{claim}
				Let $k$ be an ACFA. Suppose $h$ is a nonconstant difference polynomial in $k_{\sigma}[x_1, \dots,x_{n-1}]$. Then $h$ has a nonroot in $(k^*)^{n-1}$.
			\end{claim}
			\renewcommand\qedsymbol{$\blacksquare$}
			\begin{proof}
				We prove this claim by induction over $n-1$. Suppose $n-1=1$, and consider the difference polynomial $h(x_1)+b$, where $0 \neq b\neq -c_0$ for $c_0$ being the constant term of $h$. This is an element of $k_{\sigma}[x_1]$, which is not a monomial. Since $k$ is an ACFA, Lemma \ref{ACFAroot} implies that this polynomial has a nonzero root $\alpha_1$. This means that $h(\alpha_1)+b=0$, and consequently $h(\alpha_1)\neq 0$.\\
				We assume that the statement holds, if the number of variables is less than $n-1$.
				We prove the statement for $h$ in $k_{\sigma}[x_1, \dots,x_{n-1}]$.
				We write $h$ as a polynomial in one variable $x_{n-1}$, with coefficients from $k_{\sigma}[x_1,\dots,x_{n-2}]$. More precisely, $h$ is of the following form:
				\begin{equation*}
					h=\sum_{k=1}^{M}h_k\left(x_1,\dots,x_{n-2}\right)x^{u_{(n-1,k)}(\sigma)}_{n-1}.
				\end{equation*}	
				For each $k$, $1\leq k\leq M$, we have $h_k\in k_{\sigma}[x_1,\dots,x_{n-2}] $. By induction assumption, for each $k$, there exists $\left(\alpha_{(k,1)},\dots, \alpha_{(k,n-2)}\right)$ in $(k^*)^{n-2}$ that is a nonroot of $h_k$. Therefore, there exists a point $\left(a_1,\dots, a_{n-2}\right)$ such that $h\left(a_1,\dots, a_{n-2},x_{n-1}\right)$ is a nonzero difference polynomial in one variable $x_{n-1}$. Thus, from the first step of the induction, it has a nonroot $a_{n-1}$. This means that $\left(a_{1},\dots, a_{n-2},a_{n-1}\right)$ is a nonroot of $h$.      
			\end{proof} 
			\renewcommand\qedsymbol{$\square$}
			As explained before this claim, using the nonroot $\left(a_{1},\dots,a_{n-1}\right)$ of $h$ from the previous claim, we find a root $\left(a_{1},\dots,a_{n-1},a_n\right)\in (k^*)^n$ of $f$.
		\end{proof}
		
		\begin{defn}
			Let $(K,v)$ be a valued field. We consider a totally ordered collection 
			$\left\lbrace B_i \right\rbrace_{i\in I}$ of balls in $K$. Then $K$ is called \textit{spherically complete}, if for every such collection we have $\underset{i\in I}{\bigcap} B_i\neq \emptyset$.\\
			Hahn fields are spherically complete, see \cite{Krull} and \cite{Kaplansky}.
		\end{defn}
		\begin{assum}(\textit{General assumptions})\label{assump}
			Here are some of the general assumptions we make throughout the present work:
			\begin{itemize}
				\item $(K,\sigma)$ is a multiplicative valued difference field that is spherically complete and of characteristic zero. We also assume that $\rho$, the scaling exponent of $\sigma$ is transcendental.(The importance of this assumption is explained in Definition \ref{difftrop}.)  
				
				\item The valuation has a splitting. This splitting $\psi:\Gamma \rightarrow K$ interacts with $\sigma_{\Gamma}$ as follows:
				\begin{equation*}
					\forall a\in \Gamma: \psi\left(\sigma_{\Gamma}(a)\right)=\sigma \left(\psi(a)\right).
				\end{equation*}
				We also use the notation $\psi(a)=t^{a}$.
				\item The difference residue field is an ACFA and is of characteristic zero.
				\item The difference value group $\Gamma$ is a subgroup of $\mathbb{R}$ that is a $\mathbb{Q}(\rho)$-module.
			\end{itemize}
		\end{assum}
		\begin{exmp}
			Let $k$ be an ACFA of characteristic zero. Assume $\mathbb{R}$ is regarded as an ordered abelian group with an automorphism $\sigma_{\Gamma}$ such that
			\begin{equation*}
				\forall x\in \mathbb{R}:\,\,\, \sigma_{\Gamma}(x)=\rho \cdot x
			\end{equation*}
			where $\rho$ is a fixed positive real number which is transendental. Then $K=k((t^{\mathbb{R}}))$ satisfies the assumptions in Assumption \ref{assump}.
		\end{exmp}
		
		\subsubsection{Difference Tropical Objects}   
		Here, we introduce some difference tropical objects. For clarity, we also present some examples. The reader can find classical analogues of this material in \cite{Diane}.
		\begin{defn}\label{difftrop}
			Take a Laurent difference polynomial $f \in K_\sigma[x_1^{\pm1},\dots,x_n^{\pm1}]$. So $f$ can be written as $f(x) =\underset{u(\sigma)\in(\mathbb{Z}[\sigma])^n}{\mathlarger{\sum}} c_{u(\sigma)}x^{u(\sigma)}$. To obtain the   \textit{tropicalization of $f$}, we replace classical addition and multiplication operations with tropical ones, $\bigoplus$ and $\bigodot$, replace the coefficients with their respective valuations, and replace $\sigma$ with the induced automorphism on the value group of $K$, denoted by $\sigma_{\Gamma}$. In other words, the tropicalization of $f$ is defined as follows:
			\begin{align*}
				\mathoptrop(f)(w)&=
				\underset{u(\sigma_{\Gamma})\in(\mathbb{Z}[\sigma_{\Gamma}])^n}{\bigoplus} 
				\left( v(c_{u(\sigma)})\odot w^{u(\sigma_{\Gamma})} \right)\\
				&=	\underset{u(\sigma_{\Gamma})\in(\mathbb{Z}[\sigma_{\Gamma}])^n}{\min} 
				\left( v(c_{u(\sigma)})+u(\rho)\cdot w \right).    
			\end{align*}
			Here by $ w^{u(\sigma_{\Gamma})}$ we mean $\underset{i \in \left\lbrace 1,\dots,n\right\rbrace}{\bigodot}w_i^{u_i(\sigma_{\Gamma})}$, for $u(\sigma_{\Gamma})=\left(u_1(\sigma_{\Gamma}),\dots ,u_n(\sigma_{\Gamma})\right)$. By using Notation \ref{sigma}, we have $w_i^{u_i(\sigma_{\Gamma})}=u_i(\sigma_{\Gamma})(w_i)=u_i(\rho)\cdot w_i$.
			\par By tropicalizing a Laurent difference polynomial $f$, we obtain a tropical polynomial $\mathoptrop(f)$. Each tropical monomial appearing in $\mathoptrop(f)$ corresponds to a difference monomial in $f$. In order to keep this correspondence and avoid possible annihilation of some tropical monomials, we suppose the automorphism of $K$ is such that its scaling exponent, $\rho$ is transcendental.
		\end{defn}
		\begin{rem}
			Let $f$ be a difference polynomial in one variable. Using the notation defined in Remark \ref{sigmaJ}, the tropicalization of $ f(x)=\underset{J\in \Lambda }{\mathlarger{\sum}}c_J \boldsymbol{\sigma}^J(x)$ is of the following form:
			\begin{equation*}
				\mathoptrop(f)(w)=\underset{J\in \Lambda}{\min}\left\lbrace v(c_J)+J \boldsymbol{\sigma}_{\Gamma}(w) \right\rbrace,
			\end{equation*}
			where by $\boldsymbol{\sigma}_{\Gamma}(w)$, we mean $
			\left(w,\sigma_{\Gamma}(w), \dots, \sigma_{\Gamma}^n(w)\right)$, and $\Lambda$ is a finite subset of $\mathbb{N}^{n+1}$.
		\end{rem}
		\begin{exmp}\label{extrop}
			Let $f \in \mathbb{C}((t^{\mathbb{R}}))_{\sigma}[x_1^{\pm 1},x_2^{\pm 1}]$ be 
			\begin{equation*}
				f(x_1,x_2)=(1+t)x_1x_2^{\sigma^3}+t^2x_2^{\sigma}+1,
			\end{equation*}
			and $\rho=\pi$. The tropicalization of $f$ is 
			\begin{equation*}
				\mathoptrop(f)(w_1,w_2)=\min \lbrace w_1+\rho^3\cdot w_2, 2+\rho\cdot w_2,0 \rbrace= \min \lbrace w_1+\pi^3 w_2,2+\pi w_2,0 \rbrace.
			\end{equation*}
		\end{exmp}
		\begin{defn}\label{trophyp}
			Let $f$ be a Laurent difference polynomial in $n$ variables with coefficients from $K$. We say $w$ is a \textit{tropical root of $f$}, if in $\mathoptrop(f)(w)$ the minimum is attained at least twice.
			\par A \textit{difference tropical hypersurface}, which is denoted by $\mathoptrop(V(f))$, is the set of all tropical roots of a Laurent difference polynomial $f$. In other words, we have 
			\begin{equation*}
				\mathoptrop(V(f))=\lbrace w\in \mathbb{R}^n\,\,|\,\,\text{in}\,\, \mathoptrop(f)(w)\,\,\text{the minimum is attained at least twice}\rbrace.
			\end{equation*}
		\end{defn}
		\begin{defn}\label{diffinitial}
			Suppose $f(x) =\underset{u(\sigma)\in(\mathbb{Z}[\sigma])^n}{\mathlarger{\sum}} c_{u(\sigma)}x^{u(\sigma)}$ is a Laurent difference polynomial in $K_\sigma[x_1^{\pm1},\dots,x_n^{\pm1}]$.\\
			The \textit{initial form} of $f$ with respect to the point $w \in \mathbb{R}^n$ is a difference polynomial with coefficients in $\mathbf{k}$ that is defined as follows: 
			\begin{equation*}
				\mathopin_w(f)=\underset{\substack{u(\sigma):v(c_{u(\sigma)})+u(\rho)\cdot w\\ =\mathoptrop(f)(w)}}{\sum} \overline{c_{u(\sigma)} t^{-v(c_{u(\sigma)})}}\cdot x^{u(\Bar{\sigma})}.
			\end{equation*}
			In other words, to obtain $\mathopin_w(f)$, we consider those monomials of $\mathoptrop(f)$ which achieve the minimum at $w$. Corresponding to each of these monomials, a monomial appears in $\mathopin_w(f)$. 	 	
		\end{defn}
		\begin{rem}
			Let $f$ be a difference polynomial in one variable. If we use the notation defined in Remark \ref{sigmaJ}, then the initial form of $ f(x)=\underset{J\in \Lambda }{\mathlarger{\sum}}c_J \boldsymbol{\sigma}^J(x)$ is of the following form:
			\begin{equation*}
				\mathopin_w(f)(x)=\underset{\substack{J:v(c_J)+|J|_{\rho}\cdot w\\=\mathoptrop(f)(w)}}{\sum}\overline{t^{-v(c_J)}c_J}\boldsymbol{\bar{\sigma}}^{J}(x).
			\end{equation*}
			
		\end{rem}
		\begin{exmp}
			Take $f(x_1,x_2)=(1+t)x_1x_2^{\sigma^3}+t^2x_2^{\sigma}+1$ from Example \ref{extrop}. The initial form of $f$ with respect to  $w=(3\pi^2,\frac{-2}{\pi})$ is as follows:
			\begin{equation*}
				\mathopin_{w}(f)=\Bar{1}x_2^{\Bar{\sigma}}+\Bar{1}.
			\end{equation*}
		\end{exmp}
		From the definition of the initial form the following lemma can be easily deduced.
		\begin{lem}\label{roots}
			Suppose $f$ is a Laurent difference polynomial in $K_\sigma[x_1^{\pm1},\dots,x_n^{\pm1}]$. Then $w$ is a tropical root of $f$ if and only if $\mathopin_w(f)$ is not a monomial. 
		\end{lem}
		\begin{proof}
			$\Longrightarrow)$ Since $f$ is an element of $K_{\sigma}\left[x_1^{\pm1},\dots,x_n^{\pm1}\right]$, it is  of the form\\
			$f=\underset{u(\sigma)\in (\mathbb{Z}[\sigma])^n}{\mathlarger{\sum}}c_{u(\sigma)}x^{u(\sigma)}$. Suppose $w=(w_1,\dots,w_n)$ is a tropical root of $f$. By definition, we have
			\begin{equation*}
				\mathoptrop(f)(w)= \underset{u(\sigma) \in (\mathbb{Z}[\sigma])^n }{\min}(v(c_{u(\sigma)})+u(\rho)\cdot w),
			\end{equation*}
			and $w$ attains the minimum in at least two different tropical monomials each of which corresponds to a monomial in $\mathopin_w(f)$. Therefore $\mathopin_w(f)$ has at least two monomials, meaning that $\mathopin_w(f)$ is not a monomial.\\
			$\Longleftarrow)$ If $\mathopin_w(f)$ is not a monomial, from the definition of the initial form, each of its monomials corresponds to a tropical monomial in $\mathoptrop(f)(w)$ each of which attains the minimum. Thus, $\mathoptrop(f)(w)$ achieves the minimum for at least two different tropical monomials. This means $w$ is a tropical root of $f$.
		\end{proof}
		The two following lemmas are well-known in the context of tropical ploynomials, and the generalization that $f$ and $g$ are Laurent difference polynomials does not affect their validity. Therefore, we can state them in terms of Laurent difference polynomials.
		\begin{lem}\label{trop(fg)}
			Let $f$ and $g$ be two Laurent difference polynomials in $K_{\sigma}[x_1^{\pm1},\dots,x_n^{\pm1}]$. Suppose $w$ is an element of $\Gamma$, then we have 
			\begin{equation*}
				\mathoptrop(fg)(w)=\mathoptrop(f)(w)+\mathoptrop(g)(w).
			\end{equation*}
		\end{lem}
		\begin{shortversion}
			\begin{proof}
				As $f$ and $g$ are elements of $K_{\sigma}[x_1^{\pm1},\dots,x_n^{\pm1}]$, they can be written as\\
				$f=\underset{u \in (\mathbb{Z}[\sigma])^n}{\mathlarger{\sum}}c_u x^{u}$ and $g=\underset{u' \in (\mathbb{Z}[\sigma])^n}{\mathlarger{\sum}}d_{u'} x^{u'}$.\\
				Suppose $c_{u_i}x^{u_i}$ , $c_{u_j}x^{u_j}$ and $d_{u'_k}x^{u'_k}$, $d_{u'_l}x^{u'_l}$ are two monomials of $f$ and $g$ respectively, where $u_i+u'_k=u_j+u'_l=\nu$.\\
				
				If we consider tropicalization of $fg$, it will be the minimum of tropicalization of each monomial appearing in $fg$ such as 
				\begin{equation*}
					v(c_{u_i}d_{u'_k})+(u_i+u'_k)(\rho)\cdot w=v(c_{u_i})+v(d_{u'_k})+\nu (\rho)\cdot w
				\end{equation*} 
				and 
				\begin{equation*}
					v(c_{u_j}d_{u'_l})+(u_j+u'_l)(\rho)\cdot w=v(c_{u_j})+v(d_{u'_l})+\nu (\rho)\cdot w
				\end{equation*}
				So for such two tropical monomials we have
				\begin{align*}
					& \min\left\lbrace v(c_{u_i})+v(d_{u'_k})+\nu (\rho)\cdot w , v(c_{u_j})+v(d_{u'_l})+\nu (\rho)\cdot w \right\rbrace\\
					&=\min \left\lbrace v(c_{u_i})+v(d_{u'_k}) , v(c_{u_j})+v(d_{u'_l}) \right\rbrace +\nu (\rho)\cdot w
				\end{align*}
				Generally speaking for any fixed $\nu$ appearing in $fg$ we have
				\begin{align*}
					& \underset{u+u'= \nu}{\min} \left\lbrace v(c_{u})+ v(d_{u'})+\nu (\rho)\cdot w \right\rbrace\\
					&=  \underset{u+u'= \nu}{\min} \left\lbrace v(c_{u})+ v(d_{u'})\right\rbrace +\nu (\rho)\cdot w 
				\end{align*}		
				So if we consider $\mathoptrop(fg)(w)$ we will have 
				\begin{equation*}
					\mathoptrop(fg)(w)= \underset{\nu}{\min}\left\lbrace \underset{u+u'=\nu}{\min}\left\lbrace v(c_u)+v(d_{u'})\right\rbrace + \nu(\rho)\cdot w  \right\rbrace 
				\end{equation*}
				Secondly we have 
				\begin{equation*}
					\mathoptrop(f)(w)+\mathoptrop(g)(w)=\underset{u}{\min}\left\lbrace v(c_u)+u(\rho)\cdot w \right\rbrace+\underset{u'}{\min}\left\lbrace v(d_{u'})+u'(\rho)\cdot w \right\rbrace 
				\end{equation*}
				Suppose in $\mathoptrop(f)(w)$ the minimum is obtained at $u_i$ and in $\mathoptrop(g)(w)$ the minimum is obtained at $u'_j$.\\
				If we add up all tropical monomials of $f$ with all tropical monomials of $g$ and then take the minimum, it is obtained at $u_i+u'_j$.\\
				In other words we have
				\begin{equation}\label{trop+}
					\begin{aligned}
						\mathoptrop(f)(w)+\mathoptrop(g)(w)&=v(c_{u_i})+u_i(\rho)\cdot w+v(d_{u'_j})+u'_j(\rho)\cdot w\\
						&=\min\left\lbrace v(c_u)+u(\rho)\cdot w+v(d_{u'})+u'(\rho)\cdot w \right\rbrace\\
						&=\min \left\lbrace v(c_u)+v(d_{u'})+(u(\rho)+u'(\rho))\cdot w\right\rbrace\\
						&=\underset{\underset{u+u'=\nu}{\nu}}{\min}\left\lbrace v(c_u)+v(d_{u'})+\nu(\rho)\cdot w \right\rbrace\\
						&=\underset{\nu}{\min}\left\lbrace \underset{u+u'=\nu}{\min}\left\lbrace v(c_u)+v(d_u') \right\rbrace +\nu(\rho)\cdot w \right\rbrace\\
						&=\mathoptrop(fg)(w) 
					\end{aligned}  
				\end{equation} 
			\end{proof}
		\end{shortversion}
		\begin{lem}\label{in(fg)}
			Let $f$ and $g$ be two Laurent difference polynomials in $K_{\sigma}[x_1^{\pm1},\dots ,x_n^{\pm1}]$. Suppose $w$ is an element of $\Gamma$, then we have
			\begin{equation*}
				\mathopin_w(fg)=\mathopin_w(f)\mathopin_w(g).
			\end{equation*}
		\end{lem}
		\begin{shortversion}
			\begin{proof}
				We have $f=\underset{u(\sigma) \in (\mathbb{Z}[\sigma])^n}{\mathlarger{\sum}}c_{u(\sigma)} x^{u(\sigma)}$ and $g=\underset{u'(\sigma) \in (\mathbb{Z}[\sigma])^n}{\mathlarger{\sum}}d_{u'(\sigma)} x^{u'(\sigma)}$. 
				Then $fg$ can be written as 
				\begin{equation*}
					fg=\underset{\nu(\sigma) \in (\mathbb{Z[\sigma]})^n}{\sum}e_{\nu(\sigma)}x^{\nu(\sigma)}\,\,\text{where}\,\, e_{\nu(\sigma)}=\underset{u(\sigma)+u'(\sigma)=\nu(\sigma)}{\sum}c_{u(\sigma)} d_{u'(\sigma)}   
				\end{equation*}
				
				From Lemma \ref{trop(fg)} we have $\mathoptrop(fg)(w)=\mathoptrop(f)(w)+\mathoptrop(g)(w)$. Therefore, we have 
				\begin{equation}
					\begin{aligned}
						\mathopin_w(fg)&=\underset{\substack{\nu(\sigma):v(e_{\nu(\sigma)})+ \nu(\rho)\cdot w=\mathoptrop(fg)(w)\\=\mathoptrop(f)(w)+\mathoptrop(g)(w)}}{\sum}\overline{e_{\nu(\sigma)}t^{-v(e_{\nu(\sigma)})}}x^{\nu(\sigma)}\\
						&=\underset{\substack{\nu(\sigma):v(e_{\nu(\sigma)})+ \nu(\rho)\cdot w\\=\mathoptrop(f)(w)+\mathoptrop(g)(w)}}{\sum}\,\,\, \overline{\underset{u(\sigma)+u'(\sigma)=\nu(\sigma)}{\sum}c_{u(\sigma)}d_{u'(\sigma)}t^{-\mathoptrop(f)(w)-\mathoptrop(g)(w)+\nu(\rho)\cdot w}}\cdot x^{\nu(\sigma)}\\
						&=\underset{\substack{\nu(\sigma):v(e_{\nu(\sigma)})+ \nu(\rho)\cdot w\\=\mathoptrop(f)(w)+\mathoptrop(g)(w)}}{\sum}\,\,\, \underset{u(\sigma)+u'(\sigma)=\nu(\sigma)}{\sum}\overline{c_{u(\sigma)}d_{u'(\sigma)}t^{-\mathoptrop(f)(w)-\mathoptrop(g)(w)+(u(\rho)+u'(\rho))\cdot w}}\cdot x^{\nu(\sigma)}
					\end{aligned}
				\end{equation}
				As in $\mathopin_w(fg)$ we are focusing on those $\nu$ for which in $\mathoptrop(fg)(w)$ the minimum is obtained and each $\nu$ in $fg$ can be written as $u+u'$ where $u$ corresponds to a tropical monomial in $\mathoptrop(f)(w)$ and $u'$ corresponds to a tropical monomial in $\mathoptrop(g)(w)$ as we calculated in \ref{trop+} we know that $u$ and $u'$ obtain the minimum in $\mathoptrop(f)(w)$ and $\mathoptrop(g)(w)$ respectively.Therefore we can write
				\begin{equation*}
					\begin{aligned}
						\mathopin_w(fg)&=\underset{\substack{\nu:v(e_\nu)+\nu(\rho)(w)\\ =\mathoptrop(f)(w)+\mathoptrop(g)(w)}}{\sum}\underset{u+u'=\nu}{\sum}\overline{c_u t^{-\mathoptrop(f)(w)+u(\rho)\cdot w}d_{u'}t^{-\mathoptrop(g)(w)+u'(\rho)\cdot w}}x^{u+u'}\\&=\underset{\substack{\nu:v(e_\nu)+\nu(\rho)(w)\\ =\mathoptrop(f)(w)+\mathoptrop(g)(w)}}{\sum}\underset{u+u'=\nu}{\sum}\overline{c_u t^{-v(c_u)}d_{u'}t^{-v(d_{u'})}}x^{u}x^{u'}\\
						&=\underset{\substack{\nu:v(e_\nu)+\nu(\rho)(w)\\ =\mathoptrop(f)(w)+\mathoptrop(g)(w)}}{\sum}\underset{u+u'=\nu}{\sum}\overline{c_u t^{-v(c_u)}}\cdot x^{u}\cdot \overline{d_{u'}t^{-v(d_{u'})}}x^{u'}\\
						&=\left(\underset{\substack{u:v(c_u)+u(\rho)(w)\\ =\mathoptrop(f)(w)}}{\sum} \overline{c_u t^{-v(c_u)}}.x^{u} \right) \cdot \left(\underset{\substack{u':v(d_{u'})+u'(\rho)(w)\\ =\mathoptrop(g)(w)}}{\sum} \overline{d_{u'}t^{-v(d_{u'})}}x^{u'} \right)\\ &=\mathopin_w(f)\cdot \mathopin_w(g)
					\end{aligned}
				\end{equation*} 
			\end{proof}
		\end{shortversion}
		\subsection{Polyhedral Geometry}\label{plygeo}
		In this subsection, we provide some preliminaries on Polyhedral Geometry. Readers who are familiar with the subject can skip ahead, while those who want to know more, can refer to the main references for this subsection, which are \cite{Ziegler, Diane, Ewald}. 
		\begin{defn}\label{plh}
			A \textit{polyhedron} is a subset of $\mathbb{R}^n$ that is the intersection of finitely many closed half spaces. It is usually denoted by $P$. More precisely, it can be described as:
			\begin{equation*}
				P=\left\lbrace x \in \mathbb{R}^n\,\, | \,\, \mathcal{A}\cdot x \leq \mathcal{B}  \right\rbrace, 
			\end{equation*}
			where $\mathcal{A} \in M_{d\times n}(\mathbb{R})$ and $\mathcal{B}\in \mathbb{R}^d$.
		\end{defn}
		A specific class of polyhedra are polytopes. A \textit{polytope} is a bounded polyhedron. Below, in Definition \ref{pltp}, we see an equivalent definition of a polytope. See Lecture $1$ of \cite{Ziegler}. 
		\begin{defn}\label{conv}
			\begin{shortversion}
				Let $A$ be a subset of $\mathbb{R}^n$, $A$ is called \textit{convex}, if for any two elements $a,b \in A$, the straight line segment connecting $a$ and $b$ is also in $A$. More precisely, if for $a,b \in A$ we have $\lambda a+(1-\lambda)b \in A$ where $0 \leq \lambda \leq 1$.
			\end{shortversion}
			\par	For $U \subseteq \mathbb{R}^n$ the \textit{convex hull} of $U$, denoted by $conv(U)$, is 
			\begin{shortversion}
				defined to be the intersection of all convex sets containing $U$, or equivalently is 
			\end{shortversion}the smallest convex set containing $U$.
		\end{defn}
		\begin{rem}
			The following definition is in fact the definition of a convex polytope.
			In the literature, a polytope can be convex or nonconvex, but everywhere in this paper, by a polytope, we mean a convex polytope. Therefore, in this definition, we omit the word "convex". 
		\end{rem}
		\begin{defn}\label{pltp}
			If $U$ is a finite subset of $\mathbb{R}^n$, namely $U=\left\lbrace u_1,\dots,u_r\right\rbrace $ then $conv(U)$ is called a \textit{polytope}. In this case, it can be described as:
			\begin{equation*}
				conv(U)=\left\lbrace \sum_{i=1}^{r} \lambda_i u_i\,\,| \,\, \forall i,\,\,\,\,\, 0 \leq \lambda_i \leq1 \,\,\,\,\,and\,\,\,\,\,\sum_{i=1}^{r} \lambda_i =1 \right\rbrace. 
			\end{equation*}
		\end{defn}
		\begin{notn}
			Both polytopes and polyhedra are usually denoted by $P$, and their distinction will be clear from the context. 
		\end{notn}
		\begin{shortversion}
			\begin{rem}
				Dropping the condition $ \forall i,\,\,\, 0 \leq \lambda_i \leq1$ in Definition \ref{pltp}, we obtain the \textit{affine hull (affine span)} of the set $U$ which is denoted by $\mathrm{aff}(U)$.\\
				More generally, affine subspaces which are translate of vector(also called linear) subspaces can be described as the affine hull of a finite set of points.
			\end{rem}
		\end{shortversion}
		\begin{defn}\label{face}
			A \textit{face} $F$ of a polyhedron $P$, is a set of the form 
			\begin{equation*}
				\left\lbrace x \in P\,\, | \,\, \forall y \in P\,\,\, v\cdot x \leq v\cdot y  \right\rbrace,  
			\end{equation*}
			where $v$ is a linear functional in $\left( \mathbb{R}^n\right)^{\vee}$. In this case, we say $F$ is determined by $v$, and we use the notation $F=\mathrm{face}_v(P)$. 
			\par The \textit{dimension of a face} $F$ is defined to be the dimension of the affine span of $F$, which is  denoted by $\mathrm{aff}(F)$. In particular, the \textit{dimension of a polyhedron} equals to the dimension of its affine span. The faces of dimension $0$ are called vertices and the set of all vertices of $P$ is denoted by $V_P$. The faces of dimension $1$ are called \textit{edges}.
			A face that is not included in any other proper face is called a \textit{facet}.
			Additionally, $\emptyset$ is a face of dimension $-1$. Note that faces are topologically closed sets.
		\end{defn}
		\begin{rem}
			Since a polytope is a bounded polyhedron, it is meaningful to talk about a face of a polytope. Moreover, a polytope $P$ is a face of itself, which is of maximal dimension.
		\end{rem}
		\begin{defn}\label{Newpolytope}
			Let $f(x)=\underset{u(\sigma)\in \left(\mathbb{Z}[\sigma]\right)^n}{\sum}c_{u(\sigma)}x^{u(\sigma)}$ be a Laurent difference polynomial. We can associate a polytope to $f(x)$. Assume $\rho$ is the scaling exponent of $\sigma$, which is transcendental. For any exponent $u(\sigma)$ appearing in $f$, consider $u(\rho)\in \mathbb{R}^n$. Then this polytope is defined as follows:
			\begin{equation*}
				P=conv(U)\,\, \text{where}\,\, U=\left\lbrace u(\rho)\,\,|\,\,c_{u(\sigma)}\neq 0 \right\rbrace.
			\end{equation*}
			This polytope is called the \textit{Newton polytope} associated to $f$.
		\end{defn}
		\begin{prop}\label{prop2.3}
			Let $P \subseteq \mathbb{R}^n$ be a polytope, and $F$ be a face of $P$. Then $F$ is also a polytope, and the set of its vertices is $V_F:=V_P \cap F$. More generally, the faces of $F$ are those faces of $P$ which are contained in $F$. 
		\end{prop}
		\begin{proof}
			The reader can find the proof in Proposition 2.3 of \cite{Ziegler}.    
		\end{proof}
		\begin{defn}\label{proj}
			Let $P\subseteq \mathbb{R}^n$ and $Q \subseteq \mathbb{R}^m$ be two polytopes. Then an affine map 
			\begin{align*}
				\pi:&\mathbb{R}^n \longrightarrow \mathbb{R}^m,\\
				& x \longmapsto Ax-z,
			\end{align*}
			with $A \in  M_{m\times n}(\mathbb{R})$ and $z\in \mathbb{R}^m$, is called a \textit{projection of polytopes}, provided we have $\pi(P)=Q$.
		\end{defn}
		\begin{lem}\label{lem7.10}
			Let $\pi: P \longrightarrow Q$ be a projection of polytopes, and $F$ be a face of $Q$. Then the preimage of $F$, $\pi^{-1}(F)$, is a face of $P$.
		\end{lem}
		\begin{proof}
			The proof can be found in \cite{Ziegler}, Lemma 7.10.
		\end{proof}
		\begin{shortversion}
			\begin{defn}\label{int}
				Let $F$ be a face of a polyhedron (or a polytope $P$), and $\dim(F)=d$ for some $d \leq \dim(P)$. The \textit{interior} of $F$, denoted by $int(F)$, is the set of all the points on $F$ that are not on any face of dimension less than $d$. It can be shown that this definition coincide with the topological interior of $F$ as a set of points.
			\end{defn}
		\end{shortversion}
		\begin{defn}\label{cplx}
			A finite collection $\Sigma$ of polyhedra is called a \textit{polyhedral complex}, if it satisfies two following conditions:
			\begin{enumerate}
				\item For each $P \in \Sigma$, any face $F$ of $P$ is also an element of $\Sigma$;
				\item for any two elements $P, P^{\prime}\in \Sigma$, $P \cap P^{\prime}$ is a face of both $P$ and $P^{\prime}$.
			\end{enumerate}
			Each polyhedron appearing in a polyhedral complex $\Sigma$ is called a \textit{cell} of $\Sigma$.
			Note that, $\emptyset$ is also a cell of $\Sigma$.\\
			The \textit{dimension of a polyhedral complex} is defined to be the maximum of the dimensions of its cells.\\
			A cell of a polyhedral complex $\Sigma$ is called a \textit{facet} of $\Sigma$, if it is not a face of another cell with higher dimension.  
			We say $\Sigma$ is \textit{pure}, if all its facets are of the same dimension.
		\end{defn}
		\begin{defn}\label{Gamma}
			\sloppy A polyhedron $P=\left\lbrace x \in \mathbb{R}^n\,\, | \,\, \mathcal{A}\cdot x \leq \mathcal{B} \right\rbrace $ is called a $\left(\Gamma,\mathbb{Q}(\rho)\right)-$\textit{polyhedron}, with $\rho$ being a positive transcendental real number, if $\mathcal{A} \in M_{d \times n}(\mathbb{Q}(\rho))$ and $\mathcal{B} \in \Gamma^d $.\\
			$\Sigma$ is called a $\left(\Gamma,\mathbb{Q}(\rho)\right)-$\textit{polyhedral complex}, if each of its cells is a $\left(\Gamma,\mathbb{Q}(\rho)\right)-$polyhedron. 
		\end{defn}
		
		\begin{defn}\label{supp}
			For a polyhedral complex $\Sigma \subset \mathbb{R}^n$, its \textit{support} is denoted by  $\lvert \Sigma\rvert$, and is defined as follows:
			\begin{equation*}
				\lvert \Sigma\rvert= \left \lbrace x \in \mathbb{R}^n\,\, | \,\, \text{There is a cell in}\,\, \Sigma\,\, \text{which contains}\,\,\, x \right \rbrace.
			\end{equation*}
		\end{defn}
		\begin{defn}\label{sbd}
			Suppose $r$ vectors $u_1, \dots , u_r$ in $\mathbb{R}^n$ are given. A vector  $w=(w_1,\dots,w_r)\in \mathbb{R}^r$ is called a \textit{weight vector}.
			We consider the two following polytopes:
			\begin{align}
				&P=conv\left\lbrace u_i\,\, |\,\, 1\leq i\leq r \right\rbrace\subseteq\mathbb{R}^n,\\
				&P_w=conv\left\lbrace\left(u_i,w_i\right)\,\,|\,\,1\leq i\leq r\right\rbrace\subseteq\mathbb{R}^{n+1}.
			\end{align}
			We say $F$ is a \textit{lower face} of $P_w$ if $\mathrm{face_v(P_w)}=F$ for some
			$v=\left(v_1,\dots,v_{n+1}\right) \in \left(\mathbb{R}^{n+1}\right)^{\vee}$, with $v_{n+1}$ positive.
			\par Suppose $\pi:\mathbb{R}^{n+1} \rightarrow \mathbb{R}^n$ is the projection onto the first $n$ coordinates. Using this map, if we project all lower faces of $P_w$ to $\mathbb{R}^n$,
			we obtain the following polyhedral complex:
			\begin{equation*}
				\Sigma_w=\left\lbrace \pi(F)\,\,|\,\,F\, \text{is a lower face of}\, P_w\right\rbrace,
			\end{equation*}
			for which we have $\lvert \Sigma_w \rvert=P$. This polyhedral complex gives a subdivision of $P$. This is called a \textit{regular subdivision} of $u_1, \dots, u_r$(or a \textit{regular subdivision} of $P$) induced by $w$.
		\end{defn}
		\begin{rem}
			The projection map defined in Definition \ref{sbd} is a projection of polytopes from the face $F$ to $\pi(F)$.
		\end{rem}
		\begin{shortversion}
			\begin{rem}\label{compact}
				Let $v$ be a vector in $\mathbb{R}^n$. Then for the vector $(v,1) \in \mathbb{R}^{n+1}$, we can find a lower face $F$ in $P_w$ such that $\mathrm{face_{(v,1)}(P_w)}=F$.
			\end{rem}
			
			\begin{proof}
				Assume $\phi_v:P_w \rightarrow \mathbb{R}$ is defined by $\phi_v(x)=(v,1)\cdot x$. Clearly $\phi_v$ is a continuous map. As $P_w$ is closed and bounded, it is compact. Hence, $\phi_v$ obtains its minimum on $P_w$. In other words, there exists a nonempty subset $F$ of $P_w$ such that $\forall x \in F,\,\,\,\, (v,1)\cdot x$ obtains the minimum. This means
				\begin{center}
					$F=\left\lbrace x\in P_w\,\,| \,\, (v,1)\cdot x \leq (v,1)\cdot y \,\,\, \forall y \in P_w \right\rbrace $.
				\end{center}  
				By definition, $F$ is a face determined by $(v,1)$ which is a vector with positive last coordinate, so $F$ is a lower face.
			\end{proof}
		\end{shortversion}
		\begin{defn}\label{n-cone}
			Let $P \subseteq \mathbb{R}^n$ be a polytope. To each face $F$ of $P$, we associate the cone $\mathcal{N}_P(F)$ defined as follows:
			\begin{center}
				$\mathcal{N}_P(F)=\left\lbrace v \in (\mathbb{R}^n)^\vee \,\,|\,\,\mathrm{face_v(P)}\supseteq F   \right\rbrace $.
			\end{center}
			It is called the \textit{closed normal cone} associated to $F$.
			\par The \textit{open normal cone} associated to $F$ is also defined as 
			\begin{center}
				$\mathring{\mathcal{N}}_P(F)=\left\lbrace v \in (\mathbb{R}^n)^\vee\,\, |\,\,\mathrm{face_v(P)}=F   \right\rbrace$.    
			\end{center}
			Obviously, we have $\mathring{\mathcal{N}}_P(F) \subseteq \mathcal{N}_P(F)$.
		\end{defn}
		
		\begin{defn}\label{n-fan}
			Using the previous definition, we define the \textit{normal fan} of a polytope $P$, which is
			\begin{center}
				$\mathcal{N}_P=\left\lbrace \mathcal{N}_P(F)\,\,|\,\, F \textrm{ is a face of } P \right\rbrace $.
			\end{center}
		\end{defn}
		Note that $\mathcal{N}_P$ is a polyhedral complex.
		\begin{defn}\label{dual}
			Assume $P\subseteq \mathbb{R}^n$, $P_w\subseteq \mathbb{R}^{n+1}$ and $\pi:\mathbb{R}^{n+1} \rightarrow \mathbb{R}^n$ are as defined in Definition \ref{sbd}. We write $\tilde{\pi}$ for the restriction of $\pi$ to the following set:
			\begin{equation*}
				\lbrace (v,1)\,\,|\,\,(v,1)\in \mathbb{R}^{n+1}\rbrace.   
			\end{equation*}
			We associate the following set to each lower face $F$ of $P_w$:
			\begin{equation*}
				\Tilde{\pi}\left(\mathring{\mathcal{N}}_P(F)\right)=\left\lbrace v\in \mathbb{R}^n\,\,|\,\,(v,1)\in \mathring{\mathcal{N}}_P(F) \right\rbrace.
			\end{equation*}
			The collection of all such sets form a polyhedral complex which is called the \textit{dual complex}, and is denoted by $\Sigma_w^{\circ}$. In other words, we have 
			\begin{equation*}
				\Sigma_w^{\circ}=\left\lbrace \Tilde{\pi}\left(\mathring{\mathcal{N}}_P(F)\right)\,\,|\,\, F\, \text{is a lower face of}\,P_w \right\rbrace.
			\end{equation*}
		\end{defn}
		\begin{rem}
			Given the assumptions of Definition \ref{dual}, the dual complex $\Sigma_{w}^{\circ}$ can intuitively be understood as the complex obtained by intersecting the normal fan of $P_w$ with the hyperplane defined by fixing the last coordinate to be one.
		\end{rem}
		\begin{shortversion}
			textcolor{red}{Maybe it's a good idea to put a picture here to clarify the intuition of all the last three definitions.}
		\end{shortversion}
		\begin{defn}\label{skl}
			Let $\Sigma$ be a polyhedral complex of dimension $n$. For $k<n$,\\ the $k$-\textit{skeleton} is the polyhedral complex consisting of all cells of $\Sigma$ with dimension at most $k$.
		\end{defn}
		\begin{defn}\textbf{(Some terminology about posets)}
			\begin{itemize}
				\item By a \textit{poset} $S$, we mean a partially ordered set. 
				\item A totally ordered subset of a poset $S$ is called a \textit{chain} in $S$. If a chain has $n$ elements, its \textit{length} is defined to be $n-1$. For a chain with infinitely many elements, the length is defined to be $\infty$.
				\item A poset $S$ is called \textit{bounded}, if it has a unique minimal element which is usually denoted by $\hat{0}$ and a unique maximal element which is usually denoted by $\hat{1}$.
				\item A finite poset is called \textit{graded}, if it is bounded and all its maximal chains are of the same length.
				\item For a graded poset $S$ and an element $a$ of $S$, the \textit{rank} of $a$ is defined to be the length of the maximal chain starting with $\hat{0}$ and ending with $a$. It is denoted by $r(a)$.
				\item A \textit{lattice} is a bounded poset $S$ such that for any two elements $a,b\in S$, there is an infimum ( denoted by $a \wedge b$) and also a supremum ( denoted by $a \vee b$) in $S$.
			\end{itemize}
		\end{defn}
		\begin{defn}
			Let $P\subseteq \mathbb{R}^n$ be a polytope. The \textit{face lattice} of $P$, denoted by $L(P)$, is the poset of all faces of $P$ where the partial order is inclusion.
		\end{defn}
		\begin{shortversion}
			\begin{thm}\label{2.7}
				If $P$ is a polytope, then its face lattice, L(P), is a graded lattice of length $\dim(P)+1$. Moreover, if $F$ is a face of $P$, then we have $r(F)=\dim(F)+1$.
			\end{thm}
			\begin{proof}
				The proof can be found in \cite{Ziegler}, Theorem 2.7. 
			\end{proof}
		\end{shortversion}
		\begin{defn}\label{fpo}
			Let $\Sigma$ be a polyhedral complex. Consider the poset of all its cells, where the partial order is inclusion and denote it by $L(\Sigma)$. Then $(L(\Sigma) , \subseteq)$ is called the \textit{cell poset} of $\Sigma$.
			\par Note that, this poset does not necessarily have a unique maximal element; in this case, it is not a bounded poset.
			To turn such cell posets into a bounded poset, we add an artificial maximal element, $\hat{1}$. Define $\hat{L}(\Sigma)=L(\Sigma)\cup \hat{1}$, then $\left(\hat{L}(\Sigma), \subseteq\right)$ is a bounded poset. Thus, for a pure polyhedral complex, $\left(\hat{L}(\Sigma), \subseteq\right)$ is a graded poset.
		\end{defn}
		\begin{rem}
			Let $\Sigma$ be a pure polyhedral complex, then for any cell $F$ of the graded poset $\left(\hat{L}(\Sigma),\subseteq\right)$, the rank of $F$ is defined, and we have $r(F)=\dim(F)+1$.
		\end{rem}
		Below, we recall two lemmas which are well known in Polyhedral Geometry.
		\begin{lem}\label{perp}
			Let $P\subseteq \mathbb{R}^n$ be a full dimensional polytope (of dimension $n$) and $F$ be one of its faces. For $\mathcal{N}_P(F)$, we define $\dim \left(\mathcal{N}_P(F) \right):= \dim\left( \mathrm{aff}\left( \mathcal{N}_P(F) \right)  \right)  $. Then we have $\dim(F)+ \dim\left(\mathcal{N}_P(F) \right)=n $.
		\end{lem}
		\begin{shortversion}
			\begin{proof}
				In this proof, by $L(P)$ we mean the face lattice of $P$ excluding $\emptyset$ . Define the following map
				\begin{align*}
					\phi:L(P)\longrightarrow \mathcal{N}_P, \\
					\,\,\,\,\,F\longmapsto \mathcal{N}_P(F).
				\end{align*}
				Firstly, we prove that $\phi$ is an anti-isomorphism meaning that, it is one-to-one and for any two faces $F$ and $F^{\prime}$ in $L(P)$, if $F\subset F^{\prime}$ then $\mathcal{N}_P(F) \supset \mathcal{N}_P(F^{\prime})$.\\
				From the definition, $\mathcal{N}_P(F^{\prime})$ can be written as
				\begin{equation}\label{N_P}
					\mathcal{N}_P(F^{\prime})=\left\lbrace  v\,\, |\,\, \mathrm{face_v(P)} \supseteq F^{\prime}\right\rbrace=\left\lbrace v\,\, |\,\, \forall x'  \in  F' \,\,\, \forall y \in P \,\,\, v\cdot x'\leq v\cdot y \right\rbrace.
				\end{equation}
				Suppose $u \in \mathcal{N}_P(F^{\prime})$ and $F \subset F^{\prime}$, from \ref{N_P}, this means, for any $x \in F$ and for any $y \in P $ we have $u\cdot x \leq u\cdot y$. In other words, we have \\
				$u \in \left\lbrace v\,\, | \,\, \forall x \in F \,\,\, \forall y \in P\,\,\, v\cdot x \leq v\cdot y \right\rbrace$ and similar to \ref{N_P}, this means $ u\in \mathcal{N}_P(F)$. Hence, $\mathcal{N}_P(F) \supseteq \mathcal{N}_P(F^{\prime})$.
				\par Now, what remains is to prove that $\phi$ is one-to-one. 
				Assume for $F \neq F^{\prime}$, we have $\mathcal{N}_P(F)=\mathcal{N}_P(F')$. Take an element $v$ of $\mathring{\mathcal{N}}_P(F)$. From the definition, this means $F=\mathrm{face_v(P)}$. By assumption, $\mathcal{N}_P(F)=\mathcal{N}_P(F')$ and from \ref{n-cone} we have $\mathring{\mathcal{N}}_P(F) \subseteq \mathcal{N}_P(F)=\mathcal{N}_P(F')$.
				This implies $v \in \mathcal{N}_P(F')$ meaning that $F' \subseteq \mathrm{face_v(P)}=F$. Similarly, by taking an element $v'$ of  $\mathcal{N}(F')$, we can deduce $F \subseteq F'$. Hence, $F=F'$. Therefore, $\phi$ is one to one. 
				\par From Theorem \ref{2.7}, we know that for a given polytope, $P$, its face lattice, $L(P)$, is a graded lattice of length $\dim(P)+1$. This means that all the maximal chains in $L(P)$ have the same length, which is $\dim(P)+1$. As we have assumed that $L(P)$ is the face lattice excluding $\emptyset$, we have reduced the length of all maximal chains by one. Therefore, all maximal chains of $L(P)$ have length $\dim(P)$.\\
				Let $F_0 \subset F_1 \subset \dots \subset F_i \subset \dots \subset F_n =P$ be a maximal chain in $L(P)$. Let $F_i$ and $F_{i+1}$ be two consecutive elements of this chain. From Theorem \ref{2.7}, we know that $r(F_i)= \dim(F_i)+1$ and $r(F_{i+1})=\dim(F_{i+1})+1$. As this chain is assumed to be maximal, this means that there is no other face $\hat{F_i}$ such that $F_i \subset \hat{F_i} \subset F_{i+1}$. Therefore, we can deduce that $r(F_{i+1})=r(F_{i})+1$, which gives $\dim(F_{i+1})- \dim(F_i)=1$. In this proof, we supposed that the empty face is excluded from $L(P)$. Hence, $F_0$ is a face of dimension $0$, which is a vertex. In fact, this maximal chain of length $\dim(P)=n$, starts with a face of dimension $0$ and ends to a face of dimension $n$. This means, for any natural number $i\,\,\, 0\leq i \leq n$, we can find a face in this chain, namely $F_i$, which is of dimension $i$.\\
				Apply the map $\phi$ on this maximal chain. We obtain a chain
				\begin{center}
					$\mathcal{N}_P(F_0) \supset\mathcal{N}_P(F_1) \supset \dots \supset \mathcal{N}_P(F_i)\supset \dots \supset \mathcal{N}_P(F_n)$,
				\end{center}
				which is of length $\dim(P)$.\\
				For $\mathcal{N}_P(F_i)$ and $\mathcal{N}_P(F_{i+1})$ which are two consecutive elements of this chain, $ \mathcal{N}_P(F_i) \supset \mathcal{N}_P(F_{i+1})$ means $\mathcal{N}_P(F_{i+1})$ is a proper face of the cone $\mathcal{N}_P(F_i)$.
				Therefore, we have $\dim \mathcal{N}_P(F_i) > \dim \mathcal{N}_P(F_{i+1})$. As the length of this chain is $\dim(P)$, it implies $\dim (\mathcal{N}_P(F_i))= \dim(P)-i$. Hence, for any $i$, we have $\dim(F_i)+\dim(\mathcal{N}_P(F_i))= i+\dim(P)-i=\dim(P)=n$.
				\par Note that for $F$, an arbitrarily given face of $P$, $F$ is necessarily contained in such a maximal chain. Otherwise, as we have $F \subset P$, this is a part of a chain that is not contained in any maximal chain of the described form.\\
				Therefore, this chain is also maximal but not of the length $\dim(P)$, and we know this can not happen because $L(P)$ is graded. Finally, this means that for any face $F$ of $P$, we have $\dim(F)+ \dim(\mathcal{N}_P(F))=n$. 
			\end{proof}
		\end{shortversion}
		\begin{lem}\label{pure}
			Let $\Sigma$ be a polyhedral complex whose support is $\mathbb{R}^n$, then $\Sigma$ is pure of dimension $n$.
		\end{lem}
		\begin{shortversion}
			\begin{proof}
				Let $\mathcal{A}$ be the set of all faces of $\Sigma$ of dimension $n$. In other words, we have
				\begin{equation*}
					\mathcal{A}=\left\lbrace F\,\,|\,\, F\,\, \text{is a face in}\,\, \Sigma\,\, \text{such that}\,\, \dim(F)=n  \right\rbrace.
				\end{equation*}
				Firstly, we want to prove that $\underset{F\in \mathcal{A}}{\bigcup}int(F)$ is dense in $\mathbb{R}^n$, which means $\overline{\underset{F\in \mathcal{A}}{\bigcup}int(F)}=\mathbb{R}^n$. Equivalently, we want to prove that the topological interior of $\mathbb{R}^n \setminus \underset{F\in \mathcal{A}}{\bigcup}int(F)= \lvert \Sigma \rvert \setminus \underset{F\in \mathcal{A}}{\bigcup}int(F)$ is empty.\\
				If we define 
				\begin{equation*}
					\mathcal{B}=\left\lbrace F\,\,|\,\, F\,\,\text{is a face in}\,\,\Sigma\,\,\text{such that} \,\, \dim(F)<n\right\rbrace,
				\end{equation*}
				we have
				\begin{equation*}
					\lvert \Sigma \rvert \setminus \underset{F\in \mathcal{A}}{\bigcup}int(F)=\underset{F\in \mathcal{B}}{\bigcup}F.
				\end{equation*}
				As $\mathcal{B}$ consists of faces of dimension less than $n$, $\underset{F\in \mathcal{B}}{\bigcup}F$ can not contain a $n$-ball.Therefore the interior of $\underset{F\in \mathcal{B}}{\bigcup}F$ is empty.\\
					Hence, $\overline{\underset{F\in \mathcal{A}}{\bigcup}int(F)}=\mathbb{R}^n$. As faces are topologically closed sets, we have
					\begin{equation*}
						\mathbb{R}^n = \overline{\underset{F \in \mathcal{A}}{\bigcup}int(F)} \subseteq \underset{F \in \mathcal{A}}{\bigcup}\overline{int(F)}\subseteq \underset{F \in \mathcal{A}}{\bigcup}\overline{F}=\underset{F \in \mathcal{A}}{\bigcup}F \subseteq \mathbb{R}^n.
					\end{equation*}
					So we have 
					\begin{equation}\label{**}
						\underset{F \in \mathcal{A}}{\bigcup}F = \mathbb{R}^n
					\end{equation}
					Finally, we choose arbitrarily a facet $G$ of $\Sigma$, and we prove $\dim(G)=n$. To prove this, we assume the opposite, meaning that $\dim(G) < n$. Choose $x\in G$ such that $x\in int(G)$ or in other words, $x$ is not on a proper face of $G$. As $x$ is also a point in $\mathbb{R}^n$, \ref{**} implies that there is a face $F$ of dimension $n$ such that $x \in F$. Therefore, we have $G \cap F \neq \emptyset$. From the definition of a polyhedral complex, $G \cap F$ is a face of both $G$ and $F$. As $x$ was chosen from an improper face of $G$, we have $G \cap F = G$. This means, $G \subseteq F$, and this contradicts the assumption that $G$ is a facet. Hence, $\dim(G)=n$, and $\Sigma$ is pure of dimension $n$. 
				\end{proof}
			\end{shortversion}
			
			\section{Combinatorics of Difference Tropical Hypersurfaces}
			In this section, we describe the combinatorics of a difference tropical hypersurface.
			This is done in Proposition \ref{combi}. The proof is Similar to the classical case, as in Proposition $3.1.6$ of \cite{Diane}. All Polyhedral Geometry needed to understand this result is presented in Subsection \ref{plygeo}.
			\begin{prop}\label{combi}
				If $f=\underset{u(\sigma)}{\sum}c_{u(\sigma)}x^{u(\sigma)}$ is a Laurent difference polynomial, then its associated difference tropical hypersurface, $\mathoptrop(V(f))$, is the support of a pure $\left(\Gamma,\mathbb{Q}(\rho)\right)-$polyhedral complex of dimension $(n-1)$.\\
				More precisely, it is the $(n-1)$-skeleton of the polyhedral complex dual to the regular subdivision $\Sigma_{val}$ of $P$, where $P$ is defined as in Definition \ref{Newpolytope}, and $val$ is the weight vector given by $v\left( c_{u(\sigma)}\right) $ for $c_{u(\sigma)}\neq 0$.
				
			\end{prop}
			\begin{proof}
				We begin the proof by describing $P_{val}$, and the dual complex $ \Sigma_{val}^{\circ}$. In the first step, we prove Claim \ref{pi(F)} which allows us to show that $\mathoptrop(V(f))$ is the $(n-1)$-skeleton of the dual complex. Subsequently, we establish Claim \ref{m(v)}, and having this claim in hand, we demonstrate that $\mathoptrop(V(f))$ is a $\left(\Gamma,\mathbb{Q}(\rho)\right)-$polyhedral complex. Finally, we show that $ \lvert\Sigma_{val}^{\circ}\rvert=\mathbb{R}^n$ guarenteeing that $\mathoptrop(V(f))$ is pure.
				\par Based on Definition \ref{trophyp}, we know that $\mathoptrop(V(f))$ is the set of all points for which the minimum in $\mathoptrop(f)$ is attained at least twice.\\
				Considering the Newton polytope of $f$, we define
				\begin{center}
					$ P_{val}:=conv\left\lbrace \left( u(\rho), v(c_{u(\sigma)}) \right): c_{u(\sigma)} \neq 0  \right\rbrace \subset \mathbb{R}^{n+1} $. 
				\end{center} 
				Suppose $\pi :\mathbb{R}^{n+1}\longrightarrow \mathbb{R}^n$ is defined by $\pi\left(x_1,\dots,x_n,x_{n+1} \right)=(x_1,\dots,x_n) $.
				As we discussed in Definition \ref{sbd}, if we take
				$\pi(F)$ for all lower faces $F$ of $P_{val}$, we obtain the regular subdivision of $P$ induced by $v\left( c_{u(\sigma)}\right) $ for $ c_{u(\sigma)}\neq 0$.\\
				Taking $\tilde{\pi}\left(\mathring{\mathcal{N}}_P(F) \right) $ for all lower faces $F$ of $P_{val}$, we obtain the polyhedral complex dual to the regular subdivision of $P$ induced by $v\left( c_{u(\sigma)}\right) $ for $ c_{u(\sigma)}\neq 0$.\\
				Suppose $v=(v_1,\dots,v_n,1)\in \mathring{\mathcal{N}}_P(F)$. Then we have
				\begin{equation}	
					\mathopin_{\pi(v)}(f)=\underset{\substack{u(\sigma):v\left( c_{u(\sigma)}\right)+u(\rho)\cdot (v_1,\dots,v_n)\\=\mathoptrop(f)(v_1,\dots,v_n)}}{\sum}\overline{t^{v\left(c_{u(\sigma)}\right)}\cdot c_{u(\sigma)}} x^{u(\sigma)}.
				\end{equation}
				Note that, we can write
				\begin{center}
					$v\left( c_{u(\sigma)}\right)+u(\rho)\cdot (v_1,\dots,v_n)=1\cdot v\left( c_{u(\sigma)}\right)+u(\rho)\cdot (v_1,\dots,v_n)=v\cdot \left(u(\rho),v\left( c_{u(\sigma)}\right)\right) $.
				\end{center}
				So we have
				\begin{equation}\label{pi(v)}
					\mathopin_{\pi(v)}(f)=\underset{\substack{u(\sigma):v\cdot \left(u(\rho),v\left( c_{u(\sigma)}\right)\right)\\
							=\mathoptrop(f)(\pi(v))}}{\sum}\overline{t^{-v\left(c_{u(\sigma)} \right) }\cdot c_{u(\sigma)}} x^{u(\sigma)}.
				\end{equation}
				\begin{claim}\label{pi(F)}
					Let $u(\sigma)$ be an exponent appearing in $\mathopin_{\pi(v)}(f)$. Then $u(\rho)$ is in $\pi(F)$. Moreover, for each vertex $u(\rho)$ of $\pi(F)$, $u(\sigma)$ is an exponent appearing in $\mathopin_{\pi(v)}(f)$.     
				\end{claim}
				\renewcommand\qedsymbol{$\blacksquare$}
				\begin{proof}
					First, we suppose $u_0(\sigma)$  is an exponent appearing in $\mathopin_{\pi(v)}(f)$, then we prove that there is a vertex in $\pi(F)$ which corresponds to this exponent.\\
					As $u_0(\sigma)$ is an exponent appearing in $\mathopin_{\pi(v)}(f)$, this means
					\begin{center}
						$v\cdot \left(u_0(\rho), v\left(c_{u_0(\sigma)} \right)  \right)=\mathoptrop(f)(\pi(v)) $.   
					\end{center}
					In other words, $v\cdot\left(u_0(\rho), v\left(c_{u_0(\sigma)} \right)  \right)$ obtains the minimum among all tropical monomials of the form $v\left(c_{u(\sigma)} \right)+u(\rho)\cdot\pi(v) $.\\
					Assume $y \in P_{val} =conv\left\lbrace \left( u(\rho), v(c_{u(\sigma)}) \right): c_{u(\sigma)} \neq 0  \right\rbrace $, so $y$ can be written as:
					\begin{equation*}
						y=\underset{u(\sigma)}{\sum}\lambda_{u(\sigma)}\left(u(\rho), v\left(c_{u(\sigma)} \right)  \right),
					\end{equation*}
					such that $\forall u(\sigma), \,\,\,\,\ 0\leq\lambda_{u(\sigma)}\leq 1$ and $\underset{u(\sigma)}{\mathlarger{\sum}}\lambda_{u(\sigma)}=1$.\\
					We also have
					\begin{align}
						v\cdot y&=\underset{u(\sigma)}{\sum}\lambda_{u(\sigma)} v\cdot \left(u(\rho), v\left(c_{u(\sigma)} \right) \right)\\
						&=\lambda_{u_0(\sigma)} v\cdot \left(u_0(\rho), v\left(c_{u_0(\sigma)} \right)  \right)+\underset{u(\sigma)\neq u_0(\sigma)}{\sum}\lambda_{u(\sigma)} v\cdot\left(u(\rho), v\left(c_{u(\sigma)} \right)\right).\label{v.y} 
					\end{align}
					Since $u_0(\sigma)$ is one of the exponents appearing in $\mathopin_{\pi(v)}(f)$, we have
					\begin{equation*}
						v\cdot\left(u_0(\rho),v\left(c_{u_0(\sigma)} \right)  \right)\leq v\cdot\left(u(\rho), v(c_{u(\sigma)})\right),   
					\end{equation*}
					so \eqref{v.y} can be written as:
					\begin{equation*}
						\lambda_{u_0(\sigma)} v\cdot \left(u_0(\rho), v\left(c_{u_0(\sigma)} \right)  \right)+\underset{u(\sigma)\neq u_0(\sigma)}{\sum}\lambda_{u(\sigma)}\left( v\cdot\left(u_0(\rho), v\left(c_{u_0(\sigma)} \right)\right)+ \alpha_{u(\sigma)}\right),
					\end{equation*}
					where 
					\begin{equation*}
						\alpha_{u(\sigma)}=
						v\cdot\left(u(\rho), v(c_{u(\sigma)})\right)-v\cdot\left(u_0(\rho),v\left(c_{u_0(\sigma)} \right)  \right)\geq 0.
					\end{equation*}
					As
					\begin{equation*}
						\underset{u(\sigma)}{\sum}\lambda_{u(\sigma)}=1,
					\end{equation*}
					we have
					\begin{equation*}
						v\cdot y= v\cdot\left(u_0(\rho),v\left(c_{u_0(\sigma)} \right)  \right)+ \underset{u(\sigma)\neq u_0(\sigma)}{\sum}\lambda_{u(\sigma)}\alpha_{u(\sigma)}.
					\end{equation*}
					On the other hand, for all $u(\sigma)$, we have $\lambda_{u(\sigma)} \geq 0$.
					Hence, this gives
					\begin{equation}\label{u_0}
						v\cdot y \geq v\cdot\left(u_0(\rho),v\left(c_{u_0(\sigma)} \right)  \right).
					\end{equation}
					As $\left(u_0(\rho),v\left(c_{u_0(\sigma)} \right)  \right) \in P_{val} $ and $v \in \mathring{\mathcal{N}}_P(F) $, \eqref{u_0} implies
					$\left(u_0(\rho),v\left(c_{u_0(\sigma)} \right)  \right) \in F$. Therefore $u_0(\rho) \in \pi(F)$.\\
					\par Now, assume a vertex of $\pi(F)$ is given. We prove that there is an exponent in $\mathopin_{\pi(v)}(f)$ associated to this vertex.
					\par From Definition \ref{sbd}, a vertex of $\pi(F)$ is of the form $u(\rho)$. Consider the restriction of the projection map on the face $F$. Hence,  $\pi\big|_{F}: F \longrightarrow \pi(F)$ is a projection of polytopes.
					Therefore, by Lemma \ref{lem7.10}, the preimage of $u(\rho)$ is a face of $F$; in fact it is a vertex of $F$.\\
					To see this, assume $\pi\big|_{F}^{-1}(u(\rho))$ is a face $G$ such that $\dim(G)\geq 1$. As $G$ is a face of $F$ and $F$ is a face of $P_{val}$, from Proposition \ref{prop2.3}, we have $G=conv(\mu)$ where $\mu= G \cap U$ for $U$ being the set of vertices of $P_{val}$. As $\dim(G) \geq 1$, it contains at least two vertices, each of which is a vertex of $P_{val}$. Suppose $\left(u_i(\rho), v\left(c_{u_i(\sigma)}\right) \right)$ and $\left(u_j(\rho), v\left(c_{u_j(\sigma)}\right) \right) $ are two distinct vertices of $G$, so we have
					\begin{equation*}
						u(\rho)=\pi\big|_F\left(u_i(\rho), v\left(c_{u_i(\sigma)}\right) \right)=u_i(\rho)=\pi\big|_F\left(u_j(\rho), v\left(c_{u_j(\sigma)}\right)\right)=u_j(\rho),
					\end{equation*}
					and $u_i(\rho)=u_j(\rho)$ contradicts the assumption that $\rho$ is transcendental. Hence, $\dim(G)=0$ which means $G$ is a vertex of $F$ and by definition we have $G=\left(u(\rho), v(c_{u(\sigma)}) \right)$.\\
					Since $v=(v_1,\dots,v_n,1) \in \mathring{\mathcal{N}}_P(F)$, for all $y \in P_{val}$, we have $v\cdot y \geq v\cdot\left(u(\rho), v\left(c_{u(\sigma)} \right) \right)$. More specifically, for each vertex $\left(u^{\prime}(\rho), v\left(c_{u^{\prime}(\sigma)} \right) \right) \in P_{val} $, we have
					\begin{align*}
						&v\cdot \left(u(\rho), v\left(c_{u(\sigma)} \right) \right) \leq v\cdot \left(u^{\prime}(\rho), v\left(c_{u^{\prime}(\sigma)} \right) \right)\\
						\iff &\pi(v)\cdot u(\rho)+v\left(c_{u(\sigma)} \right) \leq \pi(v)\cdot u^{\prime}(\rho)+v\left(c_{u^{\prime}(\sigma)} \right).\\
					\end{align*}
					This last inequality means, in $\mathoptrop(f)(\pi(v))$ the minimum is achieved at $u(\sigma)$ and equivalently 
					$u(\sigma)$ appears as an exponent in $\mathopin_{\pi(v)}(f)$.
				\end{proof}
				\renewcommand\qedsymbol{$\square$}
				Now, by using Claim \ref{pi(F)}, we want to prove that $\mathoptrop(V(f))$ is the $(n-1)-$skeleton of the dual complex $ \Sigma_{val}^{\circ}$.
				Let $F$ be a nonvertex lower face of $P_{val}$, and $w\in \tilde{\pi}\left(\mathring{\mathcal{N}}_P(F) \right) $. As $F$ is a nonvertex face, $\pi(F)$ has more than one vertex. From Claim \ref{pi(F)}, we know each of these vertices corresponds to an exponent in $\mathopin_w (f)$. This means $\mathopin_w (f)$ is not a monomial, or in other words from Lemma \ref{roots} , we have $w \in \mathoptrop(V(f))$.\\
				Now, choose $w \in \mathoptrop(V(f))$.
				\begin{shortversion}
					By Remark \ref{compact},     
				\end{shortversion}
				Set $F=\mathrm{face_{(w,1)}(P_{val})}$. This means that $w \in \Tilde{\pi}\left(\mathring{\mathcal{N}}_P(F)\right)$. As we assumed, $w \in \mathoptrop(V(f))$, equivalently from Lemma \ref{roots}, this means $\mathopin_w(f)$ is not a monomial. So, at least two different exponents appear in $\mathopin_w(f)$, each of which corresponds to a point on $\pi(F)$. Hence, $\pi(F)$ has more that one vertex, and therefore $F$ is a face with more than one vertex.\\
				Finally, this gives, $w \in \mathoptrop(V(f))$ if and only if $w \in \Tilde{\pi}\left(\mathring{\mathcal{N}}_P(F)\right)$ where F is a face with more than one vertex. As $F$ has at least two vertices, we have $\dim(F)\geq 1$, and by Lemma \ref{perp}, this means $\dim\left(\mathcal{N}_p(F)\right) \leq n-1$. We know that $\mathring{\mathcal{N}}_P(F) \subseteq \mathcal{N}_p(F)$. This implies that $\dim\left(\mathring{\mathcal{N}}_P(F)\right)\leq \dim \left(\mathcal{N}_p(F)\right) \leq n-1$, meaning that both $\mathring{\mathcal{N}}_P(F)$ and $\Tilde{\pi}\left(\mathring{\mathcal{N}}_P(F)\right)$ are not full dimensional. Hence, $w \in \mathoptrop(V(f))$ if and only if $\Tilde{\pi}\left(\mathring{\mathcal{N}}_P(F)\right)$ is a cell in the $(n-1)-$skeleton of the dual complex which contains $w$.\\
				Thus, $\mathoptrop(V(f))$ is the $(n-1)-$skeleton of the dual complex. What remains is to prove that it is a pure $\left(\Gamma,\mathbb{Q}(\rho)\right)-$polyhedral complex.
				\par Suppose $v \in \mathbb{R}^{n+1}$ is given arbitrarily. We define 
				\begin{center}
					$m(v):=\inf\left\lbrace v\cdot y\,\,|\,\,y \in P_{val}  \right\rbrace  $.
				\end{center}
				Now for $m(v)$, we prove the following claim.
				\begin{claim}\label{m(v)}
					Let $V_{P_{val} }$ be the set of vertices of $P_{val}$. Then for $m(v)$ which is defined above, we have
					\begin{equation*}
						m(v)= \min\left\lbrace v\cdot u\,\,|\,\, u \in V_{P_{val} } \right\rbrace.
					\end{equation*}
				\end{claim}
				\renewcommand\qedsymbol{$\blacksquare$}
				\begin{proof}
					Put $m:= \min \left \lbrace v \cdot u \,\, |\,\, u \in V_{P_{val}} \right \rbrace$. Since $V_{P_{val}}$ is included in $P_{val}$, we have $m(v) \leq m$.\\
					Now, we prove $m \leq m(v)$.
					Given $y \in P_{val}$ arbitrarily. From the definition of $P_{val}$, we have
					\begin{equation*}
						y=\underset{u\in V_{P_{val}}}{\sum}\lambda_{u}\cdot u \,\,\,\text{where}\,\,\,\forall u\,\,\, 0\leq \lambda_{u} \leq 1\,\,\, \text{and} \underset{u\in V_{P_{val}}}{\sum}\lambda_{u}=1.
					\end{equation*}
					This gives
					\begin{equation*}
						v\cdot y=v\cdot \underset{u\in V_{P_{val}}}{\sum}\lambda_{u}\cdot u=
						\underset{u\in V_{P_{val}}}
						{\sum}\lambda_{u}(v\cdot u)\geq \underset{u\in V_{P_{val}}}{\sum}\lambda_{u}m=m.
					\end{equation*}
					Thus, $\forall y \in P_{val} \,\,\,v\cdot y \geq m$, which means $m(v) \geq m$. Hence, $m(v)=m$.    
				\end{proof}
				\renewcommand\qedsymbol{$\square$}
				By using Claim \ref{m(v)}, we prove that, this $(n-1)-$skeleton is a $\left(\Gamma,\mathbb{Q}(\rho)\right)-$polyhedral complex.
				\par Let $F$ be a lower face of $P_{val}$. For $v \in \mathring{\mathcal{N}}_P(F)$, we have
				\begin{equation}\label{m}
					\forall x\in F\,\,\, v\cdot x \leq v\cdot y\,\,\, \forall y\in P_{val}\Rightarrow \forall x\in F\,\,\, v\cdot x=m(v)=m .
				\end{equation}
				If we denote $F \cap V_{P_{val}}$ by $V_F$ \eqref{m} means, for any vertex $x\in V_F $, we have $v\cdot x=m$. In other words, $\mathring{\mathcal{N}}_P(F)$ can be written as:
				\begin{equation*}
					\mathring{\mathcal{N}}_P(F)=\left\lbrace  v\,\, |\,\,\forall x \in V_F \,\,\,\,\, \forall y \in V_{P_{val}}\,\,\, v\cdot x \leq v\cdot y \right\rbrace. 
				\end{equation*} 
				Therefore, $\Tilde{\pi}\left(\mathring{\mathcal{N}}_P(F) \right) $, which is a cell in the  $(n-1)-$skeleton of the dual complex, can be written as:
				\begin{equation}\label{pitilde}
					\begin{aligned}
						\tilde{\pi}\left(\mathring{\mathcal{N}}_P(F) \right)&=\left\lbrace w \in \mathbb{R}^n\,\, |\,\, (w,1)\in \mathring{\mathcal{N}}_P(F) \right\rbrace\\
						&=\left\lbrace w \in \mathbb{R}^n\,\, |\,\, \forall x \in V_F\,\,\,\,\, \forall y \in V_{P_{val}}\,\,\, (w,1)\cdot x \leq (w,1)\cdot y  \right\rbrace.   
					\end{aligned}
				\end{equation}
				Note that in \eqref{pitilde}, $x$ and $y$ are both vertices of $P_{val} $. Therefore, for some $i$ and $j$, they are of the following form:
				\begin{equation}\label{x,y}
					\begin{aligned}
						&x=\left(u_i (\rho), v\left(c_{u_i(\sigma)} \right)  \right), \\
						&y=\left(u_j (\rho), v\left(c_{u_j(\sigma)} \right)  \right). 
					\end{aligned}
				\end{equation}
				Using \eqref{x,y}, we can rewrite the inequality that appeared in \eqref{pitilde} as follows:
				\begin{align*}
					(w,1)\cdot x \leq (w,1)\cdot y &\Longleftrightarrow w\cdot u_i(\rho)+v\left(c_{u_i(\sigma)} \right) \leq w\cdot u_j(\rho)+v\left(c_{u_j(\sigma)} \right)\\
					&\Longleftrightarrow v\left(c_{u_i(\sigma)} \right)-v\left(c_{u_j(\sigma)} \right) \leq w\cdot\left(u_j(\rho)-u_i(\rho)\right).  
				\end{align*}
				As $v\left(c_{u_i(\sigma)} \right)-v\left(c_{u_j(\sigma)} \right) \in \Gamma$ and $\left(u_j(\rho)-u_i(\rho)\right) \in {\mathbb{Q}(\rho)}^n$, we conclude $\tilde{\pi}\left(\mathring{\mathcal{N}}_P(F) \right) $ is a $\left(\Gamma,\mathbb{Q}(\rho)\right)-$polyhedron. Therefore, the $(n-1)-$skeleton of the dual complex is a $\left(\Gamma,\mathbb{Q}(\rho)\right)-$polyhedral complex. To complete the proof, we show that it is also pure.
				\par To do so, arbitrarily choose a facet $F$ of this $(n-1)-$skeleton; we prove $F$ is of dimension $n-1$. Assume the opposite. In fact, we assume $\dim (F)=m < n-1$.\\
				Suppose $v \in \mathbb{R}^n$ is given. 
				\begin{shortversion}
					From Remark \ref{compact},    
				\end{shortversion}
				Set $G=\mathrm{face_{(v,1)}(P_{val})}$. In other words $v \in \tilde{\pi}\left(\mathring{\mathcal{N}}_p(G) \right) $ which is a cell of the dual complex $ \Sigma_{val}^{\circ}$. Therefore, $v \in \lvert \Sigma_{val}^{\circ} \rvert \subseteq \mathbb{R}^n$. 
				Hence, $\lvert \Sigma_{val}^{\circ} \rvert= \mathbb{R}^n$. From Lemma \ref{pure}, $\Sigma_{val}^{\circ}$ is pure of dimension $n$, meaning that, there exists a facet $F^{\prime}$ in $\Sigma_{val}^{\circ}$ that is of dimension $n$ and contains $F$. Comparing $r(F)$ and $r(F^{\prime})$, we can find a cell $F^{\prime\prime}$ of the $(n-1)$-skeleton which contains $F$. But this is a contradiction, because $F$ is assumed to be a facet. Hence, $dim(F)=n-1$, and the $(n-1)$-skeleton is pure. 
				Thus, $\mathoptrop\left(V(f) \right) $ is a pure $\left(\Gamma,\mathbb{Q}(\rho)\right)$-polyhedral complex.\\
				\begin{shortversion}
					Consider the face poset of $\Sigma_{val}^{\circ}$, and the rank function on the elements of $L(\Sigma_{val}^{\circ})$. We have
					\begin{center}
						$r(F^{\prime})=\dim(F^{\prime})+1=n+1$ and $r(F)=\dim(F)+1=m+1< n$.
					\end{center}
					Suppose $\emptyset \subset \dots \subset F \subset \dots \subset F^{\prime}$ is a maximal chain. Note that, $r(F)=m+1<n+1=r(F^{\prime})$. Since we have $m<n-1$, there exists a face $F^{\prime\prime}$ in this chain such that $m+1 < r(F^{\prime\prime})< n+1$. This gives $\dim(F^{\prime\prime})< n$. Thus, $F^{\prime\prime}$ is a cell in the $(n-1)$-skeleton which contains $F$. But this is a contradiction, because $F$ is assumed to be a facet in the $(n-1)$-skeleton. Hence, $dim(F)=n-1$, and the $(n-1)$-skeleton is pure. 
					Thus, $\mathoptrop\left(V(f) \right) $ is a pure $\left(\Gamma,\mathbb{Q}(\rho)\right)$-polyhedral complex.
					
				\end{shortversion}
			\end{proof}
			\section{Difference Newton Lemma}
			This section provides an essential tool for proving the Difference Kapranov Theorem: the Difference Newton Lemma. For the notation consult \ref{sigmabold}. It states 
			\begin{thm}\label{DiffNewton}
				(Difference Newton Lemma)
				Let $K$ be a multiplicative valued difference field of characteristic zero that is spherically complete.
				Assume $\mathbf{k}$, the difference residue field of $K$, is an ACFA and the difference value group $\Gamma$ of $K$ is a $\mathbb{Q}(\rho)$-module where $\rho$, the scaling exponent of $\sigma$, is transcendental. Given $f\in K_{\sigma}[x]$ is not constant and suppose $b\in K$ such that $f(b)\neq 0$.\\
				We define $\varepsilon := \underset{ \underset{ |J|\geq 1}{J}}{\max}\,\,\varepsilon _J$, where
				\begin{equation*}
					\varepsilon_J :=\frac{1}{|J|_{\rho}}\left(v(f(b)) - v(f_{(J)}(b)) \right).   
				\end{equation*}
				There exists a root $a \in K $ of $f$ such that $v(a-b)=\varepsilon$.
			\end{thm}
			From this section on, we keep these assumptions on $K$.
			\par To prove the Difference Newton Lemma, firstly, we need some definitions. The main reference for these definitions is \cite{Van den Dries}.
			\begin{defn}
				Let $K$ be a valued field, and $\left(a_{\rho}\right)$ be a sequence of elements in $K$. The sequence $\left(a_{\rho}\right)$ is called \textit{well-indexed}, if the set of indices is well-ordered without maximal element.
				\par A well-indexed sequence $\left(a_{\rho}\right)$ is called \textit{pseudoconvergent} to a point $a$ if and only if
				\begin{equation*}
					\exists \rho_0\,\,\,\,\, \forall\delta,\rho:\,\,\, \delta>\rho>\rho_0 \longrightarrow v(a-a_{\delta})>v(a-a_{\rho}).
				\end{equation*}
				In this case, we use the notation $a_{\rho}\rightsquigarrow a$, and $a$ is called a \textit{pseudolimit} of this sequence.
			\end{defn}
			\begin{defn}
				A well-indexed sequence is called \textit{pseudocauchy} (abbreviated as pc-sequence) if and only if
				\begin{equation*}
					\exists \rho_0\,\,\,\,\, \forall\tau,\delta,\rho:\,\,\, \tau>\delta>\rho>\rho_0 \longrightarrow v(a_{\tau}-a_{\delta})>v(a_{\delta}-a_{\rho}).
				\end{equation*}
			\end{defn}
			\par Lemma \ref{l1} states that if a point $b$ is not a root of a difference polynomial $f$, we can find a point $a$ in its neighbourhood which is a better estimation of a possible root.  
			\begin{lem}\label{l1}
				Suppose
				$f\in K_{\sigma} [x]$
				is a nonconstant difference polynomial.
				Let
				$b$ be an element of $K$ which is not a root of
				$f$. We define 
				$\varepsilon := \underset{ \underset{ \lvert J\rvert \geq 1}{J}}{\max}\,\,\varepsilon _J $
				where 
				\begin{equation*}
					\varepsilon_J :=\frac{1}{\lvert J\rvert_{\rho}}\left(v(f(b)) - v(f_{(J)}(b)) \right).    
				\end{equation*}
				
				Then
				\begin{enumerate}
					\item
					there exists 
					$a \in K$ 
					such that $v(a-b) = \varepsilon$ and 
					$v(f(a))> v(f(b))$.
					\item
					for any point 
					$a\in K$
					with the properties in $(1)$, we have  
					$\varepsilon < \varepsilon'$
					where 
					$\varepsilon' := \underset{ \underset{ |J|\geq 1}{J}}{\max}\,\,\varepsilon' _J $
					and
					$\varepsilon'_J :=\frac{1}{|J|_{\rho}}\left(v(f(a)) - v(f_{(J)}(a)) \right) $.
				\end{enumerate}
			\end{lem}
			\begin{proof}
				\begin{enumerate}
					\item
					Note that $\varepsilon$ 
					is an element of the difference value group 
					$\Gamma$.
					Fix any
					$a'\in K$ 
					such that
					$v(a'-b)= \varepsilon$.\\
					Assume  
					$I$ 
					is a multi-index for which
					$\lvert I\rvert \geq 1$.\\ 
					Then we have
					\begin{equation}\label{non-neg}
						\begin{aligned}
							v\left(\frac{f_{(I)} (b) \boldsymbol{\sigma} ^I (a'-b)}{f(b)}\right) &= v(f_{(I)}(b))+v\left(\boldsymbol{\sigma} ^I (a'- b)\right)-v(f(b))\\
							&=v(f_{(I)} (b))+|I|_{\rho}v(a'-b)-v(f(b))\\
							&=-\lvert I \rvert _{\rho} \varepsilon_I +\lvert I \rvert _{\rho} \varepsilon\\
							&=\lvert I \rvert _{\rho} \left(\varepsilon - \varepsilon_I \right),
						\end{aligned}
					\end{equation}
					which is nonnegative. Obviously, if $J$ is a multi-index for which $\lvert J\rvert \geq 1$ and $\varepsilon_J = \varepsilon$, then we have
					\begin{equation*}
						v \left(\frac{f_{(J)}(b)\boldsymbol{\sigma} ^J (a'- b)}{f(b)}\right) = \lvert J \rvert_{\rho}\left(\varepsilon-\varepsilon_J\right)=0,   
					\end{equation*}
					which means
					$\frac{f_{(J)} (b) \boldsymbol{\sigma} ^J (a'- b)}{f(b)} \notin \mathcal{M}$, with $\mathcal{M}$ being the maximal ideal of the valuation ring. \\
					Hence, we can define the following nonconstant difference polynomial in $\mathbf{k}_{\sigma}[x]$:
					\begin{equation*}
						\phi (x) := 1+ \sum\limits_{\underset{|I|\geq 1}{I}} \overline{\left(\frac{f_{(I)} (b) \boldsymbol{\sigma} ^I (a'- b)}{f(b)}\right)}\boldsymbol{\bar{\sigma}}^{I}(x).
					\end{equation*}
					Since 
					$\mathbf{k}$ is an ACFA, and
					$\phi (x)$ is not monomial, by Lemma \ref{ACFAroot}, $\phi(x)$
					has a nonzero root 
					$\bar{u}$. This means that 
					if $u$ is a lift of $\bar{u}$ in the valuation ring then $v(u)=0$.\\
					Define $a:= (a'-b)u + b$.
					We have
					\begin{equation*}
						v(a-b) = v((a'-b)u) = v(a'-b) + v(u) = \varepsilon + 0 = \varepsilon. 
					\end{equation*}
					Consider the Taylor expansion of 
					$f(a)$
					around the point 
					$b$.
					It is
					\begin{equation*}
						f(a)= f(b) + \sum\limits_{\underset{|I|\geq 1}{I}} f_{(I)} (b) \boldsymbol{\sigma} ^I (a - b),    
					\end{equation*}
					from which we have
					\begin{equation}\label{div}
						\frac{f(a)}{f(b)} =  1 + \sum\limits_{\underset{|I|\geq 1}{I}}\frac{f_{(I)} (b) \boldsymbol{\sigma} ^I (a - b)}{f(b)}.
					\end{equation}Similar to what we did in \eqref{non-neg}, we can see that
					\begin{equation*}
						\forall I\,\,\,\, \lvert I \rvert \geq 1\,\,:\,\,v\left(\frac{f_{(I)} (b) \boldsymbol{\sigma} ^I (a-b)}{f(b)}\right) \geq 0. 
					\end{equation*}
					Taking the residue of both sides in \eqref{div}, 
					and substituting \\
					$a=(a'- b)u+ b$ gives
					\begin{equation*}
						\overline{\frac{f(a)}{f(b)}} = 1 + \sum\limits_{\underset{|I|\geq 1}{I}} \overline{\frac{f_{(I)} (b) \boldsymbol{\sigma} ^I (a - b)}{f(b)}} = 1 + \sum\limits_{\underset{|I|\geq 1}{I}} \overline{\frac{f_{(I)} (b) \boldsymbol{\sigma} ^I (a'- b)}{f(b)}}\bar{\boldsymbol{\sigma}}^I(\bar{u}) = \phi (\bar{u}) = 0.
					\end{equation*}
					This implies
					$v(f(a))> v(f(b))$.
					\item
					The proof of this part mainly consists of proof by contradiction. For this purpose, we need two technical steps which enable us to obtain the contradiction. 
					\begin{shortversion}
						Before discussing these steps,  We recall that the lexicographical order is defined as follows:
						\begin{equation*}
							I \geq_{lex} J \Leftrightarrow \exists k, \,\,\,\,\, 0\leq k \leq n\,\,\,\,\, \textrm{such that}\,\,\,\,\, i_k \geq j_k, \,\,\,\,\, \textrm{and}\,\,\,\,\, \forall l < k, \,\,\,\,\,
							i_l= j_l.
						\end{equation*}
					\end{shortversion}
					\begin{itemize}
						\item \textbf{Step 1}: We assume
						$\varepsilon = \varepsilon_I$
						where
						$I$ 
						is the maximal multi-index, with respect to lexicographical order, for which $\varepsilon_I$ attains the maximum.
						Then we prove
						$v(f_{(I)}(a))=v(f_{(I)}(b))$. 
						\item \textbf{Step 2}: If for some multi-index $J$, we have $\varepsilon=\varepsilon_J$, then we have 
						$\underset{I}{\min}\left(\lvert I \rvert_{\rho}\varepsilon + v(f_{(I)}(b))\right) = \lvert J\rvert_{\rho}\varepsilon + v(f_{(J)}(b))$.
					\end{itemize}
					\textbf{Step 1}: Let $I$ is as assumed above. Then $v(f_{(I)}(a))=v(f_{(I)}(b))$.
					\renewcommand\qedsymbol{$\blacksquare$}
					\begin{proof}
						
						For an arbitrary nonzero multi-index 
						$L \in \mathbb{N}^{n+1}$, we have $I<_{lex}I+L$. As we have assumed $I$ is the maximal multi-index for which $\varepsilon= \varepsilon_I$, we then have 
						$\varepsilon_{I+L} < \varepsilon _I = \varepsilon $. This means
						\begin{align}
							\begin{split}
								&\frac{1}{|I+L|_{\rho}}\left( v(f(b)) - v(f_{(I+L)}(b)\right)< \frac{1}{|I|_{\rho}}\left(v(f(b)) - v(f_{(I)}(b))\right)\\
								\Leftrightarrow& v(f(b)) - v(f_{(I+L)}(b) < \frac{|I+L|_{\rho}}{|I|_{\rho}}\left(v(f(b)) - v(f_{(I)}(b))\right)\\
								\Leftrightarrow&  - v(f_{(I+L)}(b))< \left(1+ \frac{|L|_{\rho}}{|I|_{\rho}}\right)\left(v(f(b)) - v(f_{(I)}(b))\right) - v(f(b))\\
								\Leftrightarrow&  - v(f_{(I+L)}(b))< - v(f_{(I)}(b)) + |L|_{\rho}\varepsilon_I
							\end{split}\\
							\begin{split}\label{I max}
								\Leftrightarrow&   v\left(f_{(I+L)}(b)\boldsymbol{\sigma} ^L (a-b)\right) > v(f_{(I)}(b)).
							\end{split}
						\end{align}
						Consider  the Taylor expansion of 
						$f_{(I)}(a)$
						around the point 
						$b$. It is 
						\begin{equation*}
							f_{(I)}(a) = f_{(I)}(b) + \sum\limits_{L\neq 0} f_{(I)(L)}(b) \boldsymbol{\sigma} ^L (a-b).    
						\end{equation*}
						This gives
						\begin{equation*}
							v(f_{(I)}(a)) \geq \min\limits_{L\neq 0} \left\lbrace v(f_{(I)}(b)) , v(f_{(I)(L)}(b) \boldsymbol{\sigma} ^L (a-b)) \right\rbrace.
						\end{equation*}
						Clearly, we have $v(f_{(I)(L)}(b) \boldsymbol{\sigma} ^L (a-b))=v(f_{(I + L)}(b) \boldsymbol{\sigma} ^L (a-b))$. Thus, \eqref{I max} implies
						\begin{equation*}
							v(f_{(I)}(a)) \geq \min \left\lbrace v(f_{(I)}(b)) , v(f_{(I + L)}(b) \boldsymbol{\sigma} ^L (a-b)) \right\rbrace = v(f_{(I)}(b)).    
						\end{equation*}
						which gives 
						$v(f_{(I)}(a))= v(f_{(I)}(b))$.
					\end{proof}
					\textbf{Step 2}:If for a multi-index $J$, we have
					$\varepsilon = \varepsilon_J$, then we have
					\begin{equation*}
					\underset{I}{\min}\left(\lvert I \rvert_{\rho}\varepsilon + v(f_{(I)}(b))\right) = \lvert J\rvert_{\rho}\varepsilon + v(f_{(J)}(b)).
					\end{equation*}  
					
					\begin{proof}
						Suppose
						$I$
						is arbitrarily chosen. We want to show that
						\begin{equation*}
							\lvert J\rvert_{\rho}\varepsilon + v(f_{(J)}(b))\leq \lvert I \rvert_{\rho}\varepsilon + v(f_{(I)}(b)).
						\end{equation*}
						We have 
						\begin{align*}
							\varepsilon_I \leq \varepsilon = \varepsilon_J &\Rightarrow  \lvert I \rvert_{\rho}\varepsilon_I + v(f_{(I)}(b)) \leq 
							\lvert I\rvert_{\rho}\varepsilon + v(f_{(I)}(b)) \\
							&\Rightarrow v(f(b)) \leq \lvert I\rvert_{\rho}\varepsilon + v(f_{(I)}(b))\\
							&\Rightarrow \lvert J\rvert_{\rho} \varepsilon_J + v(f_{(J)}(b)) \leq \lvert I\rvert_{\rho} \varepsilon + v(f_{(I)}(b)).
						\end{align*}
						From
						$\varepsilon_J = \varepsilon $, it follows that  
						\begin{equation*}
							\underset{I}{\min}\left( \lvert I \rvert_{\rho}\varepsilon + v(f_{(I)}(b))\right)= \lvert J \rvert_{\rho} \varepsilon + v(f_{(J)}(b)).   
						\end{equation*}
					\end{proof}
					\renewcommand\qedsymbol{$\square$}
					Returning to the main statement, we want to prove $\varepsilon < \varepsilon'$.
					We prove this by contradiction. Suppose 
					$\varepsilon \geq \varepsilon'$, where 
					$\varepsilon = \varepsilon _I$
					(and 
					$I$
					is the maximal multi-index attaining the maximum).
					Also assume
					$\varepsilon'= \varepsilon'_{J}$.
					Thus, we have
					\begin{equation*}
						\varepsilon \geq \varepsilon' \Rightarrow \lvert I\rvert_{\rho} \varepsilon + v(f_{(I)}(b)) \geq \lvert I\rvert_{\rho} \varepsilon' + v(f_{(I)}(b)).   
					\end{equation*}
					From what we proved in step 1, for such
					$I$ 
					we have 
					$v(f_{(I)}(b)) = v(f_{(I)}(a))$.
					Thus,
					$\lvert I\rvert_{\rho} \varepsilon + v(f_{(I)}(b)) \geq \lvert I\rvert_{\rho} \varepsilon'+ v(f_{(I)}(a))$.
					Using Step 2 and the assumption that $\varepsilon= \varepsilon_I$, we can write
					\begin{align*}
						v(f(b)) &\stackrel{\varepsilon= \varepsilon_I}{\eqdef} \lvert I\rvert_{\rho} \varepsilon+ v(f_{(I)}(b))\\ 
						&\stackrel{Step1}{\geq}\lvert I\rvert_{\rho} \varepsilon'+ v(f_{(I)}(a)\\ 
						&\stackrel{Step2}{\geq} |J|_{\rho} \varepsilon'+ v(f_{(J)}(a))\\
						&\stackrel{\varepsilon'= \varepsilon'_{J}}{\eqdef}v(f(a)).
					\end{align*}
					This gives $v(f(b)) \geq v(f((a))$
					which contradicts the condition in part $(1)$. Hence 
					$\varepsilon < \varepsilon'$.
				\end{enumerate}
			\end{proof}
			As we proved in Lemma \ref{l1}, if $b$ is not a root of $f$, we can find a better estimation of a possible root around it. In Proposition \ref{l2}, by using Lemma \ref{l1}, we build a pc-sequence. Assuming that the field $K$ is spherically complete, this pc-sequence necessarily has a pseudolimit, which is a root of $f$. This implies the main result of this section which is Theorem \ref{DiffNewton}.
			\begin{prop}\label{l2}
				Let
				$f\in K_{\sigma}[x]$
				be a nonconstant difference polynomial and assume that
				$b\in K $ is not a root of $f$.
				Define $\varepsilon := \underset{ \underset{ \lvert J\rvert \geq 1}{J}}{\max}\,\,\varepsilon _J $
				where 
				\begin{equation*}
					\varepsilon_J :=\frac{1}{\lvert J\rvert_{\rho}}\left(v(f(b)) - v(f_{(J)}(b)) \right).    
				\end{equation*}
				Similarly, for a point $a_{\eta}$, we define $\varepsilon _{(\eta,J)}$ and also  
				$\varepsilon_{\eta} := \underset{ \underset{ \lvert J\rvert \geq 1}{J}}{\max}\,\,\varepsilon _{(\eta,J)}$.\\ 
				Then there exists a pc-sequence 
				$\{a _{\eta}\}$ 
				in 
				$K$ with the following properties:
				\begin{enumerate}
					\item
					$a_0 = b$ 
					and 
					$\{a _{\eta}\}$
					has a pseudolimit 
					$a\in K$, such that $f(a)=0$ and\\
					$v(a-b)=\varepsilon$;	
					\item
					$\{v(f(a _{\eta}))\}$	
					is strictly increasing;
					\item
					$v(a_{\eta'} - a_{\eta}) = \varepsilon_{\eta}$ 
					whenever 
					$\eta < \eta'$;
					\item
					For any
					$\eta < \eta'$
					we have 
					$\varepsilon _{\eta }< \varepsilon _{\eta'}$. 
				\end{enumerate}
			\end{prop}
			\begin{proof}
				We prove by contradiction. Suppose no such pc-sequence exists.
				Assume 
				$a_0=b$.
				For some ordinal
				$\lambda >0$, by using Lemma \ref{l1}, we inductively construct a sequence $\{a_{\eta}\} _{\eta < \lambda}$ such that for all $\eta$, $a_{\eta}$ is not a root of $f$, and it satisfies properties $(2)$,$(3)$ and $(4)$.
				Then we use transfinite recursion to extend this sequence to an arbitrarily long sequence which arises a contradiction.\\
				From properties $(2)$, $(3)$ and $(4)$,
				we conclude that for all $\eta < \eta'<\eta''< \lambda$ we have
				\begin{equation*}
					v(a_{\eta} - a_{\eta'}) = \varepsilon _{\eta} <\varepsilon _{\eta'}= v(a_{\eta'} - a_{\eta''}),     
				\end{equation*}
				
				which means 
				$\{a_{\eta}\} _{\eta < \lambda}$
				is a pc-sequence.
				We discuss two different possible cases for 
				$\lambda$:
				\begin{enumerate}[label=(\roman*)]
					\item
					$\lambda$
					is a successor ordinal which means it can be written as 
					$\lambda = \mu + 1 $.\\
					By Lemma \ref{l1}, for 
					$f\in K_{\sigma}[x]$
					and
					$a_{\mu} \in K$,
					there exists
					$a_{\lambda} \in K$
					such that 
					\begin{itemize}
						\item 
						$v(a_{\lambda} - a_{\mu}) = \varepsilon _{\mu}$
						and
						$v(f(a_{\lambda})) > v(f(a_{\mu}))$;
						\item $\varepsilon_{\mu} < \varepsilon_{\lambda}$.
					\end{itemize}
					If $a_{\lambda}$ is a root of $f$, these two properties imply $\lbrace v(a_{\lambda}-a_{\eta})\rbrace$ is eventually increasing, which means $a_{\lambda}$ is a pseudolimit of $\{a_{\eta}\} _{\eta < \lambda}$ and we are done. Otherwise, 
					we extend the sequence to 
					$\{a_{\eta}\} _{\eta < \lambda +1}$
					with the same properties. 
					\item
					Let 
					$\lambda$
					is a limit ordinal. As $K$ is spherically complete, $\{a_{\eta}\} _{\eta < \lambda}$ as a pc-sequence has a pseudolimit $a_{\lambda}$. Assume $a_{\lambda}$ is not a root of $f$.
					We want to check whether  
					$a_{\lambda}$
					has the properties $(2)$,$(3)$ and $(4)$. If so, we extend the sequence to 
					$\{a_{\eta}\} _{\eta < \lambda +1}$.\\
					We start with checking $(3)$. Since 
					$a_{\eta} \rightsquigarrow a_{\lambda}$, 
					by definition we have
					\begin{equation} \label{pseudo}
						\exists \eta _0\,\,\, \textrm{such that}\,\,\, \forall\,\, \eta' , \eta ;\,\,\,\eta' > \eta > \eta _0 \Rightarrow
						v(a_{\lambda} - a_{\eta'}) > v(a_{\lambda} - a_{\eta}).
					\end{equation}
					Thus, we have
					\begin{equation}\label{eta}
						\begin{aligned}
							\lambda > \eta +1 > \eta > \eta _0 \Rightarrow 
							v(a_{\eta +1} - a_{\eta}) &= v(a_{\eta +1} - a_{\lambda} + a_{\lambda} - a_{\eta})\\
							&\geq \min\{ v(a_{\eta +1} - a_{\lambda}) ,  v(a_{\eta} - a_{\lambda})\} \\
							&=v(a_{\eta} - a_{\lambda}).
						\end{aligned}
					\end{equation}
					Since the inequality in \eqref{pseudo} is strict, therefore $v(a_{\eta +1} - a_{\lambda}) \neq  v(a_{\eta} - a_{\lambda})$. Then, we have
					$ v(a_{\eta} - a_{\lambda}) =  v(a_{\eta +1} - a_{\eta})= \varepsilon_{\eta }$.\\
					If
					$\gamma$
					is such that 
					$\gamma \leq \eta _0 <\eta$, we can write 
					\begin{equation}\label{gamma}
						\begin{aligned}
							v(a_{\lambda} - a_{\gamma}) =  v(a_{\lambda} - a_{\eta} + a_{\eta} - a_{\gamma})
							&\geq 
							\min\{ v(a_{\lambda} - a_{\eta}) ,  v(a_{\eta} - a_{\gamma})\}\\
							&= \min\{\varepsilon_{\eta} , \varepsilon_{\gamma}\}\\
							&= \varepsilon _{\gamma}.
						\end{aligned}
					\end{equation}
					Therefore, \eqref{eta} and \eqref{gamma} imply that
					$a_{\lambda}$
					satisfies property $(3)$ for the sequence.\\
					We continue by checking $(2)$. Consider the Taylor expansion of 
					$f(a_{\lambda})$
					around the point 
					$a_{\eta}$. This is 
					\begin{equation*}
						f(a_{\lambda}) = f(a_{\eta}) + \sum\limits_{\underset{|J| \geq 1}{J}} f_{(J)}(a_{\eta}) \boldsymbol{\sigma} ^J (a_{\lambda} - a_{\eta}).
					\end{equation*}
					By definition,
					$\varepsilon_{\eta} = \underset{ \underset{ \lvert J\rvert \geq 1}{J}}{\max}\,\,\varepsilon _{(\eta,J)} $
					where 
					$\varepsilon_{(\eta,J)} :=\frac{1}{\lvert J\rvert_{\rho}}\left(v(f(a_{\eta})) - v(f_{(J)}(a_{\eta})) \right) $.\\ 
					For an arbitrary $(\eta,J)$ we have 
					\begin{equation*}
						\varepsilon_{(\eta,J)} \leq \varepsilon_{\eta } \Rightarrow \frac{1}{\lvert J\rvert_{\rho}}\left(v(f(a_{\eta}) - v(f_{(J)}(a_{\eta}))\right) \leq \varepsilon _{\eta}.
					\end{equation*}
					This gives
					\begin{align*}
						v(f(a_{\eta})) &\leq \lvert J\rvert_{\rho} \varepsilon_{\eta } +  v(f_{(J)}(a_{\eta}))\\
						&= \lvert J\rvert_{\rho}v(a_{\eta +1} - a_{\eta}) + v(f_{(J)}(a_{\eta}))\\
						&=|J|_{\rho}v(a_{\lambda} - a_{\eta}) + v(f_{(J)}(a_{\eta}))\\
						&= v(\boldsymbol{\sigma} ^J (a_{\lambda} - a_{\eta})\cdot f_{(J)}(a_{\eta})). 
					\end{align*}
					In the Taylor expansion of $f(a_{\lambda})$, take the valuation of both sides. This yields
					\begin{equation}\label{lambda}
						v(f(a_{\lambda})) \geq \min\limits_{\underset{|J| \geq 1}{J}} \{v(f(a_{\eta})),v(f_{(J)}(a_{\eta}) \boldsymbol{\sigma} ^J (a_{\lambda} - a_{\eta}))  \} = v(f(a_{\eta})).
					\end{equation}
					But the equality can not happen. To see this, suppose $v(f(a_{\lambda}))=v(f(a_{\eta}))$ for some $\eta$. As $\lambda$ is a limit ordinal, we have $\lambda > \eta+1 >\eta$. From the second property of the sequence, this means 
					\begin{equation*}
						v(f(a_{\lambda}))=v(f(a_{\eta}))< v(f(a_{\eta+1})),
					\end{equation*}
					and this contradicts \eqref{lambda}. Hence
					$v(f(a_{\lambda})) > v(f(a_{\eta}))$.\\
					Finally, we check $(4)$. Apply Lemma \ref{l1} for the difference polynomial $f$, a nonroot $a_{\eta}$ and the point
					$a_{\lambda}$ which satisfies the properties of the first part. Thus, the second part of this lemma implies
					$\varepsilon_{\eta } < \varepsilon_{\lambda}$.\\
					Therefore,
					$a_{\lambda}$ has all the properties of the sequence  
					$\{a_{\eta}\}_{\eta < \lambda }$ which enables us to add $a_{\lambda}$ to the sequence.\\
					This means, we can build an arbitrarily long sequence $ \{a_{\eta}\}_{\eta < \lvert K\rvert^+ }$ such that for all $\eta$, $a_{\eta}$ is not a root of $f$, and it satisfies properties $(2)$,$(3)$ and $(4)$. This is a contradiction since all 
					$a_{\eta}$
					are distinct elements in 
					$K$. Hence, for some $\lambda$, $a_{\lambda}$ is a root of $f$.
				\end{enumerate}
			\end{proof}
			\section{Lifting roots and the Difference Kapranov Theorem}
			This section consists of two subsections. We will see the final result which is the Difference Kapranov Theorem in the second subsection. In the first one, we prove Proposition \ref{main} which is the main ingredient to prove the Difference Kapranov Theorem. This is done gradually. In the first place, we are going to prove a simpler version of Proposition \ref{main}, namely for the case where $f$ is a difference polynomial in one variable.
			\par The general assumption of this section is that $K$ is a multiplicative valued difference field of characteristic zero, and is spherically complete. A particular setting is given by Example \ref{Hahn}. We also assume that the valuation has a splitting and $\rho$ the scaling exponent of $\sigma$ is transendental. Besides, we assume that the difference value group $\Gamma$ of $K$  is a divisible subgroup of $\mathbb{R}$ that is a $\mathbb{Q}(\rho)$-module.  
			The difference residue field of $K$ is also assumed to be an ACFA, and of characteristic zero.\\
			\subsection{Liffting Roots}
			We gradually work towards the goal of this part, which is Proposition \ref{main} . First, we will prove a simpler version of this proposition for the case where $f$ is a difference polynomial in one variable.
			
			\begin{lem}\label{1case}
				Suppose $f \in K_{\sigma}\left[x \right] $ is a difference polynomial. Assume $w\in \Gamma$ and $\mathopin_{w}(f)$ is not a monomial. Let $\bar{\alpha}$ be a nonzero root of $\mathopin_{w}(f)$ in the difference residue field $\mathbf{k}$. Then $f$ has a root $a\in K$ such that $v(a)=w$ and $\overline{t^{-w}a}=\bar{\alpha}$.
			\end{lem}
			\begin{proof}
				Choose $\alpha$ as a representative of $\bar{\alpha}$ and set $b=t^w \alpha \in K$. Then clearly we have $v(b)=v(t^w \alpha)=v(t^w)+v(\alpha)=w+0=w$. In addition, $\overline{t^{-w}b}=\overline{t^{-w}\cdot t^{w} \alpha}=\bar{\alpha}$.
				Applying Theorem \ref{DiffNewton}, for $f \in K_{\sigma}[x]$ and $b\in K$, there exists a root $a\in K$ such that $v(a-b)=\varepsilon$ (where $\varepsilon$ is as defined in the same theorem ). It suffices to prove that this root $a$ satisfies the desired properties.\\
				\begin{claim}
					Let $a \in K$ be a root of $f$ which is obtained by applying Theorem \ref{DiffNewton}, and $\varepsilon$ be as defined in this theorem. If $\varepsilon > w$, then this root satisfies the following properties:
					\begin{itemize}
						\item $v(a)=w$;
						\item $\overline{t^{-w}a}=\bar{\alpha}$.
					\end{itemize}
				\end{claim}
				\renewcommand\qedsymbol{$\blacksquare$}
				\begin{proof}
					We write
					\begin{equation*}
						v(a)=v(a-b+b)\geq \min \left\lbrace v(a-b),v(b) \right\rbrace=\min\left\lbrace \varepsilon,w \right\rbrace=w. 
					\end{equation*}
					Since the minimum is attained uniquely, we have $v(a)=w$.\\
					Moreover, we have 
					\begin{align*}
						\varepsilon - w >0 &\Rightarrow v(a-b)+v(t^{-w})>0\\
						&\Leftrightarrow v\left( t^{-w}a-t^{-w}b \right)>0 \\
						&\Leftrightarrow \overline{t^{-w}a}=\overline{t^{-w}b}=\bar{\alpha}.
					\end{align*}
					This means both conditions are satisfied by the root $a$.   
				\end{proof}
				\renewcommand\qedsymbol{$\square$}
				Now, it suffices to prove that $\varepsilon > w$.\\
				As $f\in K_{\sigma}[x]$, it is of the form $f(x)=\underset{J\in \Lambda}{\sum}c_J \boldsymbol{\sigma}^J(x)$ where $\Lambda$ is a finite subset of $\mathbb{N}^{n+1}$. The tropicalization of $f$ at $w$ is
				\begin{align*}\label{trop(w)}
					\mathoptrop(f)(w)&=\underset{J\in \Lambda}{\min}\left\lbrace v(c_J)+J \boldsymbol{\sigma}_{\Gamma}(w) \right\rbrace\\   &=\underset{J\in \Lambda}{\min}\left\lbrace v(c_J) + (j_0,j_1,\dots ,j_n)\cdot (w,\rho\cdot w,\dots,\rho^n\cdot w)\right\rbrace\\
					&=\underset{J\in \Lambda}{\min}\left\lbrace v(c_J) + \lvert J \rvert_{\rho}\cdot w\right\rbrace.
				\end{align*}
				
				
				By the assumptions, $\mathopin_w(f)$ is not a monomial. Suppose $\Delta \subseteq \Lambda$ consists of all those multi-indices whose corresponding monomials appear in $\mathopin_w(f)$. Obviously $\Delta$ has more than one element.\\
				From the definition of $\mathopin_w(f)$, for any $I \in \Delta$, we can write
				\begin{equation}\label{initial}
					\mathoptrop(f)(w)=\underset{J\in \Lambda}{\min}\left\lbrace v(c_J) + \lvert J \rvert_{\rho}\cdot w\right\rbrace=v(c_I)+\lvert I \rvert_{\rho}\cdot w.
				\end{equation}
				On the other hand, we have 
				\begin{equation*}
					f(b)=\underset{J \in \Lambda}{\sum}c_J \boldsymbol{\sigma}^J(b) .
				\end{equation*}
				This implies 
				\begin{equation}\label{val}
					v(f(b))\geq 
					\underset{J\in \Lambda}{\min}\left\lbrace v(c_J)+ \lvert J \rvert _{\rho}v(b)\right\rbrace
					. 
				\end{equation}
				What we obtained in \eqref{val} and \eqref{initial} together give
				\begin{equation}\label{delta}
					\forall I\in \Delta, \,\,\,\,\,v(f(b))\geq \underset{J\in \Lambda}{\min}\left\lbrace v(c_J)+ \lvert J \rvert _{\rho}\cdot w\right\rbrace =v(c_I)+ \lvert I \rvert _{\rho}\cdot w
					.
				\end{equation}
				By $J \geq I$, we mean, for all $r$ such that $0\leq r \leq n$, we have $j_r \geq i_r $. If we also consider
				Definition \ref{f_J}, then
				for any $I\in \Delta$, we have
				\begin{align*}
					f_{(I)}(b)=P_{(I)}(\boldsymbol{\sigma}(b))=&\frac{\partial^{\lvert I\rvert}P\left(b,\sigma(b),\dots ,\sigma^n(b) \right) }{\partial x_0^{i_0} \partial x_1^{i_1}\dots \partial x_n^{i_n}}\cdot \frac{1}{i_0!i_1!\dots i_n!}\\
					=&\left.\frac{\partial^{\lvert I \rvert} \underset{J \in \Lambda}{\sum}c_J \boldsymbol{x}^J}{\partial x_0^{i_0} \partial x_1^{i_1}\dots \partial x_n^{i_n}}\cdot \frac{1}{i_0!i_1!\dots i_n!}\right|_{\boldsymbol{x}=\boldsymbol{\sigma}(b)}\\
					=&\left.\underset{J \geq I}{\sum}c_J \frac{\partial^{\lvert I \rvert} \boldsymbol{x}^J}{\partial x_0^{i_0} \partial x_1^{i_1}\dots \partial x_n^{i_n}}\cdot \frac{1}{i_0!i_1!\dots i_n!}\right|_{\boldsymbol{x}=\boldsymbol{\sigma}(b)}\\
					=&\left.\underset{J \geq I}{\sum}c_J\left( \frac{J!}{(J-I)!}\cdot \boldsymbol{x}^{(J-I)}\cdot\frac{1}{I!} \right)\right|_{\boldsymbol{x}=\boldsymbol{\sigma}(b)}\\ 
					=&\left.\underset{J \geq I}{\sum}c_J\left( \binom{J}{I}\cdot \boldsymbol{x}^{(J-I)} \right)\right|_{\boldsymbol{x}=\boldsymbol{\sigma}(b)}.
				\end{align*}
				
				Here by $\binom{J}{I}$ we mean $\binom{j_0}{i_0}\binom{j_1}{i_1}\cdots \binom{j_n}{i_n}$. Denoting this coefficient by $\alpha_J$
				we have
				\begin{equation*}
					f_{(I)}(b)=\underset{J \geq I}{\sum}c_J\cdot \alpha_J\cdot \boldsymbol{\sigma}^{(J-I)}(b).
				\end{equation*}
				Let $I_m$ be the maximal multi-index in $\Delta$ with respect to lexicographical order.\\
				Then we have 
				\begin{equation}\label{mderiv}
					f_{(I_m)}(b)=\underset{J \geq I_m}{\sum}c_J\cdot \alpha_J\cdot \boldsymbol{\sigma}^{(J-I_m)}(b).
				\end{equation}
				Since $I_m$ is the greatest element in $\Delta$, there is no element of $\Delta$ appearing in $f_{(I_m)}(b)$. Since $\alpha_J$ is a natural number, we have $v(\alpha_J)=0$. Therefore, \eqref{mderiv} and \eqref{delta} imply
				\begin{align*}
					v(f_{(I_m)}(b)) &\geq \underset{J\geq I_m}{\min}\left\lbrace v(c_J)+v(\alpha_J)+\lvert J-I_m \rvert_{\rho}v(b) \right\rbrace\\
					&= \underset{J\geq I_m}{\min}\left\lbrace v(c_J)+\lvert J \rvert_{\rho}w\right\rbrace - \lvert I_m \rvert_{\rho}w\\
					&=v(c_{I_m})+ \lvert I_m \rvert_{\rho}w- \lvert I_m \rvert_{\rho}w=v(c_{I_m}).
				\end{align*}
				
				On the right side of the above inequality, the minimum is attained only once. Hence, 
				\begin{equation}\label{vmderiv}
					v(f_{(I_m)}(b))=v(c_{I_m}).
				\end{equation}
				As $I_m \in \Delta$, \eqref{delta} gives 
				\begin{equation}\label{valfb}
					v(f(b))\geq v(c_{I_m})+ \lvert I_m \rvert_{\rho}w.
				\end{equation}
				From \eqref{vmderiv} together with \eqref{valfb}, we have
				\begin{equation*}
					v(f(b))-v(f_{(I_m)}(b)) \geq \lvert I_m\rvert_{\rho}\cdot w.
				\end{equation*}
				Multiply both sides by $\dfrac{1}{\lvert I_m\rvert_{\rho}}$ and consider the definition of $\varepsilon_{I_m}$. Thus, 
				we have 
				\begin{equation*}
					\varepsilon_{I_m}=\dfrac{1}{\lvert I_m \rvert_{\rho}}\left(v(f(b))-v(f_{(I_m)}(b)) \right) \geq w.
				\end{equation*}
				Therefore, we obtain
				\begin{equation*}
					\varepsilon=\underset{ \underset{ \lvert J\rvert \geq 1}{J}}{\max}\,\,\varepsilon_{J} \geq \varepsilon_{I_m} \geq w.
				\end{equation*}
				Thus far, we have $\varepsilon \geq w$. We show that the equality cannot occur. We prove this in two steps.
				\begin{itemize}
					\item \textbf{Step 1}: We have
					\begin{equation}\label{initialalpha}
						\mathopin_w(f)(\bar{\alpha})=\underset{I \in \Delta}{\sum}\overline{t^{-v(c_I)}c_I}\boldsymbol{\bar{\sigma}}^{I}(\bar{\alpha}).
					\end{equation}
					From the assumptions, $\bar{\alpha}$ is a root of $\mathopin_w(f)$, and \eqref{initialalpha} is zero.\\
					On the other hand, we take $f(b)=\underset{J\in \Lambda}{\sum}c_J\boldsymbol{\sigma}^J(b)$, and divide the both sides of this equation by $t^{v\left(c_{I_m}\boldsymbol{\sigma}^{I_m}(b) \right) }$. Since \eqref{valfb} implies that
					$v\left(\dfrac{f(b)}{t^{v\left(c_{I_m}\boldsymbol{\sigma}^{I_m}(b) \right) }}\right)$ is nonnegative, we can consider the residue of both sides. This gives
					\begin{equation}\label{res}
						\overline{\dfrac{f(b)}{t^{v\left(c_{I_m}\boldsymbol{\sigma}^{I_m}(b) \right) }}}=\overline{\underset{J \in \Lambda}{\sum}t^{-v\left(c_{I_m}\boldsymbol{\sigma}^{I_m}(b) \right) }c_J\boldsymbol{\sigma}^{J}(b)}.
					\end{equation}
					In this step, we mainly prove that 
					\begin{equation}\label{equation}
						\overline{\dfrac{f(b)}{t^{v\left(c_{I_m}\boldsymbol{\sigma}^{I_m}(b) \right) }}}=\mathopin_w(f)(\bar{\alpha})=0.   
					\end{equation}
					\item \textbf{Step 2}: In this final step, we use \eqref{equation} to show that $\varepsilon > w$. 
				\end{itemize}
				\textbf{Step 1:} We prove \eqref{equation}.
				\renewcommand\qedsymbol{$\blacksquare$}
				\begin{proof}  
					Suppose $J \in \Lambda$, the valuation of the J-th summand appearing on the right side of \eqref{res}, is
					\begin{align*}\label{val}
						v\left(t^{-v\left(c_{I_m}\boldsymbol{\sigma}^{I_m}(b) \right) }c_J\boldsymbol{\sigma}^{J}(b) \right)
						&=-v\left(c_{I_m}\boldsymbol{\sigma}^{I_m}(b) \right)+v(c_J)+v(\boldsymbol{\sigma}^J(b))\\
						&=-v(c_{I_m})-\lvert I_m\rvert_{\rho}w+v(c_J)+\lvert J \rvert_{\rho}w.   
					\end{align*}
					As we discussed before, $v(c_{I_m})+\lvert I_m\rvert_{\rho}w$ attains the minimum. This means that\\ $v\left(t^{-v\left(c_{I_m}\boldsymbol{\sigma}^{I_m}(b) \right) }c_J\boldsymbol{\sigma}^{J}(b) \right)$ is nonnegative.
					Also note that for all $J \in \Lambda \setminus \Delta$, we have
					\begin{equation*}
						v\left(t^{-v\left(c_{I_m}\boldsymbol{\sigma}^{I_m}(b) \right) }c_J\boldsymbol{\sigma}^{J}(b) \right) > 0.
					\end{equation*}
					Hence, \eqref{res} can be written as: 
					\begin{equation}\label{rsep}
						\overline{\dfrac{f(b)}{t^{v\left(c_{I_m}\boldsymbol{\sigma}^{I_m}(b) \right) }}}=\underset{I\in \Delta}{\sum}\overline{t^{-v\left(c_{I_m}\boldsymbol{\sigma}^{I_m}(b) \right) }c_I\boldsymbol{\sigma}^{I}(b)}.
					\end{equation} 
					We want to show that each summand appearing on the right side of \eqref{rsep} coincides with the corresponding summand appearing in $\mathopin_w(f)(\bar{\alpha})$.\\
					Suppose $I \in \Delta$ and $\overline{t^{-v\left(c_{I_m}\boldsymbol{\sigma}^{I_m}(b) \right) }c_I\boldsymbol{\sigma}^{I}(b)}$ is the corresponding summand in \eqref{rsep}. 
					From the definition of $\boldsymbol{\sigma}$ , it can be written as:
					\begin{equation}\label{last}
						\begin{aligned}
							\overline{t^{-v\left(c_{I_m}\boldsymbol{\sigma}^{I_m}(b) \right) }c_I\boldsymbol{\sigma}^{I}(b)}\cdot\dfrac{\boldsymbol{\sigma}^{I}(\bar{\alpha})}{\boldsymbol{\sigma}^{I}(\bar{\alpha})}&= \overline{\dfrac{t^{-v\left(c_{I_m}\boldsymbol{\sigma}^{I_m}(b) \right) }c_I\boldsymbol{\sigma}^{I}(b)}{\boldsymbol{\sigma}^{I}(t^{-w}b)}}\cdot \boldsymbol{\sigma}^{I}(\bar{\alpha})\\
							&=\overline{t^{-v(c_{I_m})-\lvert I_m \rvert_{\rho}w+\lvert I \rvert_{\rho}w}c_I}\cdot \boldsymbol{\sigma}^{I}(\bar{\alpha}).
						\end{aligned}
					\end{equation}
					Since $v(c_{I_m})+\lvert I_m \rvert_{\rho}w = v(c_I)+\lvert I \rvert_{\rho}w$, \eqref{last} equals
					$\overline{t^{-v(c_I)}c_I}\boldsymbol{\sigma}^{I}(\bar{\alpha})$ which is the I-th term of $\mathopin_w(f)(\bar{\alpha})$.
					This yields
					\begin{equation*}
						\overline{\dfrac{f(b)}{t^{v\left(c_{I_m}\boldsymbol{\sigma}^{I_m}(b) \right) }}}=\mathopin_w(f)(\bar{\alpha})=0.
					\end{equation*}
				\end{proof}
				\textbf{Step 2:} We prove $\varepsilon > w$.
				\begin{proof} 
					From Step 1, we have $\overline{\dfrac{f(b)}{t^{v\left(c_{I_m}\boldsymbol{\sigma}^{I_m}(b) \right) }}}=0$.
					This means that 
					\begin{equation*}
						v\left(\dfrac{f(b)}{t^{v\left(c_{I_m}\boldsymbol{\sigma}^{I_m}(b) \right) }} \right)>0.   
					\end{equation*}
					Equivalently, this gives
					\begin{equation}\label{neq}
						v(f(b))>v\left(c_{I_m}\boldsymbol{\sigma}^{I_m}(b) \right)=v(c_{I_m})+\lvert I_m \rvert_{\rho}w. 
					\end{equation}
					From the equality in \eqref{vmderiv}, \eqref{neq} can be written as:
					\begin{equation*}
						v(f(b))-v(f_{(I_m)}(b))>\lvert I_m \rvert_{\rho}w.
					\end{equation*}
					Dividing both sides of this inequality by $\lvert I_m \rvert_{\rho}$ gives $\varepsilon_{I_m}>w$. Therefore, $\varepsilon > w$.
				\end{proof}
				\renewcommand\qedsymbol{$\square$}
			\end{proof}
			In Lemma \ref{n=1}, we generalize Lemma \ref{1case} to a Laurent difference polynomial in one variable.
			\begin{lem}\label{n=1}
				Suppose $f \in K_{\sigma}[x^{\pm1}]$ is a Laurent difference polynomial. Assume that $w \in \Gamma$ and $\mathopin_w(f)$ is not a monomial. Let $\bar{\alpha}$ be a nonzero root of $\mathopin_w(f)$ in the difference residue field $\mathbf{k}$. Then f has a root $a \in K$ such that $v(a)=w$ and $\overline{t^{-w}a}=\bar{\alpha}$.  
			\end{lem}
			\begin{proof}
				If $f \in K_{\sigma}[x^{\pm1}]$, from Remark \ref{Jmax}, there exists $g(x)\in K_{\sigma}[x]$, such that 
				\begin{equation*}
					g(x)=f(x)\cdot \boldsymbol{\sigma}^{J_{\max}}(x).    
				\end{equation*}
				By Lemma \ref{in(fg)}, we have
				\begin{equation*}
					\mathopin_w(g)=\mathopin_w(f\cdot \boldsymbol{\sigma}^{J_{\max}})=\mathopin_w(f)\cdot \mathopin_w(\boldsymbol{\sigma}^{J_{\max}}).
				\end{equation*}
				By the assumptions $\mathopin_w(f)$ is not a monomial. Hence, $\mathopin_w(g)$ is not a monomial. Moreover,
				\begin{equation*}
					\mathopin_w(g)(\bar{\alpha})=\mathopin_w(f)(\bar{\alpha})\cdot \mathopin_w(\boldsymbol{\sigma}^{J_{\max}})(\bar{\alpha})=0,
				\end{equation*}
				which means $\bar{\alpha}$ is a root of $\mathopin_w(g)$. Therefore, we can apply Lemma \ref{1case} to obtain a root $a \in K$ for $g$ such that
				\begin{itemize}
					\item $v(a)=w$ and \item $\overline{t^{-w}a}=\bar{\alpha}$. 
				\end{itemize}
				Hence, $f(a)\boldsymbol{\sigma}^{J_{\max}}(a)=g(a)=0$. As $\bar{\alpha} \in \mathbf{k}^*$ and $\overline{t^{-w}a}=\bar{\alpha}$, we have $a \neq 0$.
				This means $\boldsymbol{\sigma}^{J_{\max}}(a)\neq 0$ which implies $f(a)=0$. We also have $v(a)=w$ and
				$\overline{t^{-w}a}=\bar{\alpha}$.
			\end{proof}
			Moving another step forward to obtain Proposition \ref{main}, we prove a similar statement in Proposition \ref{3.1.5} for a Laurent difference polynomial $f$ in $n$ variables, with an extra assumption on $f$.
			\begin{prop}\label{3.1.5}
				Let $f \in K_{\sigma}[x_1^{\pm1},\dots,x_n^{\pm1}]$ be a Laurent difference polynomial with different $\sigma$-powers of $x_n$ in its different monomials (see Notation \ref{sigma}). Assume $\underline{w}=\left(w_1, \dots, w_n\right) \in \Gamma^n$ such that $\mathopin_{\underline{w}}(f)$ is not a monomial. Suppose $\bar{\alpha}=\left(\bar{\alpha}_1, \dots, \bar{\alpha}_n\right)$ is a root of $\mathopin_{\underline{w}}(f)$ in $(\mathbf{k^*})^n$. Then there exists an element $y=\left(y_1, \dots, y_n\right)$  in $(K^*)^n$ that is a root of $f$, and satisfies the following conditions:
				\begin{itemize}
					\item $v(y)=\underline{w}$,
					\item for all $i,\,\,\,1\leq i \leq n$ we have $\overline{t^{-w_i}\cdot y_i}=\bar{\alpha}_i$.
				\end{itemize}  
			\end{prop}
			\begin{proof}
				First, we look at the case $n=1$. By Lemma \ref{n=1}, there exists a root $y$ for $f$ which satisfies the desired conditions. 
				The difference polynomial $f$ can be regarded as a difference polynomial in one variable $x_n$ with coefficients in $ K_{\sigma}[x_1^{\pm1},\dots,x_{n-1}^{\pm1}]$.\\
				To see this better, we use the following notation.
				\begin{notn}\label{notation}
					Given $\underline{u}(\sigma):=(u_1(\sigma), \dots ,u_{n-1}(\sigma),u_n(\sigma)) \in \left(\mathbb{Z}[\sigma]\right)^{n}$. Set $u(\sigma):=(u_1(\sigma), \dots ,u_{n-1}(\sigma))$. This means that for a monomial $c_{\underline{u}(\sigma)}x^{u(\sigma)}x_n^{u_n(\sigma)}$, we have $c_{\underline{u}(\sigma)} \in K$, and by $x^{u(\sigma)}$, we mean $x_1^{u_1(\sigma)}\cdots x_{n-1}^{u_{n-1}(\sigma)}$.    
				\end{notn}
				By using this notation each monomial of $f$ takes the form $c_{\underline{u}(\sigma)}x^{u(\sigma)}x_n^{u_n(\sigma)}$.
				We regard $f$ in this manner, and we use it to define a one-variable Laurent difference polynomial $g$. Subsequently, we apply the case $n=1$ to find a root for $g$. This root ultimately yields a root of $f$ that satisfies the desired conditions.\\
				Suppose for each $i$, $1 \leq i \leq n-1$, $w_i$ and $\Bar{\alpha}_i\neq 0$ are given. Choose $\alpha_i$ to be a representative of $\Bar{\alpha}_i$. Define $y_i=t^{w_i}\alpha_i$.
				Then, for each $i$, \,\,\, $1\leq i \leq n-1$, $y_i$ satisfies the following conditions:
				\begin{itemize}
					\item $v(y_i)=w_i$,
					\item $\overline{t^{-w_i}y_i}= \bar{\alpha}_i$. 
				\end{itemize}
				Note that, $\alpha_i \neq 0$, and clearly $y_i$ is nonzero.\\ 
				Since, for each $i$, we have $y_i \in K^*$, and all monomials in $f$ have different $\sigma$-powers of $x_n$, $g(x_n)=f(y_1,\dots,y_{n-1},x_n)$ is a nonzero Laurent difference polynomial. Now, we will find $y_n$ such that $y=\left(y_1,\dots,y_n\right)$ is a root of $f$ satisfying the intended conditions. \\
				Using Notation \ref{notation}, $g(x_n)$ can be written as:
				\begin{equation*}
					g(x_n)=\underset{\underline{u}(\sigma)}{\sum}d_{\underline{u}(\sigma)}x_n^{u_n(\sigma)},
				\end{equation*}
				where $d_{\underline{u}(\sigma)}=c_{\underline{u}(\sigma)}y^{u(\sigma)}$ and $u_n(\sigma) \in \mathbb{Z}[\sigma]$ which is of the form $u_n(\sigma)=\sum_{j_n=1}^{m_n} a_{j_n} \sigma^{j_n}$. So the tropicalization of each monomial of $g(x_n)$ is as follows:
				\begin{equation}\label{tropicmon}
					\begin{aligned}
						\mathoptrop\left(d_{\underline{u}(\sigma)}x_n^{u_n(\sigma)}\right)(w_n)
						&=v(d_{\underline{u}(\sigma)})+\mathoptrop\left(\prod_{j_n=1}^{m_n}\left( \sigma^{j_n}(x_n)\right) ^{a_{j_n}} \right)(w_n)\\
						&=v(d_{\underline{u}(\sigma)})+\sum_{j_n=1}^{m_n}a_{j_n}\sigma_{\Gamma}^{j_n}(w_n)\\
						&=v(c_{\underline{u}(\sigma)})+v(y^{u(\sigma)})+u_n(\rho)\cdot w_n.
					\end{aligned}
				\end{equation}
				As $u(\sigma) \in (\mathbb{Z}[\sigma])^{n-1} $, and for all $i$ with $1\leq i \leq n-1$, we have $v(y_i)=w_i$, we can write
				\begin{equation*}
					v(y^{u(\sigma)})=v(\prod_{i=1}^{n-1}y_i^{u_i(\sigma)})=\sum_{i=1}^{n-1} u_i(\rho)\cdot w_i.
				\end{equation*}
				Hence, from \eqref{tropicmon}, we have 
				\begin{align*}
					\mathoptrop\left(d_{\underline{u}(\sigma)}x_n^{u_n(\sigma)}\right)(w_n)
					&=v(c_{\underline{u}(\sigma)})+ \sum_{i=1}^{n} u_i(\rho)\cdot w_i.
				\end{align*}
				Suppose for $\underline{w}=(w_1,\dots,w_{n-1},w_n)$, for all $i$ with $1\leq i \leq n-1$, we have $w_i=v(y_i)$, and $w_n$ is a variable. Then
				\begin{equation}\label{tropmonog}
					\mathoptrop\left(d_{\underline{u}(\sigma)}x_n^{u_n(\sigma)}\right)(w_n)
					=v(c_{\underline{u}(\sigma)})+\underline{u}(\rho)\cdot \underline{w} ,
				\end{equation} 
				which is a monomial of $\mathoptrop(g)(w_n)$.\\
				On the other hand, consider the tropicalization of a monomial $c_{\underline{u}(\sigma)}x^{u(\sigma)}x_n^{u_n(\sigma)}$ in $f(x_1,\dots,x_{n-1},x_n)=\sum c_{\underline{u}(\sigma)}x^{u(\sigma)}x_n^{u_n(\sigma)}$, where for each $i$, we have $u_i(\sigma) \in \mathbb{Z}[\sigma]$,  i.e $u_i(\sigma)=\sum_{j_i=1}^{m_i} a_{j_i} \sigma^{j_i}$. It is of the following form:
				\begin{equation}\label{tropmonf}
					\begin{aligned}
						\mathoptrop\left(c_{\underline{u}(\sigma)}x^{u(\sigma)}x_n^{u_n(\sigma)}\right)(\underline{w})&=v(c_{\underline{u}(\sigma)})+\sum_{i=1}^{n}\mathoptrop\left(x_i^{u_i(\sigma)}\right)(w_i)\\
						&=v(c_{\underline{u}(\sigma)})+\sum_{i=1}^{n}\sum_{j_i=1}^{m_i}a_{j_i}\sigma_{\Gamma}^{j_i}(w_i)\\
						&=v(c_{\underline{u}(\sigma)})+\underline{u}(\rho)\cdot \underline{w}.
					\end{aligned}
				\end{equation}
				
				The equations in \eqref{tropmonog} and \eqref{tropmonf} together imply $\mathoptrop\left(d_{\underline{u}(\sigma)}x_n^{u_n(\sigma)}\right)(w_n)=\mathoptrop\left(c_{\underline{u}(\sigma)}x^{u(\sigma)}x_n^{u_n(\sigma)}\right)(\underline{w})$. So in general, $\mathoptrop(g)(w_n)=\mathoptrop(f)(\underline{w})$.\\
				Besides, we have
				\begin{align*}
					\mathopin_{w_n}(g)&=\underset{\underline{u}(\sigma):v(d_{\underline{u}(\sigma)})+u_n(\rho)\cdot w_n=\mathoptrop(g)(w_n)}{\sum}\overline{t^{-v(d_{\underline{u}(\sigma)})}d_{\underline{u}(\sigma)}}\cdot x_n^{u_n(\sigma)}\\
					&=\underset{\underline{u}(\sigma):v(c_{\underline{u}(\sigma)})+u(\rho)\cdot w+u_n(\rho)\cdot w_n=\mathoptrop(g)(w_n)}{\sum}\overline{t^{-v(c_{\underline{u}(\sigma)})}c_{\underline{u}(\sigma)}t^{-u(\rho)\cdot w}y^{u(\sigma)}}\cdot x_n^{u_n(\sigma)},
				\end{align*}
				where by $w$ we mean $(w_1,\dots ,w_{n-1})$.\\
				Since $v\left(t^{-u(\rho)\cdot w}y^{u(\sigma)} \right)=0$, and also $v(t^{-v(c_{\underline{u}(\sigma)})}c_{\underline{u}(\sigma)})=0$, we have\\
				
				\begin{equation}\label{in}
					\begin{aligned}
						\mathopin_{w_n}(g)&=\underset{\underline{u}(\sigma):v(c_{\underline{u}(\sigma)})+\underline{u}(\rho)\cdot \underline{w}=\mathoptrop(f)(\underline{w})}{\sum}\overline{t^{-v(c_{\underline{u}(\sigma)})}c_{\underline{u}(\sigma)}}\,\,\,\overline{t^{-u(\rho)\cdot w}y^{u(\sigma)} }\cdot x_n^{u_n(\sigma)}\\
						&=\underset{\underline{u}(\sigma):v(c_{\underline{u}(\sigma)})+\underline{u}(\rho)\cdot \underline{w}=\mathoptrop(f)(\underline{w})}{\sum}\overline{t^{-v(c_{\underline{u}(\sigma)})}c_{\underline{u}(\sigma)}}\,\,\,\overline{\prod_{i=1}^{n-1}t^{-u_i(\rho).w_i}y_i^{u_i(\sigma)}}\cdot x_n^{u_n(\sigma)}\\
						&=\underset{\underline{u}(\sigma):v(c_{\underline{u}(\sigma)})+\underline{u}(\rho)\cdot \underline{w}=\mathoptrop(f)(\underline{w})}{\sum}\overline{t^{-v(c_{\underline{u}(\sigma)})}c_{\underline{u}(\sigma)}}\,\,\,\overline{\prod_{i=1}^{n-1}(t^{-w_i}y_i)^{u_i(\sigma)}}\cdot x_n^{u_n(\sigma)}.
					\end{aligned}
				\end{equation}
				As for all $i,\,\,\, 1 \leq i \leq n-1$, $v(y_i)=w_i$,  we have 
				\begin{align*}
					v\left(\left(t^{-w_i}y_i \right)^{u_i(\sigma)} \right)&=v\left(\prod_{j_i=1}^{m_i} \sigma^{j_i}\left( \left(t^{-w_i}y_i\right)^{a_{j_i}}\right) \right)\\
					&=\sum_{j_i=1}^{m_i}\rho^{j_i}a_{j_i}v(t^{-w_i}y_i)=0, 
				\end{align*}
				that allows us to write
				\begin{equation}\label{resi}
					\overline{\prod_{i=1}^{n-1}\left(t^{-w_i}y_i \right)^{u_i(\sigma)}}= \prod_{i=1}^{n-1}\overline{\left(t^{-w_i}y_i \right)^{u_i(\sigma)}}.
				\end{equation}
				The automorphism $\sigma$ on $K$ induces an automorphism $\Bar{\sigma}$ on the residue field, such that for all $\bar{x}\in \mathbf{k}$, we have $\Bar{\sigma}(\bar{x})=\overline{\sigma(x)}$. By abuse of notation, we also denote $\Bar{\sigma}$ by $\sigma$ so \eqref{resi} can be written as: 
				\begin{equation*}
					\prod_{i=1}^{n-1}\overline{\left(t^{-w_i}y_i \right)^{u_i(\sigma)}}
					=\prod_{i=1}^{n-1}\left(\overline{t^{-w_i}y_i} \right)^{u_i(\sigma)}.
				\end{equation*}
				Finally, \eqref{in} can be written as:
				\begin{equation}\label{inf=ing}
					\begin{aligned}
						\mathopin_{w_n}(g)(x_n)&=\underset{\underline{u}(\sigma):v(c_{\underline{u}(\sigma)})+\underline{u}(\rho)\cdot \underline{w}=\mathoptrop(f)(\underline{w})}{\sum}\overline{t^{-v(c_{\underline{u}(\sigma)})}c_{\underline{u}(\sigma)}}\,\,\,\prod_{i=1}^{n-1}(\overline{t^{-w_i}y_i})^{u_i(\sigma)}\cdot x_n^{u_n(\sigma)}\\
						&=\underset{\underline{u}(\sigma):v(c_{\underline{u}(\sigma)})+\underline{u}(\rho)\cdot \underline{w}=\mathoptrop(f)(\underline{w})}{\sum}\overline{t^{-v(c_{\underline{u}(\sigma)})}c_{\underline{u}(\sigma)}}\,\,\,\prod_{i=1}^{n-1}\bar{\alpha}_i^{u_i(\sigma)}\cdot x_n^{u_n(\sigma)}\\
						&=\mathopin_{\underline{w}}(f)(\bar{\alpha}_1,\dots,\bar{\alpha}_{n-1},x_n).
					\end{aligned}
				\end{equation}
				By the assumptions, $\bar{\alpha}$ is a root of $\mathopin_{\underline{w}}(f)$. This means
				\begin{equation*}
					\mathopin_{w_n}(g)(\bar{\alpha}_n)=\mathopin_{\underline{w}}(f)(\bar{\alpha}_1,\dots,\bar{\alpha}_{n-1},\bar{\alpha}_n)=0.
				\end{equation*}
				Note that $f$ is a Laurent difference polynomial with different $\sigma$-powers of $x_n$ in its different monomials.
				From the definition of the initial forms, as $\mathopin_{\underline{w}}(f)$ is not a monomial, therefore $\mathopin_{\underline{w}}(f)(\bar{\alpha}_1,\dots ,\bar{\alpha}_{n-1},x_n)$ is not a monomial.
				From \eqref{inf=ing}, it is implied that $\mathopin_{w_n}(g)$ is not a monomial. \\
				To sum up, we have $g$ as a Laurent difference polynomial in one variable $x_n$ and that $w_n \in \Gamma$ such that $\mathopin_{w_n}(g)$ is not a monomial. We also know that $\bar{\alpha}_n \in \mathbf{k}^*$ is a root of $\mathopin_{w_n}(g)$.
				From the case $n=1$, there exists a point $y_n \in K^*$ such that 
				\begin{itemize}
					\item $0=g(y_n)=f(y_1,\dots,y_{n-1},y_n)$,
					\item $v(y_n)=w_n$,
					\item $\overline{t^{-w_n}y_n}=\bar{\alpha}_n$.
				\end{itemize}
				Hence, $(y_1,\dots,y_n)$ is a root of $f$ which satisfies the desired properties of the statement.
			\end{proof}
			In Proposition \ref{3.1.5}, we imposed a condition on $f$. Lemma \ref{phistar} shows that this does not lead to loss of generality. More precisely, even if $f$ is an arbitrary Laurent difference polynomial, associated to $f$, one can find a Laurent difference polynomial $g$ with the desired condition.
			\begin{lem}\label{phistar}
				Let $f$ be a Laurent difference polynomial in $K_{\sigma}[x_1^{\pm1},\dots,x_n^{\pm1}]$. For a natural number $l$, we define the following automorphism: 
				\begin{equation*}
					\phi^{*}_l:K_{\sigma}[x_1^{\pm1},\dots,x_n^{\pm1}] \longrightarrow K_{\sigma}[x_1^{\pm1},\dots,x_n^{\pm1}],
				\end{equation*}
				such that, for all $i$ with $1 \leq i \leq n-1$, we have $\phi^*_l(x_i):=x_ix_n^{l^{i}}$. We also define  $\phi^*_l(x_n):=x_n$, and $\phi^*_l(x^{\sigma}):=(\phi^*_l(x))^{\sigma}$ which means that $\phi^*_l(f\left(x_1,\dots,x_n\right)):=f\left(x_1x_n^{l^{1}},\dots ,x_{n-1}x_n^{l^{n-1}},x_n\right)$.
				Then, for a large enough $l$, $g:=\phi^*_l(f)$ is a Laurent difference polynomial with different $\sigma$-powers of $x_n$ in its different monomials.  
			\end{lem}
			\begin{proof}
				Using Notation \ref{notation}, for a monomial $x^{u(\sigma)}x_n^{u_n(\sigma)}$, we have 
				\begin{align*}
					\phi^*_l\left( x^{u(\sigma)}x_n^{u_n(\sigma)}\right)
					&=\prod_{i=1}^{n-1} \left( x_ix_n^{l^i}\right)^{u_i(\sigma)} x_n^{u_n(\sigma)}\\
					&= x^{u(\sigma)} x_n^{u_n(\sigma) +\sum_{i=1}^{n-1}u_i(\sigma)l^i}.
				\end{align*}
				Suppose $x^{u(\sigma)}x_n^{u_n(\sigma)}$ and $x^{u'(\sigma)}x_n^{u'_n(\sigma)}$ are two monomials. We have 
				\begin{equation*}
					\phi^*_l\left( x^{u(\sigma)}x_n^{u_n(\sigma)}\right)= x^{u(\sigma)} x_n^{u_n(\sigma) +\sum_{i=1} ^{n-1}u_i(\sigma)l^i},
				\end{equation*}
				and also
				\begin{equation*}
					\phi^*_l\left( x^{u'(\sigma)}x_n^{u'_n(\sigma)}\right) = x^{u'(\sigma)} x_n^{u_n(\sigma) +\sum_{i=1} ^{n-1}u'_i(\sigma)l^i}.
				\end{equation*}
				Consider the following equality
				\begin{equation*}
					u_n(\sigma) + \sum_{i=1} ^{n-1} u_i(\sigma)l^i = u_n(\sigma) + \sum_{i=1} ^{n-1} u'_i(\sigma)l^i.
				\end{equation*}
				This gives
				\begin{equation*}
					\sum_{i=1} ^{n-1}\left( u_i(\sigma)- u'_i(\sigma)\right)l^i=0.
				\end{equation*}
				Regarding $\sum_{i=1} ^{n-1}\left( u_i(\sigma)- u'_i(\sigma)\right)l^i$ as a polynomial in $\mathbb{Z}[\sigma][l]$, it has finitely many roots in $\mathbb{N}$. If we choose $l$ to be a natural number greater than all of them, then the image of these two monomials under $\phi^*_l$ would be two monomials with different $\sigma$-powers of $x_n$.\\
				Hence, for a large enough natural number $l$, $\phi^*_l$ maps $f$ to a Laurent difference polynomial with different $\sigma$-powers of $x_n$ in its different monomials.
			\end{proof}
			\begin{lem}\label{rootofg}
				Let $f$ be a Laurent difference polynomial in $K_{\sigma}[x_1^{\pm1},\dots,x_n^{\pm1}]$. Suppose
				$\underline{w}=(w_1,\dots,w_n) \in \Gamma^n$, and $\bar{\alpha}=(\bar{\alpha}_1,\dots,\bar{\alpha}_n) \in (\mathbf{k}^*)^n$ are given.\\
				Assume for a large enough natural number $l$, $\phi^*_l$ is defined as in Lemma \ref{phistar}. If  $y'=(y'_1,\dots,y'_n)$ is a root of $g=\phi^*_l(f)$ with the following properties:
				\begin{itemize}
					\item $\forall i,\,\,\,\,\, 1 \leq i\leq n-1:\,\, v(y'_i)=w_i-l^{i}w_n:=w^{\prime}_i$ and $\overline{t^{-w_i+l^{i}w_n}y'_i}=\bar{\alpha}_i \bar{\alpha}_n^{-l^{i}}$ ,
					\item $v(y'_n)=w_n$ and $\overline{t^{-w_n}y'_n}=\bar{\alpha}_n$.
				\end{itemize}
				Then there exists a root $y=(y_1,\dots,y_n)$ for $f$ satisfying the following conditions:
				\begin{itemize}
					\item $v(y)=\underline{w}$,
					\item $\forall i,\,\, 1 \leq i \leq n:\,\, \overline{t^{-w_i}y_{i}}=\bar{\alpha}_i$.
				\end{itemize} 
			\end{lem}
			\begin{proof}	
				Let $y'=(y'_1,\dots,y'_n)$ be a root of $\phi^*_l(f)$ such that $y'_1,\dots,y'_n \in K^*$, and it satisfies the following conditions:
				\begin{itemize}
					\item $\forall i, \,\, 1 \leq i\leq n-1:\,\, v(y'_i)=w'_i=w_i-l^{i}w_n$ and $\overline{t^{-w_i+l^{i}w_n}y'_i}=\bar{\alpha}_i \bar{\alpha}_n^{-l^{i}}$,
					\item $v(y'_n)=w_n$ and $\overline{t^{-w_n}y'_n}=\bar{\alpha}_n$.
				\end{itemize}
				Define $y=(y_1,\dots,y_n)$ as follows: 
				\begin{equation*}
					\forall i, \,\,\, 1 \leq i \leq n-1\,\,\, y_i=y' _i {y'_{n}}^{l^{i}}\,\,\,\textrm{and}\,\,\,\, y_n = y'_n.
				\end{equation*}
				Since $\phi^*_l(f)(y')=0$, from the definition of $\phi^*_l$, we have
				\begin{equation*}
					0=\phi^*_l(f)(y'_1, \dots, y'_n)=f\left(y'_1{y'_{n}}^{l}, \dots, y'_{n-1}{y'_{n}}^{l^{n-1}}, y'_n\right)=f(y),
				\end{equation*}

			which means $y$ is a root of $f$.
			For the root $y$, we have 
			
			\begin{equation}\label{a}
				\begin{aligned}
					v(y)&=(v(y_1) , \dots , v(y_{n-1}),v(y_n))\\ 
					&=(v(y'_1) + l^1v(y'_n), \dots , v(y'_{n-1}) + l^{n-1}v(y'_n),v(y'_n))=\underline{w}.
				\end{aligned}
			\end{equation}
			We assumed for the root $y'$ of $\phi^*_l(f)$ that 
			\begin{equation*}
				\overline{t^{-w_i + l^iw_n}y'_i} = \bar{\alpha}_i\bar{\alpha}_n^{-l_i}\,\,  \textrm{for all}\,\,\, i, \,\,\, 1\leq i\leq n-1.
			\end{equation*}
			So for all $i$ with $1\leq i \leq n-1$, by an easy computation, we obtain
			\begin{equation}\label{b}
				\begin{aligned}
					\bar{\alpha}_i 
					= \overline{t^{-w_i + l^iw_n - w_nl^i}y'_i{y'_n}^{l^i}}
					=\overline{t^{-w_i}y_i}.
				\end{aligned}
			\end{equation}
			By \eqref{a} and \eqref{b}, the root $y$ of $f$ has the intended properties.
		\end{proof}
		\begin{lem} \label{initialroot}
			Let $f$ be a Laurent difference polynomial in $K_{\sigma}[x_1^{\pm1},\dots,x_n^{\pm1}]$. Suppose
			$\underline{w}=(w_1,\dots,w_n) \in \Gamma^n$, and $\bar{\alpha}=(\bar{\alpha}_1,\dots,\bar{\alpha}_n) \in (\mathbf{k}^*)^n$ are given. Assume that $\mathopin_{\underline{w}}f$ is not a monomial and we have $\mathopin_{\underline{w}}f(\bar{\alpha})=0$ .\\
			If for a large enough natural number $l$, we define $\phi^*_l$ as in Lemma \ref{phistar}, then for $g=\phi^*_l(f)$, we have $\mathopin_{\underline{w}'}(g)$ is not a monomial, where $\underline{w}'=(w_1-l^1w_n,\dots,w_{n-1}-l^{(n-1)}w_n,w_n)$, and $\bar{\alpha}'=(\bar{\alpha}_1\bar{\alpha}_n^{-l^1},\dots,\bar{\alpha}_{n-1}\bar{\alpha}_n^{-l^{(n-1)}},\bar{\alpha}_n)$ is a root of $\mathopin_{\underline{w}'}(g)$.
		\end{lem}
		\begin{proof}
			We start by proving that $\mathoptrop(f)(\underline{w})=\mathoptrop\left(\phi^*_l(f) \right)(\underline{w}')$. To do so, we look at the tropicalization of each monomial in each of them.
			The tropicalization of the monomial $c_{\underline{u}(\sigma)}x^{u(\sigma)}x_n^{u_n(\sigma)}$ of $f$ is
			\begin{equation*}
				\mathoptrop(c_{\underline{u}(\sigma)} x^{u(\sigma)}x_n^{u_n(\sigma)})(\underline{w})=v(c_{\underline{u}(\sigma)})+\underline{u}(\rho)\cdot \underline{w},
			\end{equation*}
			and the tropicalization of the corresponding monomial in $\phi^*_l(f)$ at the point $\underline{w}'$ is
			\begin{align*}
				&\mathoptrop \left( \phi^*_l \left( c_{\underline{u}(\sigma)}x^{u(\sigma)}x_n^{u_n(\sigma)} \right) \right)(\underline{w}')=\mathoptrop\left(c_{\underline{u}(\sigma)}x^{u(\sigma)}x_n^{u_n(\sigma)+\sum_{i=1}^{n-1}u_i(\sigma)l^{i}} \right)(\underline{w}')\\
				&=v(c_{\underline{u}(\sigma)}) + \underline{u}(\rho)\cdot \underline{w}'
				+ \sum_{i=1}^{n-1}l^iu_i(\rho)\cdot w'_n\\   
				&=v(c_{\underline{u}(\sigma)}) + (u_1(\rho) , \dots ,u_{n-1}(\rho),u_n(\rho) )(w_1 - l^1w_n, \dots , w_{n-1}- l^{(n-1)}w_n,w_n )\\
				&+\sum_{i=1}^{n-1}l^iu_i(\rho)\cdot w_n\\
				&=v(c_{\underline{u}(\sigma)})+ \sum_{i=1}^{n-1}u_i(\rho)w_i - l^iu_i(\rho )w_n+u_n(\rho)w_n+\sum_{i=1}^{n-1}l^iu_i(\rho)\cdot w_n\\
				&=v(c_{\underline{u}(\sigma)})+ \underline{u}(\rho)\cdot \underline{w} = \mathoptrop(c_{\underline{u}(\sigma)}x^{u(\sigma)}x_n^{u_n(\sigma)})(\underline{w}).
			\end{align*}
			Hence, we have
			\begin{equation*}
				\mathoptrop(f)(\underline{w})=\mathoptrop\left(\phi^*_l(f) \right)(\underline{w}').
			\end{equation*}	   
			Suppose $\mathopin_{\underline{w}}(f)$ is not a monomial. Let
			\begin{equation*}
				\overline{t^{-v(c_{\underline{u}}(\sigma))}\cdot c_{\underline{u}(\sigma)}}x^{u(\sigma)}x_n^{u_n(\sigma)}   
			\end{equation*}
			be one of its monomials. From the definition of $\mathopin_{\underline{w}}(f)$, it follows that
			\begin{equation*}
				\mathoptrop\left(c_{\underline{u}(\sigma)}x^{u(\sigma)}x_n^{u_n(\sigma)} \right)(\underline{w}) 
			\end{equation*}
			attains the minimum. Since $\mathoptrop(f)(\underline{w})=\mathoptrop\left(\phi^*_l(f) \right)(\underline{w}')$, we conclude that,
			\begin{equation*}
				\mathoptrop\left( \phi^*_l(c_{\underline{u}(\sigma)}x^{u(\sigma)}x_n^{u_n(\sigma)}) \right)(\underline{w}')  
			\end{equation*}
			attains the minimum in $\mathoptrop\left(\phi^*_l(f) \right)(\underline{w}') $. Therefore, the corresponding monomial appears in $\mathopin_{\underline{w}'}\left(\phi^*_l (f)\right)$. Thus, if $\mathopin_{\underline{w}}(f)$ has more than one monomial, so does $\mathopin_{\underline{w}'}\left(\phi^*_l (f)\right)$.
			More precisely, this corresponding monomial in
			$\mathopin_{\underline{w}'}\left(\phi^*_l(f)\right)$ is
			\begin{equation*}
				\overline{t^{-v(c_{\underline{u}(\sigma)})}c_{\underline{u}(\sigma)}}\phi^*_l(x^{u(\sigma)}x_n^{u_n(\sigma)})=\overline{t^{-v(c_{\underline{u}(\sigma)})}c_{\underline{u}(\sigma)}}x^{u(\sigma)}x_n^{u_n(\sigma)+\sum_{i=1}^{n-1}u_i(\sigma)l^{i}}.    
			\end{equation*}
			Assume $\bar{\alpha}$ is a root of $\mathopin_{\underline{w}}(f)$. By straight forward computation, we obtain
			\begin{align*}
				&\left( \overline{t^{-v(c_{\underline{u}(\sigma)})}c_{\underline{u}(\sigma)}}x^{u(\sigma)}x_n^{u_n(\sigma)+\sum_{i=1}^{n-1}u_i(\sigma)l^{i}}\right)(\bar{\alpha}')\\
				&=\left( \overline{t^{-v(c_{\underline{u}(\sigma)})}c_{\underline{u}(\sigma)}} x^{u(\sigma)}x_n^{u_n(\sigma)}\right) (\bar{\alpha}).
			\end{align*}
			Hence, we have 
			$\mathopin_{\underline{w}'}(\phi ^{\ast} _l(f))(\bar{\alpha}' ) = \mathopin_{\underline{w}}(f)(\bar{\alpha})=0 $ which means $\bar{\alpha}'$ is a root of 
			$\mathopin _{\underline{w}'}(\phi ^{\ast} _l(f))$.
		\end{proof}
		
		By combining Proposition \ref{3.1.5} with the three preceding lemmas, we derive the main result of this subsection, which is presented in the following proposition. 
		
		\begin{prop}\label{main}
			Let $f \in K_{\sigma}[x_1^{\pm1},\dots,x_n^{\pm1}]$ be a Laurent difference polynomial, and $\underline{w}=(w_1, \dots, w_n) \in \Gamma^n$ such that $\mathopin_{\underline{w}}(f)$ is not a monomial. Suppose $\bar{\alpha}$ is a root of $\mathopin_{\underline{w}}(f)$ in $(\mathbf{k^*})^n$. Then there exists an element $y$ in $(K^*)^n$ which is a root of $f$, and satisfies the following conditions:
			\begin{itemize}
				\item $v(y)=\underline{w}$,
				\item $\forall i,\,\,\,1\leq i \leq n:\,\, \overline{t^{-w_i}\cdot y_i}=\bar{\alpha}_i$.
			\end{itemize}
			
		\end{prop}
		\begin{proof}
			By Lemma \ref{phistar},
			there exists a Laurent difference polynomial $g$ corresponding to $f$ with different $\sigma$-powers of $x_n$in its different monomials. Lemma \ref{initialroot} guarantees that the assumptions of Proposition \ref{3.1.5} hold for the Laurent difference polynomial $g$, $\underline{w}'$ and $\Bar{\alpha}'$ (where $\underline{w}'$ and $\Bar{\alpha}'$ are defined as in Lemma \ref{initialroot}). Hence, from Proposition \ref{3.1.5}, we find a root $y'$ for $g$ which satisfies the following conditions:
			\begin{itemize}
				\item $v(y')=\underline{w}'$,
				\item $\forall i,\,\,\,1\leq i \leq n:\,\, \overline{t^{-w'_i}.y'_i}=\bar{\alpha}'_i$.
			\end{itemize} 
			Finally, Lemma \ref{rootofg} implies that 
			corresponding to $y'$, there exist a root $y$ for $f$ which satisfies the desired conditions.
			
		\end{proof}
		\subsection{The Difference Kapranov Theorem}
		In this subsection by using Proposition \ref{main}, we prove the difference version of Kapranov's theorem which is the main result of this paper.
		\begin{thm} \label{diffkap}{(The Difference Kapranov Theorem)}
			
			Let $K$ be a multiplicative valued difference field of characteristic zero that is spherically complete.
			Assume its difference residue field is an ACFA of characteristic zero and the scaling exponent, $\rho$ is transcendental.\\
			Let the difference value group $\Gamma$ of $K$ be a subgroup of $\mathbb{R}$ that is a $\mathbb{Q}(\rho)$-module.\\
			Suppose $f \in K_{\sigma}\left[x_1^{\pm1},\dots,x_n^{\pm1}\right]$ is a Laurent difference polynomial.The following sets coincide:
			\begin{enumerate}
				\item $\mathoptrop(V(f)) \subseteq \mathbb{R}^n$ which is the difference tropical hypersurface associated to $f$;
				\item the set of all the points $w \in \mathbb{R}^n$ for which the initial form $\mathopin_w(f)$ is not a monomial;
				\item the closure of the set $A=\left\lbrace (v(y_1),\dots,v(y_n)):(y_1,\dots,y_n)\in V(f) \right\rbrace $ in $\mathbb{R}^n$.
			\end{enumerate} 
			
		\end{thm}
		\begin{proof}
			$(1)=(2)$: As it is defined in Definition \ref{trophyp}, $\mathoptrop(V(f))$ is the set of all tropical roots of $f$. Therefore by Lemma \ref{roots}, the sets in (1) and (2) are equal.\\
			$(3)\subseteq (1)$: From Proposition \ref{combi}, we know that $\mathoptrop(V(f))$ is the support of a polyhedral complex, or equivalently it is the union of some polyhedra. As $f$ has finitely many monomials, $\mathoptrop(V(f))$ is the union of finitely many polyhedra, each of which is closed so $\mathoptrop(V(f))$ is closed.\\
			Let $(v(y_1),\dots,v(y_n)) \in A$. From the definition of $A$, $y=(y_1,\dots,y_n)\in (K^*)^{n}$ is a root of $f$. This means for any monomial $c_{u(\sigma)}x^{u(\sigma)}$ of $f$, where $c_{u(\sigma)} \neq 0$, we have $v(c_{u(\sigma)}y^{u(\sigma)}) < v(f(y))=v(0)= \infty$; in fact
			\begin{equation*}
				v\left(\underset{u(\sigma) \in (\mathbb{Z}[\sigma])^n}{\sum}c_{u(\sigma)}y^{u(\sigma)}  \right) \neq \underset{\substack{u(\sigma) \in (\mathbb{Z}[\sigma])^n \\ c_{u(\sigma)}\neq 0}}{\min}\left\lbrace v(c_{u(\sigma)})+ u(\rho)\cdot v(y) \right\rbrace. 
			\end{equation*}
			So there exist two indices $u(\sigma)$ and $u'(\sigma)$ for which we have
			\begin{equation*}
				v(c_{u(\sigma)})+u(\rho)\cdot v(y)=v(c_{u'(\sigma)})+u'(\rho)\cdot v(y),    
			\end{equation*}
			and they acheive the minimum, or in other words $\mathoptrop(f)$ attains its minimum in $v(y)$ at least twice. This means $v(y)\in \mathoptrop(V(f))$, so $\left\lbrace v(y): y \in V(f) \right\rbrace \subseteq \mathoptrop(V(f))$. As $\mathoptrop(V(f))$ is closed, we have 
			$\overline{\left\lbrace v(y): y \in V(f) \right\rbrace } \subseteq \mathoptrop(V(f))$. It follows that the set in (3) is contained in the set in (1).\\
			$(1)\subseteq (3)$: Let $w \in \mathoptrop(V(f))\cap \Gamma^n$. Since $\mathoptrop(f)$ attains its minimum at $w$ at least twice, we deduce from the equality of the sets in (1) and (2) that $\mathopin_w(f)$ is not a monomial. Besides, $\mathopin_w(f)$ is a Laurent difference polynomial with coefficients in an ACFA field. Therefore, by Theorem \ref{alpharoot} $\mathopin_w(f)$ has a root $\Bar{\alpha} \in (\mathbf{k^*})^n$. So by Proposition \ref{main}, there exists $y \in (K^{*})^n$ such that $f(y)=0$, or equivalently $y \in V(f)$ and $v(y)=w$. This means that $ \mathoptrop(V(f))\cap \Gamma^n \subseteq A$. Moreover, $\mathoptrop(V(f))\cap \Gamma^n$ is dense in $\mathoptrop(V(f))$. To see this, consider Proposition \ref{combi}. From this result, $\mathoptrop(V(f))$ is the support of a pure $\left(\Gamma,\mathbb{Q}(\rho)\right)-$polyhedral complex of dimension $(n-1)$. Suppose $P$ is a facet of this polyhedral complex. So it is a $\left(\Gamma,\mathbb{Q}(\rho)\right)-$polyhedron of dimension $(n-1)$. Define the projection map $\pi$ on $P$, which takes each point to its first $(n-1)$ coordinates. This map is also bijective. Since  the interior of $\pi(p)$ in $\mathbb{R}^{n-1}$ is nonempty, we have $\overline{\pi(p)^{\circ}}=\pi(p)$. From this, it is not difficult to show that $\pi(p)\cap \Gamma^{n-1}$ is dense in $\pi(p)$. Using the bijection, $P\cap \Gamma^n$ is also dense in $P$. Consequently $\mathoptrop(V(f))\cap \Gamma^n$ is dense in $\mathoptrop(V(f))$. Thus, we have 
			\begin{equation*}
				\mathoptrop(V(f))=\overline{\mathoptrop(V(f))\cap \Gamma^n}\subseteq \overline{A}.
			\end{equation*}
			Hence, $\mathoptrop(V(f))$ is included in the set in (3) and this implies the equality of the sets in (1) and (3).

		\end{proof}
	\end{fullversion}

	\eject \pdfpagewidth=14in \pdfpageheight=7in
	\clearpage

\end{document}